\documentclass[12pt]{amsart}

\usepackage{epsfig}
\usepackage{xspace}
\usepackage{layout}
\usepackage{latexsym}
\usepackage{amsmath}
\usepackage{amsthm}
\usepackage{amssymb}
\usepackage{amsfonts}
\usepackage{amsxtra}     

\usepackage{verbatim}
\usepackage{longtable}

\usepackage[all, dvips]{xy}
\usepackage{amsmath}
\usepackage{amstext}
\usepackage{amsfonts}
\usepackage{amscd}
\usepackage{mathrsfs}
\usepackage{latexsym}
\usepackage{amssymb}

\newcommand{\id}{\operatorname{id}}

\newcommand{\ie}{\textit{i.e., }}

\newcommand{\bbB}{\mathbb{B}}

\newcommand{\bbD}{\mathbb{D}}
\newcommand{\bbE}{\mathbb{E}}

\newcommand{\bbN}{\mathbb{N}}
\newcommand{\bbP}{\mathbb{P}}
\newcommand{\bbQ}{\mathbb{Q}}
\newcommand{\bbR}{\mathbb{R}}
\newcommand{\bbS}{\mathbb{S}}
\newcommand{\bbT}{\mathbb{T}}

\newcommand{\bfB}{\mathbf{B}}
\newcommand{\bfC}{\mathbf{C}}
\newcommand{\bfD}{\mathbf{D}}

\newcommand{\bfJ}{\mathbf{J}}
\newcommand{\bfL}{\mathbf{L}}

\newcommand{\bfP}{\mathbf{P}}
\newcommand{\bfQ}{\mathbf{Q}}
\newcommand{\bfR}{\mathbf{R}}
\newcommand{\bfS}{\mathbf{S}}
\newcommand{\bfT}{\mathbf{T}}

\newcommand{\bfV}{\mathbf{V}}

\newcommand{\op}{\text{\rm op}}

\newcommand{\Obj}{\text{\rm Obj}\;}

\newcommand{\Ho}{\text{\rm Ho}\;}

\newcommand{\colim}{\text{\rm colim}\;}
\newcommand{\im}{\text{\it im}\;}

\newcommand{\dom}{\text{\it dom}\;}

\newcommand{\Cat}{\text{\sf Cat}}

\newcommand{\Sd}{\text{\rm Sd}}
\newcommand{\Ex}{\text{\rm Ex}}

\newcommand{\Cen}{\text{\rm {\bf Center}}}
\newcommand{\Out}{\text{\rm {\bf Outer}}}
\newcommand{\Comp}{\text{\rm {\bf Comp}}}

\newcommand{\Sym}{\text{\rm Sym}}

\newcommand{\co}{\colon\thinspace}

\newcommand{\tb}[1]{\phantom{\sum^\Sigma_\Sigma} #1 \phantom{\sum^\Sigma_\Sigma}}
\newcommand{\lr}[1]{\hspace{.5mm}#1\hspace{.5mm}}
\newcommand{\halo}[1]{\phantom{e}\underset{\phantom{e}}{\overset{\phantom{e}}{#1}}\phantom{e}}

\begin{document}

\newtheorem{thm}{Theorem}[section]
\newtheorem{conj}[thm]{Conjecture}
\newtheorem{lem}[thm]{Lemma}
\newtheorem{cor}[thm]{Corollary}
\newtheorem{prop}[thm]{Proposition}
\newtheorem{rem}[thm]{Remark}

\theoremstyle{definition}
\newtheorem{defn}[thm]{Definition}
\newtheorem{examp}[thm]{Example}
\newtheorem{notation}[thm]{Notation}
\newtheorem{rmk}[thm]{Remark}

\theoremstyle{remark}

\makeatletter
\renewcommand{\maketag@@@}[1]{\hbox{\m@th\normalsize\normalfont#1}}%
\makeatother
\renewcommand{\labelenumi}{(\roman{enumi})}
\renewcommand{\labelenumii}{(\alph{enumii})}
\renewcommand{\theenumi}{(\roman{enumi})}

\def\square{\hfill ${\vcenter{\vbox{\hrule height.4pt \hbox{\vrule width.4pt
height7pt \kern7pt \vrule width.4pt} \hrule height.4pt}}}$}

\newenvironment{pf}{{\it Proof:}\quad}{\square \vskip 12pt}
\date{\today}

\title[Thomason Structure on {\bf nFoldCat}]{A Thomason Model Structure on the Category of Small $n$-fold Categories}
\author{Thomas M. Fiore and Simona Paoli}
\address{Thomas M. Fiore \\
Department of Mathematics and Statistics \\
University of Michigan-Dearborn \\
4901 Evergreen Road \\
Dearborn, MI 48128 \\
USA} \email{tmfiore@umd.umich.edu}
\address{Simona Paoli \\
Department of Mathematics \\
Penn State Altoona \\
3000 Ivyside Park \\
Altoona, PA 16601-3760 \\
USA} \email{sup24@psu.edu}


\begin{abstract}
We construct a cofibrantly generated Quillen model structure on the
category of small $n$-fold categories and prove that it is Quillen
equivalent to the standard model structure on the category of
simplicial sets. An $n$-fold functor is a weak equivalence if and
only if the diagonal of its $n$-fold nerve is a weak equivalence of
simplicial sets. This is an $n$-fold analogue to Thomason's Quillen
model structure on $\mathbf{Cat}$. We introduce an $n$-fold
Grothendieck construction for multisimplicial sets, and prove that
it is a homotopy inverse to the $n$-fold nerve. As a consequence, we
completely prove that the unit and counit of the adjunction between
simplicial sets and $n$-fold categories are natural weak
equivalences.
\end{abstract}

\maketitle

\tableofcontents
\section{Introduction}
An $n$-fold category is a higher and wider categorical structure
obtained by $n$ applications of the internal category construction.
In this paper we study the homotopy theory of $n$-fold categories.
Our main result is Theorem \ref{maintheoremsummary}. Namely, we have
constructed a cofibrantly generated model structure on the category
of small $n$-fold categories in which an $n$-fold functor is a weak
equivalence if and only if its nerve is a diagonal weak equivalence.
This model structure is Quillen equivalent to the usual model
structure on the category of simplicial sets, and hence also
topological spaces. Our main tools are model category theory, the
$n$-fold nerve, and an $n$-fold Grothendieck construction for
multisimplicial sets. Notions of nerve and versions of the
Grothendieck construction are very prominent in homotopy theory and
higher category theory, as we now explain. The Thomason model
structure on $\mathbf{Cat}$ is also often present, at least
implicitly.

The Grothendieck nerve of a category and the Grothendieck
construction for functors are fundamental tools in homotopy theory.
Theorems A and B of Quillen \cite{quillenI}, and Thomason's theorem
\cite{thomasonhocolimit} on Grothendieck constructions as models for
certain homotopy colimits, are still regularly applied decades after
their creation. Functors with nerves that are weak equivalences of
simplicial sets feature prominently in these theorems. Such functors
form the weak equivalences of Thomason's model structure on {\bf
Cat} \cite{thomasonCat}, which is Quillen equivalent to {\bf SSet}.
Earlier, Illusie \cite{illusieII} proved that the nerve and the
Grothendieck construction are homotopy inverses. Although the nerve
and the Grothendieck construction are not adjoints\footnote{In fact, the Grothendieck
construction is not even homotopy equivalent to $c$, the left adjoint to the nerve, as follows. For any simplicial set $X$, let $\Delta/X$ denote the Grothendieck construction on $X$. Then $N(\Delta/\partial\Delta[3])$ is homotopy equivalent to $\partial\Delta[3]$ by Illusie's result. On the other hand, $Nc\partial\Delta[3]=Nc\Delta[3]=\Delta[3]$, since $cX$ only depends on 0-,1-, and 2-simplices. Clearly, $\partial\Delta[3]$ and $\Delta[3]$ are not homotopy equivalent, so the Grothendieck construction
is not naturally homotopy equivalent to $c$. }, the equivalence
of homotopy categories can be realized by adjoint functors
\cite{fritschlatch1}, \cite{fritschlatch2}, \cite{thomasonCat}.
Related results on homotopy inverses are found in
\cite{latchuniqueness}, \cite{lee}, and \cite{waldhausen}. More
recently, Cisinski \cite{cisinskiasterisque} has proved two
conjectures of Grothendieck concerning this circle of ideas (see
also \cite{jardinesummary}).

On the other hand, notions of nerve play an important role in
various definitions of $n$-category \cite{leinstersurvey}, namely
the definitions of Simpson \cite{simpsonmodel}, Street
\cite{street}, and Tamsamani \cite{tamsamani}, as well as in the
theory of quasi-categories developed by Joyal \cite{joyalNotes},
\cite{joyalVolume1}, \cite{joyalVolume2},  and also Lurie
\cite{lurieStableInfinity}, \cite{lurieHigherToposTheory}. For
notions of nerve for bicategories, see for example work of Duskin
and Lack-Paoli \cite{duskinI}, \cite{duskinII}, \cite{lackpaoli},
and for left adjoints to singular functors in general also
\cite{gabrielzisman} and \cite{kellyenriched}. Fully faithful
cellular nerves have been developed for higher categories in
\cite{bergercellular}, together with characterizations of their
essential images. Nerve theorems can be established in a very
general context, as proved by Leinster and Weber in
\cite{leinsternerves} and \cite{weber}, and discussed in
\cite{ncategorycafenerves}. As an example, Kock proves in
\cite{kocktrees} a nerve theorem for polynomial endofunctors in
terms of trees.

Model category techniques are only becoming more important in the
theory of {\it higher} categories. They have been used to prove
that, in a precise sense, simplicial categories, Segal categories,
complete Segal spaces, and quasi-categories are all equivalent
models for $(\infty,1)$-categories \cite{bergnersurvey},
\cite{bergnerthreemodels}, \cite{bergnersimplicialcategories},
\cite{joyaltierneyquasisegal}, \cite{rezkhomotopytheory}, and
\cite{toenaxiomatization}. In other directions, although the
cellular nerve of \cite{bergercellular} does not transfer a model
structure from cellular sets to $\omega$-categories, it is proved in
\cite{bergercellular} that the homotopy category of cellular sets is
equivalent to the homotopy category of $\omega$-categories. For
this, a Quillen equivalence between {\it cellular spaces} and {\it
simplicial $\omega$-categories} is constructed. There is also the
work of Simpson and Pellisier \cite{pellissier},
\cite{simpsonmodel}, and \cite{simpsonhigher}, developing model
structures on $n$-categories for the purpose of $n$-stacks, and also
a model structure for $(\infty, n)$-categories.

In low dimensions several model structures have already been
investigated. On {\bf Cat}, there is the categorical structure of
Joyal-Tierney \cite{joyaltierney}, \cite{rezkcat}, as well as the
topological structure of Thomason \cite{thomasonCat},
\cite{cisinskiThomasonFix}. A model structure on pro-objects in {\bf
Cat} appeared in \cite{golasinski}, \cite{golasinskiprotranslation},
\cite{golasinskipro}. The articles
\cite{heggiehomotopycofibrations}, \cite{heggietensorproduct},
\cite{heggiehomotopycolimits}  and are closely related to the
Thomason structure and the Thomason homotopy colimit theorem. More
recently, the Thomason structure on {\bf Cat} was proved in Theorem
5.2.12 of \cite{cisinskiasterisque} in the context of Grothendieck
test categories and fundamental localizers. The homotopy categories
of spaces and categories are proved equivalent in \cite{hoyo}
without using model categories.

On {\bf 2-Cat} there is the categorical structure of \cite{lack2Cat}
and \cite{lackBiCat}, as well as the Thomason structure of
\cite{worytkiewicz2Cat}. Model structures on {\bf 2FoldCat} have
been studied in \cite{fiorepaolipronk1} in great detail. The
homotopy theory of 2-fold categories is very rich, since there are
numerous ways to view 2-fold categories: as internal categories in
{\bf Cat}, as certain simplicial objects in {\bf Cat}, or as
algebras over a 2-monad. In \cite{fiorepaolipronk1}, a model
structure is associated to each point of view, and these model
structures are compared.

However, there is another way to view 2-fold categories not treated
in \cite{fiorepaolipronk1}, namely as certain bisimplicial sets.
There is a natural notion of fully faithful double nerve, which
associates to a 2-fold category a bisimplicial set. An obvious
question is: does there exist a Thomason-like model structure on
{\bf 2FoldCat} that is Quillen equivalent to some model structure on
bisimplicial sets via the double nerve? Unfortunately, the left
adjoint to double nerve is homotopically poorly behaved as it
extends the left adjoint $c$ to ordinary nerve, which is itself
poorly behaved. So any attempt at a model structure must address
this issue.

Fritsch, Latch, and Thomason \cite{fritschlatch1},
\cite{fritschlatch2}, \cite{thomasonCat} noticed that the composite
of $c$ with second barycentric subdivision $\Sd^2$ is much better
behaved than $c$ alone. In fact, Thomason used the adjunction
$c\Sd^2 \dashv \Ex^2N$ to construct his model structure on {\bf
Cat}. This adjunction is a Quillen equivalence, as the right adjoint
preserves weak equivalences and fibrations by definition, and the
unit and counit are natural weak equivalences.

Following this lead, we move to simplicial sets via $\delta^\ast$
(restriction to the diagonal) in order to correct the homotopy type
of double categorification using $\Sd^2$. Moreover, our method of
proof works for $n$-fold categories as well, so we shift our focus
from $2$-fold categories to general $n$-fold categories. In this
paper, we construct a cofibrantly generated model structure on {\bf
nFoldCat} using the fully faithful $n$-fold nerve, via the
adjunction below,
\begin{equation} \label{mainadjunctionintro}
\xymatrix@C=4pc{\mathbf{SSet} \ar@{}[r]|{\perp}
\ar@/^1pc/[r]^-{\Sd^2} &  \ar@/^1pc/[l]^-{\Ex^2} \mathbf{SSet}
\ar@{}[r]|{\perp} \ar@/^1pc/[r]^-{\delta_!} &
\ar@/^1pc/[l]^-{\delta^\ast} \mathbf{SSet^n} \ar@{}[r]|{\perp}
\ar@/^1pc/[r]^-{c^n} & \ar@/^1pc/[l]^-{N^n} \mathbf{nFoldCat}}
\end{equation}
and prove that the unit and counit are weak equivalences. Our method
is to apply Kan's Lemma on Transfer of Structure. First we prove
Thomason's classical theorem in Theorem \ref{CatCase}, and then use
this proof as a basis for the general $n$-fold case in Theorem
\ref{MainModelStructure}. We also introduce an $n$-fold Grothendieck
construction in Definition \ref{nfoldGrothendieck}, prove that it is
homotopy inverse to the $n$-fold nerve in Theorems \ref{rhowe} and
\ref{lambdawe}, and conclude in Proposition \ref{unitcounitwe} that
the unit and counit of the adjunction (\ref{mainadjunctionintro})
are natural weak equivalences. The articles \cite{fritschlatch1} and
\cite{fritschlatch2} proved in a different way that the unit and
counit of the classical Thomason adjunction
$\mathbf{SSet}\dashv\mathbf{Cat}$ are natural weak equivalences.

Recent interest in $n$-fold categories has focused on the $n=2$
case. In many cases, this interest stems from the fact that 2-fold
categories provide a good context for incorporating two types of
morphisms, and this is useful for applications. For example, between
rings there are ring homomorphisms and bimodules, between
topological spaces there are continuous maps and parametrized
spectra as in \cite{maysigurdsson}, between manifolds there are
smooth maps and cobordisms, and so on. In this direction, see for
example \cite{grandisdouble1}, \cite{fiore1}, \cite{fiore2},
\cite{mortondouble}, \cite{shulmanonquillenfunctors},
\cite{shulmanframed}. Classical work on 2-fold categories,
originally introduced by Ehresmann as {\it double categories},
includes \cite{ehresmannone}, \cite{ehresmanntwo},
\cite{ehresmannthree}, \cite{ehresmannfour}, \cite{ehresmann2},
\cite{ehresmann}. The theory of double categories is now
flourishing, with many contributions by Brown-Mosa, Grandis-Par\'e,
Dawson-Par\'e-Pronk, Dawson-Par\'e, Fiore-Paoli-Pronk, Shulman, and
many others. To mention only a few examples, we have
\cite{brownmosa99}, \cite{grandisdouble1}, \cite{grandisdouble2},
\cite{grandisdouble3}, \cite{grandisdouble4},
\cite{dawsonparepronkpaths}, \cite{dawsonparefreedouble},
\cite{fiorepaolipronk1}, \cite{shulmanonquillenfunctors}, and
\cite{shulmanframed}.

There has also been interest in general $n$-fold categories from
various points of view. Connected homotopy $(n+1)$-types are
modelled by $n$-fold categories internal to the category of groups
in \cite{lodayfinitelymany}, as summarized in the survey paper
\cite{paoliinternalstructures}. Edge symmetric $n$-fold categories
have been studied by Brown, Higgins, and others for many years now,
for example \cite{brownhigginsgroupoidscrossedcomplexes},
\cite{brownhigginsgroupoidscubicalTcomplexes},
\cite{brownhigginscubes}, and \cite{brownhigginstensor}. There are
also the more recent {\it symmetric weak cubical categories} of
\cite{grandiscospans3} and  \cite{grandiscospans1}. The homotopy
theory of cubical sets has been studied in \cite{jardineCubical}.

The present article is the first to consider a Thomason structure on
the category of $n$-fold categories. Our paper is organized as
follows. Section \ref{nfoldcategories} recalls $n$-fold categories,
introduces the $n$-fold nerve $N^n$ and its left adjoint $n$-fold
categorification $c^n$, and describes how $c^n$ interacts with
$\delta_!$, the left adjoint to precomposition with the diagonal. In
Section \ref{barycentric} we recall barycentric subdivision,
including explicit descriptions of $\Sd^2\Lambda^k[m]$,
$\Sd^2\partial\Delta[m]$, and $\Sd^2\Delta[m]$. More importantly, we
present a decomposition of the poset $\bfP\Sd \Delta[m]$ into the
union of three posets $\Comp$, $\Cen$, and $\Out$ in Proposition
\ref{upcloseddecomposition}, as pictured in Figure
\ref{subdivisionfigure} for $m=2$ and $k=1$. Though Section
\ref{barycentric} may appear technical, the statements become clear
after a brief look at the example in Figure \ref{subdivisionfigure}.
This section is the basis for the verification of the pushout axiom
\ref{KanCorollaryiv} of Corollary \ref{KanCorollary}, completed in
the proofs of Theorems \ref{CatCase} and \ref{MainModelStructure}.

Sections \ref{retractionsection} and \ref{pushoutsection} make
further preparations for the verification of the pushout axiom.
Proposition \ref{deformationretract} gives a deformation retraction
of $|N(\Comp \cup \Cen)|$ to part of its boundary, see Figure
\ref{subdivisionfigure}. This deformation retraction finds
application in equation \eqref{QPinclusion}. The highlights of
Section \ref{pushoutsection} are Proposition
\ref{nervecommuteswithpushout} and Corollary
\ref{nervecommuteswithcolimitdecomposition} on the commutation of
nerve with certain colimits of posets. Proposition
\ref{nervecommuteswithpushout} on commutation of nerve with certain
pushouts finds application in equation \eqref{QPinclusion}. Other
highlights of Section \ref{pushoutsection} are Proposition
\ref{colimitdecomposition}, Proposition
\ref{simplicial_colimitdecomposition}, and Corollary
\ref{cor:specific_colimit_decompositions} on the expression of
certain posets (respectively their nerves) as a colimit of two
ordinals (respectively two standard simplices).  Section
\ref{Thomasonsection} pulls these results together and quickly
proves the classical Thomason theorem.

Section \ref{sectionnfolddecompositions} proves the $n$-fold
versions of the results in Sections \ref{barycentric},
\ref{retractionsection}, and \ref{pushoutsection}. The $n$-fold
version of Proposition \ref{colimitdecomposition} on colimit
decompositions of certain posets is Proposition
\ref{colimitdecompositionnfold}. The $n$-fold version of Corollary
\ref{nervecommuteswithcolimitdecomposition} on the commutation of
nerve with certain colimits of posets is Proposition
\ref{diagonaldecomposition}. The $n$-fold version of the deformation
retraction in Proposition \ref{deformationretract} is Corollary
\ref{nfolddeformationretract}. The $n$-fold version of Proposition
\ref{nervecommuteswithpushout} on commutation of nerve with certain
pushouts is Proposition \ref{nfoldnervecommuteswithpushout}.
Proposition \ref{PushoutDescription} displays a calculation of a
pushout of double categories, and the diagonal of its nerve is
characterized in Proposition \ref{pushoutsimplexdescription}.

Section \ref{Thomasonnfoldsection} pulls together the results of
Section \ref{sectionnfolddecompositions} to prove the Thomason
structure on {\bf nFoldCat} in Theorem \ref{MainModelStructure}. In
the last section of the paper, Section \ref{unitcounitsection}, we
introduce a Grothendieck construction for multisimplicial sets and
prove that it is a homotopy inverse for $n$-fold nerve in Theorems
\ref{rhowe} and \ref{lambdawe}. As a consequence, we have in
Proposition \ref{unitcounitwe} that the unit and counit are weak
equivalences.

We have also included an appendix on the Multisimplicial
Eilenberg-Zilber Lemma.

{\bf Acknowledgments:} Thomas M. Fiore and Simona Paoli thank the
Centre de Recerca Matem\`{a}tica in Bellaterra (Barcelona) for its
generous hospitality, as it provided a fantastic working environment
and numerous inspiring talks. The CRM Research Program on Higher
Categories and Homotopy Theory in 2007-2008 was a great inspiration
to us both.

We are indebted to Myles Tierney for suggesting to use the weak
equivalence $\xymatrix@1{N(\Delta/X) \ar[r] & X}$ and the Weak
Equivalence Extension Theorem \ref{weakequivalenceextension} of
Joyal-Tierney \cite{joyaltierneysimplicial} in our proof that the
unit and counit of (\ref{nfoldcatadjunction}) are weak equivalences.
We also thank Andr\'e Joyal and Myles Tierney for explaining aspects
of Chapter 6 of their book \cite{joyaltierneysimplicial}.

We thank Denis-Charles Cisinski for explaining to us his proof that
the unit and counit are weak equivalences in the Thomason structure
on $\mathbf{Cat}$, as this informed our Section
\ref{unitcounitsection}. We also thank Dorette Pronk for several
conversations related to this project.

We also express our gratitude to an anonymous referee who made many
excellent suggestions.

Thomas M.~Fiore was supported at the University of Chicago by NSF
Grant DMS-0501208. At the Universitat Aut\`{o}noma de Barcelona he
was supported by grant SB2006-0085 of the Spanish Ministerio de
Educaci\'{o}n y Ciencia under the Programa Nacional de ayudas para
la movilidad de profesores de universidad e investigadores
espa$\tilde{\text{n}}$oles y extranjeros. Simona Paoli was supported
by Australian Postdoctoral Fellowship DP0558598 at Macquarie
University. Both authors also thank the Fields Institute for its
financial support, as this project began at the 2007 Thematic
Program on Geometric Applications of Homotopy Theory at the Fields
Institute.

\section{$n$-Fold Categories} \label{nfoldcategories}

In this section we quickly recall the inductive definition of
$n$-fold category, present an equivalent combinatorial definition of
$n$-fold category, discuss completeness and cocompleteness of
$\mathbf{nFoldCat}$, introduce the $n$-fold nerve $N^n$, prove the
existence of its left adjoint $c^n$, and recall the adjunction
$\delta_! \dashv \delta^\ast$.

\begin{defn} \label{defn:nfold_category_inductive}
A {\it small $n$-fold category}
$\mathbb{D}=(\mathbb{D}_0,\mathbb{D}_1)$ is a category object in the
category of small $(n-1)$-fold categories. In detail, $\mathbb{D}_0$
and $\mathbb{D}_1$ are $(n-1)$-fold categories equipped with
$(n-1)$-fold functors
$$\xymatrix@C=3pc{\mathbb{D}_1 \times_{\mathbb{D}_0} \mathbb{D}_1
\ar[r]^-\circ & \mathbb{D}_1 \ar@/^1pc/[r]^s \ar@/_1pc/[r]_t &
\ar[l]|{\lr{u}} \mathbb{D}_0 }$$ that satisfy the usual axioms of a
category. We denote the category of $n$-fold categories by
$\mathbf{nFoldCat}$.
\end{defn}

Since we will always deal with small $n$-fold categories, we leave
off the adjective ``small''. Also, all of our $n$-fold categories
are strict. The following equivalent combinatorial definition of
$n$-fold category is more explicit than the inductive definition.
The combinatorial definition will only be needed in a few places, so
the reader may skip the combinatorial definition if it appears more
technical than one's taste.

\begin{defn} \label{defn:nfold_category_combinatorial}
The data for an {\it $n$-fold category $\bbD$} are
\begin{enumerate}
\item \label{nsets}
sets $\bbD_{\epsilon}$, one for each $\epsilon \in \{0,1\}^n$,
\item  \label{nsourcetarget}
for every $1 \leq i \leq n$ and $\epsilon'\in \{0,1\}^n$ with
$\epsilon_i'=1$ we have {\it source} and {\it target} functions
$$\xymatrix{s^i, t^i \co \bbD_{\epsilon'} \ar[r] & \bbD_\epsilon}$$
where $\epsilon\in \{0,1\}^n$ satisfies $\epsilon_i=0$ and
$\epsilon_j=\epsilon_j'$ for all $j \neq i$ (for ease of notation we
do not include $\epsilon'$ in the notation for $s^i$ and $t^i$,
despite the ambiguity),
\item \label{nunit}
for every $1 \leq i \leq n$ and $\epsilon, \epsilon'\in \{0,1\}^n $
with $\epsilon_i=0$, $\epsilon_i'=1$, and $\epsilon_j=\epsilon_j'$
for all $j \neq i$, we have a {\it unit} $\xymatrix@1{u^i \co
\bbD_\epsilon \ar[r] & \bbD_{\epsilon'}}$,
\item
for every $1 \leq i \leq n$ and $\epsilon, \epsilon'\in \{0,1\}^n $
with $\epsilon_i=0$, $\epsilon_i'=1$, and $\epsilon_j=\epsilon_j'$
for all $j \neq i$, we have a {\it composition}
$$\xymatrix{\bbD_{\epsilon'} \times_{\bbD_{\epsilon}} \bbD_{\epsilon'} \ar[r]^-{\circ^i} & \bbD_{\epsilon'}}.$$
\end{enumerate}
To form an {\it $n$-fold category}, these data are required to satisfy the following axioms.
\begin{enumerate}
\item \label{nsourcetargetcompatibility}
{\it Compatibility of source and target:} for all $1 \leq i \leq n$
and all $1 \leq j \leq n$,
$$s^i s^j = s^j s^i$$
$$t^i t^j = t^j t^i$$
$$s^i t^j = t^j s^i$$
whenever these composites are defined.
\item \label{nunitcompatibility}
{\it Compatibility of units with units:} for all $1 \leq i \leq n$
and all $1 \leq j \leq n$,
$$u^iu^j=u^ju^i$$
whenever these composites are defined.
\item \label{nunitsourcetargetcompatibility}
{\it Compatibility of units with source and target}: for all $1 \leq
i \leq n$ and all $1 \leq j \leq n$,
$$s^iu^j=u^js^i$$
$$t^iu^j=u^jt^i$$
whenever these composites are defined.
\item
{\it Categorical structure:} for every $1 \leq i \leq n$ and
$\epsilon, \epsilon'\in \{0,1\}^n $ with $\epsilon_i=0$,
$\epsilon_i'=1$, and $\epsilon_j=\epsilon_j'$ for all $j \neq i$,
the diagram in $\mathbf{Set}$
$$\xymatrix@C=3pc{\mathbb{D}_{\epsilon'} \times_{\mathbb{D}_\epsilon} \mathbb{D}_{\epsilon'}
\ar[r]^-{\circ^i} & \mathbb{D}_{\epsilon'} \ar@/^1pc/[r]^{s^i} \ar@/_1pc/[r]_{t^i} &
\ar[l]|{\lr{u^i}} \mathbb{D}_\epsilon }$$
is a category.
\item \label{nspecificinterchangelaw}
{\it Interchange law:} For every $i \neq j$ and every $\epsilon \in
\{0,1\}^n$ with $\epsilon_i=1=\epsilon_j$, the compositions
$\circ^i$ and $\circ^j$ can be interchanged, that is, if $w,x,y,z\in
\bbD_\epsilon$, and
$$t^i(w)=s^i(x), \; \; t^i(y)=s^i(z)$$
$$t^j(w)=s^j(y), \; \; t^j(x)=s^j(z),$$
$$\xymatrix@R=1pc@C=1pc{& & \\ & \ar[r]^i \ar[d]_j & \\ & &}
\;\;\;\;\;
\xymatrix{\ar@{-}[r] \ar@{-}[d] \ar@{}[dr]|w & \ar@{-}[r] \ar@{-}[d] \ar@{}[dr]|x & \ar@{-}[d] \\
\ar@{-}[r] \ar@{-}[d] \ar@{}[dr]|y & \ar@{-}[r] \ar@{-}[d]
\ar@{}[dr]|z & \ar@{-}[d] \\ \ar@{-}[r] & \ar@{-}[r] & }$$ then $(z
\circ^j y) \circ^i (x \circ^j w)=(z \circ^i x) \circ^j (y \circ^i
w)$.
\end{enumerate}
We define $\vert \epsilon \vert$ to be the number of 1's in
$\epsilon$, that is
\begin{equation*} \vert \epsilon \vert:=\vert
\{1\leq i \leq n \mid \epsilon_i=1\} \vert =\sum_{i=1}^n \epsilon_i.
\end{equation*} If
$k=\vert \epsilon \vert$, an element of $\bbD_\epsilon$ is called a
{\it $k$-cube}.
\end{defn}

\begin{rmk}
If $\bbD_\epsilon=\bbD_{\epsilon'}$ for all $\epsilon, \epsilon' \in
\{0,1\}^n$ with $\vert \epsilon \vert = \vert \epsilon' \vert$, then
the data \ref{nsets}, \ref{nsourcetarget}, \ref{nunit} satisfying
axioms \ref{nsourcetargetcompatibility}, \ref{nunitcompatibility},
\ref{nunitsourcetargetcompatibility} are an $n$-truncated {\it
cubical complex} in the sense of Section 1 of
\cite{brownhigginscubes}. Compositions and the interchange law are
also similar. The situation of \cite{brownhigginscubes} is {\it edge
symmetric} in the sense that $\bbD_\epsilon=\bbD_{\epsilon'}$ for
all $\epsilon, \epsilon' \in \{0,1\}^n$ with $\vert \epsilon \vert =
\vert \epsilon' \vert$, and the $\vert \epsilon \vert$ compositions
on $\bbD_\epsilon$ coincide with the $\vert \epsilon' \vert$
compositions on $\bbD_{\epsilon'}$. In the present article we study
the non-edge-symmetric case, in the sense that we do {\it not}
require $\bbD_\epsilon$ and $\bbD_{\epsilon'}$ to coincide when
$\vert \epsilon \vert = \vert \epsilon' \vert$, and hence, the
$\vert \epsilon \vert$ compositions on $\bbD_\epsilon$ are not
required to be the same as the $\vert \epsilon' \vert$ compositions
on $\bbD_{\epsilon'}$.
\end{rmk}

\begin{rmk}
The generalized interchange law follows from the interchange law in
\ref{nspecificinterchangelaw}. For example, if we have eight
compatible 3-dimensional cubes arranged as a 3-dimensional cube,
then all possible ways of composing these eight cubes down to one
cube are the same.
\end{rmk}

\begin{prop}
The inductive notion of $n$-fold category in Definition
\ref{defn:nfold_category_inductive} is equivalent to the
combinatorial notion of $n$-fold category in Definition
\ref{defn:nfold_category_combinatorial} in the strongest possible
sense: the categories of such are equivalent.
\end{prop}
\begin{pf}
For $n=1$ the categories are clearly the same. Suppose the proposition
holds for $n-1$ and call the categories $\mathbf{(n-1)FoldCat(ind)}$ and
$\mathbf{(n-1)FoldCat(comb)}$. Then internal categories in
 $\mathbf{(n-1)FoldCat(ind)}$ are equivalent to internal categories in
 $\mathbf{(n-1)FoldCat(comb)}$, while internal categories in
 $\mathbf{(n-1)FoldCat(comb)}$ are the same as $\mathbf{nFoldCat(comb)}$.
\end{pf}

A 2-fold category, that is, a category object in {\bf Cat}, is
precisely a {\it double category} in the sense of Ehresmann. A
double category consists of a set $\bbD_{00}$ of {\it objects}, a
set $\bbD_{01}$ of {\it horizontal morphisms}, a set $\bbD_{10}$ of
{\it vertical morphisms}, and a set $\bbD_{11}$ of {\it squares}
equipped with various sources, targets, and associative and unital
compositions satisfying the interchange law. Several homotopy
theories for double categories were considered in
\cite{fiorepaolipronk1}.

\begin{examp}
There are various standard examples of double categories. To any
category, one can associate the double category of commutative
squares. Any 2-category can be viewed as a double category with
trivial vertical morphisms or as a double category with trivial
horizontal morphisms. To any 2-category, one can also associate the
double category of {\it quintets}: a square is a square of morphisms
inscribed with a 2-cell in a given direction.
\end{examp}

\begin{examp}
In nature, one often finds {\it pseudo double categories}. These are
like double categories, except one direction is a bicategory rather
than a 2-category (see \cite{grandisdouble1} for a more precise
definition). For example, one may consider 1-manifolds,
2-cobordisms, smooth maps, and appropriate squares. Another example
is rings, bimodules, ring maps, and twisted equivariant maps. For
these examples and more, see \cite{grandisdouble1}, \cite{fiore2},
and other articles on double categories listed in the introduction.
\end{examp}

\begin{examp}
Any $n$-category is an $n$-fold category in numerous ways, just like
a 2-category can be considered as a double category in several ways.
\end{examp}

An important method of constructing $n$-fold categories from $n$
ordinary categories is the external product, which is compatible
with the external product of simplicial sets. This was called the
{\it square product} on page 251 of \cite{ehresmannone}.

\begin{defn} \label{externalproduct}
If $\bfC_1, \ldots, \bfC_n$ are small categories, then the {\it
external product} $\bfC_1 \boxtimes \cdots \boxtimes \bfC_n$ is an
$n$-fold category with object set $\Obj \bfC_1 \times \cdots \times
\Obj \bfC_n$. Morphisms in the $i$-th direction are $n$-tuples
$(f_1, \dots, f_n)$ of morphisms in $\bfC_1 \times \cdots \times
\bfC_n$ where all but the $i$-th entry are identities. Squares in the
$ij$-plane are $n$-tuples where all entries are identities except the
$i$-th and $j$-th entries, and so on. An $n$-cube is an $n$-tuple of
morphisms, possibly all non-identity morphisms.
\end{defn}

\begin{prop}
The category $\mathbf{nFoldCat}$ is locally finitely presentable.
\end{prop}
\begin{pf}
We prove this by induction. The category $\mathbf{Cat}$ of small
categories is known to be locally finitely presentable (see for
example \cite{gabrielulmer}). Assume $\mathbf{(n-1)FoldCat}$ is
locally finitely presentable. The category $\mathbf{nFoldCat}$ is
the category of models in $\mathbf{(n-1)FoldCat}$ for a sketch with
finite diagrams. Since $\mathbf{(n-1)FoldCat}$ is locally finitely
presentable, we conclude from Proposition 1.53 of
\cite{adamekrosicky1994} that $\mathbf{nFoldCat}$ is also locally
finitely presentable.
\end{pf}

\begin{prop}
The category $\mathbf{nFoldCat}$ is complete and cocomplete.
\end{prop}
\begin{pf}
Completeness follows quickly, because $\mathbf{nFoldCat}$ is a
category of algebras. For example, the adjunction between $n$-fold graphs and $n$-fold categories
is monadic by the Beck Monadicity Theorem. This means that the algebras for the induced monad
are precisely the $n$-fold categories.

The category $\mathbf{nFoldCat}$ is cocomplete because $\mathbf{nFoldCat}$ is locally finitely presentable.
\end{pf}

The colimits of certain $k$-fold subcategories are the $k$-fold
subcategories of the the colimit. To prove this, we introduce some
notation.

\begin{notation}
Let $\leq$ denote the lexicographic order on $\{0,1\}^n$, and let
$\overline{k}\in \{0,1\}^n$ with $k=|\overline{k}|$. The forgetful
functor $$\xymatrix{U_{\overline{k}}\co \mathbf{nFoldCat} \ar[r] &
\mathbf{kFoldCat}}$$ assigns to an $n$-fold category $\bbD$ the
$k$-fold category consisting of those sets $\bbD_\epsilon$ with
$\epsilon \leq \overline{k}$ and all the source, target, and
identity maps of $\bbD$ between them. If we picture $\bbD$ as an
$n$-cube with $\bbD_\epsilon$'s at the vertices and source, target,
identity maps on the edges, then the $k$-fold subcategory
$U_{\overline{k}}(\bbD)$ is a $k$-face of this $n$-cube. For
example, if $n=2$ and $k=1$, then $U_{\overline{k}}(\bbD)$ is either
the horizontal or vertical subcategory of the double category
$\bbD$.
\end{notation}

\begin{prop} \label{prop:forgetful_admits_right_adjoint}
The forgetful functor $\xymatrix@1{U_{\overline{k}}\co
\mathbf{nFoldCat} \ar[r] & \mathbf{kFoldCat}}$ admits a right
adjoint $R_{\overline{k}}$, and thus preserves colimits: for any
functor $F$ into $\mathbf{nFoldCat}$ we have
$$U_{\overline{k}}\left( \colim F \right) = \colim U_{\overline{k}} F.$$
\end{prop}
\begin{pf}
For a $k$-fold category $\bbE$, the $n$-fold category
$R_{\overline{k}}\bbE$ has
$U_{\overline{k}}R_{\overline{k}}\bbE=\bbE$, in particular the
objects of $R_{\overline{k}}\bbE$ are the same as the objects of
$\bbE$. The other cubes are defined inductively. If $k_i=0$, then
$R_{\overline{k}}\bbE$ has a unique morphism (1-cube) in direction
$i$ between any two objects. Suppose the $j$-cubes of
$R_{\overline{k}}\bbE$ have already been defined, that is
$\left(R_{\overline{k}}\bbE\right)_\epsilon$ has been defined for
all $\epsilon$ with $|\epsilon|=j$. For any $\epsilon$ with
$|\epsilon|=j+1$ and $\epsilon \nleq \overline{k}$, there is a
unique element of $\left(R_{\overline{k}}\bbE\right)_\epsilon$ for
each boundary of $j$-cubes.

The natural bijection
$$\mathbf{kFoldCat}(U_{\overline{k}}\bbD,\bbE)\cong \mathbf{nFoldCat}(\bbD,R_{\overline{k}} \bbE)$$
is given by uniquely extending $k$-fold functors defined on
$U_{\overline{k}}\bbD$ to $n$-fold functors into $R_{\overline{k}}
\bbE$.
\end{pf}

We next introduce the $n$-fold nerve functor, prove that it admits a
left adjoint, and also prove that an $n$-fold natural transformation
gives rise to a simplicial homotopy after pulling back along the
diagonal.

\begin{defn}
The {\it $n$-fold nerve} of an $n$-fold category $\bbD$ is the
multisimplicial set $N^n\bbD$ with $\overline{p}$-simplices
$$(N^n\bbD)_{\overline{p}}:=Hom_{\mathbf{nFoldCat}}([p_1] \boxtimes \cdots \boxtimes
[p_n],\bbD).$$ A $\overline{p}$-simplex is a $\overline{p}$-array of
composable $n$-cubes.
\end{defn}

\begin{rmk}
The $n$-fold nerve is the same as iterating the nerve construction
for internal categories $n$ times.
\end{rmk}

\begin{examp}
The $n$-fold nerve is compatible with external products: $N^n(\bfC_1
\boxtimes \cdots \boxtimes \bfC_n)=N\bfC_1 \boxtimes \cdots
\boxtimes N\bfC_n$. In particular,
$$N^n([m_1]\boxtimes \cdots \boxtimes [m_n])=\Delta[m_1] \boxtimes \cdots \boxtimes \Delta[m_n]=\Delta[m_1, \ldots,
m_n].$$
\end{examp}

\begin{prop}
The functor $\xymatrix@1{N^n\co \mathbf{nFoldCat} \ar[r] &
\mathbf{SSet^n}}$ is fully faithful.
\end{prop}
\begin{pf}
We proceed by induction. For $n=1$, the usual nerve functor is fully
faithful.

Consider now $n>1$, and suppose $$\xymatrix@1{N^{n-1}\co
\mathbf{(n-1)FoldCat} \ar[r] & \mathbf{SSet^{n-1}}}$$ is fully
faithful. We have a factorization
$$\xymatrix{\Cat(\mathbf{(n-1)FoldCat}) \ar@/^2pc/[rr]^{N^n} \ar[r]_-N &
\left[ \Delta^{\op}, \mathbf{(n-1)FoldCat} \right]
\ar[r]_-{N^{n-1}_*} & \left[ \Delta^{\op}, \mathbf{SSet^{n-1}}
\right], }$$ where the brackets mean functor category. The functor
$N$ is faithful, as $(NF)_0$ and $(NF)_1$ are $F_0$ and $F_1$.  It
is also full, for if $\xymatrix@1{F'\co N\bbD \ar[r] & N\bbE}$, then
$F'_0$ and $F'_1$ form an $n$-fold functor with nerve $F'$
(compatibility of $F'$ with the inclusions $\xymatrix{e_{i,i+1}\co
[1] \ar[r] & [m]}$ determines $F'_m$ from $F'_0$ and $F'_1$).

The functor $N^{n-1}_*$ is faithful, since it is faithful at every
degree by hypothesis. If $\xymatrix@1{(G_m')_m \co (N^{n-1}\bbD_m)_m
\ar[r] & (N^{n-1}\bbE_m)_m}$ is a morphism in $\left[ \Delta^{\op},
\mathbf{SSet^{n-1}} \right]$, there exist $(n-1)$-fold functors
$G_m$ such that $N^{n-1}G_m=G_m'$, and these are compatible with the
structure maps for $(\bbD_m)_m$ and $(\bbE_m)_m$ by the faithfulness
of $N^{n-1}$. So $N^{n-1}_*$ is also full.

Finally, $N^n=N^{n-1} \circ N$ is a composite of fully faithful
functors.

This proposition can also be proved using the Nerve Theorem 4.10 of \cite{weber}. For a
direct proof in the case $n=2$, see \cite{fiorepaolifundamental}.
\end{pf}

\begin{prop}
The $n$-fold nerve functor $N^n$ admits a left adjoint $c^n$ called
fundamental $n$-fold category or $n$-fold categorification.
\end{prop}
\begin{pf}
The functor $N^n$ is defined as the singular functor associated to
an inclusion. Since $\mathbf{nFoldCat}$ is cocomplete, a left
adjoint to $N^n$ is obtained by left Kan extending along the Yoneda
embedding. This is the Lemma from Kan about singular-realization
adjunctions.
\end{pf}

\begin{examp} \label{examp:categorification_and_external_products}
If $X_1, \ldots, X_n$ are simplicial sets, then $$c^n(X_1\boxtimes
\cdots \boxtimes X_n)=cX_1\boxtimes \cdots \boxtimes cX_n$$ where
$c$ is ordinary categorification. The symbol $\boxtimes$ on the left
means external product of simplicial sets, and the symbol
$\boxtimes$ on the right means external product of categories as in
Definition \ref{externalproduct}. For a proof in the case $n=2$, see
\cite{fiorepaolifundamental}.
\end{examp}

Since the nerve of a natural transformation is a simplicial
homotopy, we expect the diagonal of the $n$-fold nerve of an
$n$-fold natural transformation to be a simplicial homotopy.

\begin{defn} \label{defn_nfold_nat_transf}
An {\it $n$-fold natural transformation} $\xymatrix@1{\alpha \co F
\ar@{=>}[r] & G }$ between $n$-fold functors $\xymatrix@1{F,G\co
\bbD \ar[r] & \bbE }$ is an $n$-fold functor
$$\xymatrix{\alpha\co \bbD \times [1]^{\boxtimes n} \ar[r] & \bbE}$$
such that $\alpha\vert_{\bbD \times \{0\}}$ is $F$ and
$\alpha\vert_{\bbD \times \{1\}}$ is $G$.
\end{defn}

Essentially, an $n$-fold natural transformation associates to an
object an $n$-cube with source corner that object, to a morphism in
direction $i$ a square in direction $ij$ for all $j\neq i$ in $1
\leq j \leq n$, to an $ij$-square a 3-cube  in direction $ijk$ for
all $k \neq i,j$ in $1 \leq k \leq n$ etc, and these are
appropriately functorial, natural, and compatible.

\begin{examp} \label{examp:n_naturaltransfs_yield_an_nfold_naturaltransf}
If $\xymatrix@1{\alpha_i\co\bfC_i \times [1] \ar[r] & \bfC_i' }$ are
ordinary natural transformations between ordinary functors for $1
\leq i \leq n$, then $\alpha_1 \boxtimes \cdots \boxtimes \alpha_n$
is an $n$-fold natural transformation because of the isomorphism
$$(\bfC_1 \times [1]) \boxtimes \cdots \boxtimes (\bfC_n \times [1])
\cong (\bfC_1 \boxtimes \cdots \boxtimes \bfC_n) \times ([1]
\boxtimes \cdots \boxtimes [1]).$$
\end{examp}

\begin{prop} \label{nfoldnat_gives_simplicial_homotopy}
Let $\xymatrix@1{\alpha \co \bbD \times [1]^{\boxtimes n} \ar[r] &
\bbE}$ be an $n$-fold natural transformation. Then
$(\delta^*N^n\alpha)\circ (1_{\delta^*N^n \bbD} \times d)$ is a
simplicial homotopy from $\delta^*(N^n \alpha\vert_{\bbD \times
\{0\}})$ to $\delta^*(N^n \alpha\vert_{\bbD \times \{1\}})$.
\end{prop}
\begin{pf}
We have the diagonal of the $n$-fold nerve of $\alpha$
$$\xymatrix@1@C=4pc{\delta^*(N^n \bbD) \times
\delta^*(N^n[1]^{\boxtimes n}) \ar[r]^-{\delta^*N^n \alpha} &
\delta^*N^n \bbE},$$ which we then precompose with $1_{\delta^*N^n
\bbD}\times d$ to get
$$\xymatrix@1@C=4pc{(\delta^*N^n \bbD) \times \Delta[1]
\ar[r]^-{1_{\delta^*N^n \bbD} \times d} & \delta^*(N^n \bbD) \times
\Delta[1]^{\times n} \ar[r]^-{\delta^*N^n \alpha} & \delta^*N^n
\bbE}.$$
\end{pf}

Lastly, we consider the behavior of $c^n$ on the image of the left
adjoint $\delta_!$. The diagonal functor
$$\xymatrix{\delta \co\Delta \ar[r] & \Delta^n}$$
$$[m] \mapsto ([m],\ldots, [m])$$
induces $\xymatrix@1{\delta^\ast\co \mathbf{SSet^n} \ar[r] &
\mathbf{SSet}}$ by precomposition. The functor $\delta^\ast$ admits
both a left and right adjoint by Kan extension. The left adjoint
$\delta_!$ is uniquely characterized by two properties:
\begin{enumerate}
\item
$\delta_!(\Delta[m])=\Delta[m,\ldots,m]$,
\item
$\delta_!$ preserves colimits.
\end{enumerate}
Thus,
\begin{equation*}
\delta_!X=\delta_!(\underset{\Delta[m] \rightarrow X}{\colim}
\Delta[m])=\underset{\Delta[m] \rightarrow
X}{\colim}\delta_!\Delta[m]=\underset{\Delta[m] \rightarrow
X}{\colim} \Delta[m,\ldots,m]
\end{equation*}
where the colimit is over the simplex category of the simplicial set
$X$. Further, since $c^n$ preserves colimits, we have
\begin{equation*}
c^n\delta_!X=\underset{\Delta[m] \rightarrow X}{\colim}
c^n\Delta[m,\ldots,m]=\underset{\Delta[m] \rightarrow X}{\colim}
[m]\boxtimes \cdots \boxtimes [m].
\end{equation*}

Clearly, $c^n \delta_! \Delta[m]=[m] \boxtimes \cdots \boxtimes
[m]$. The calculation of $c^n \delta_! \Sd^2 \Delta[m]$ and $c^n
\delta_! \Sd^2 \Lambda^k[m]$ is not as simple, because external
product does not commute with colimits. We will give a general
procedure of calculating the $n$-fold categorification of nerves of
certain posets in Section \ref{sectionnfolddecompositions}.

\section{Barycentric Subdivision and Decomposition of $\bfP\Sd\Delta[m]$}\label{barycentric}

\begin{figure}
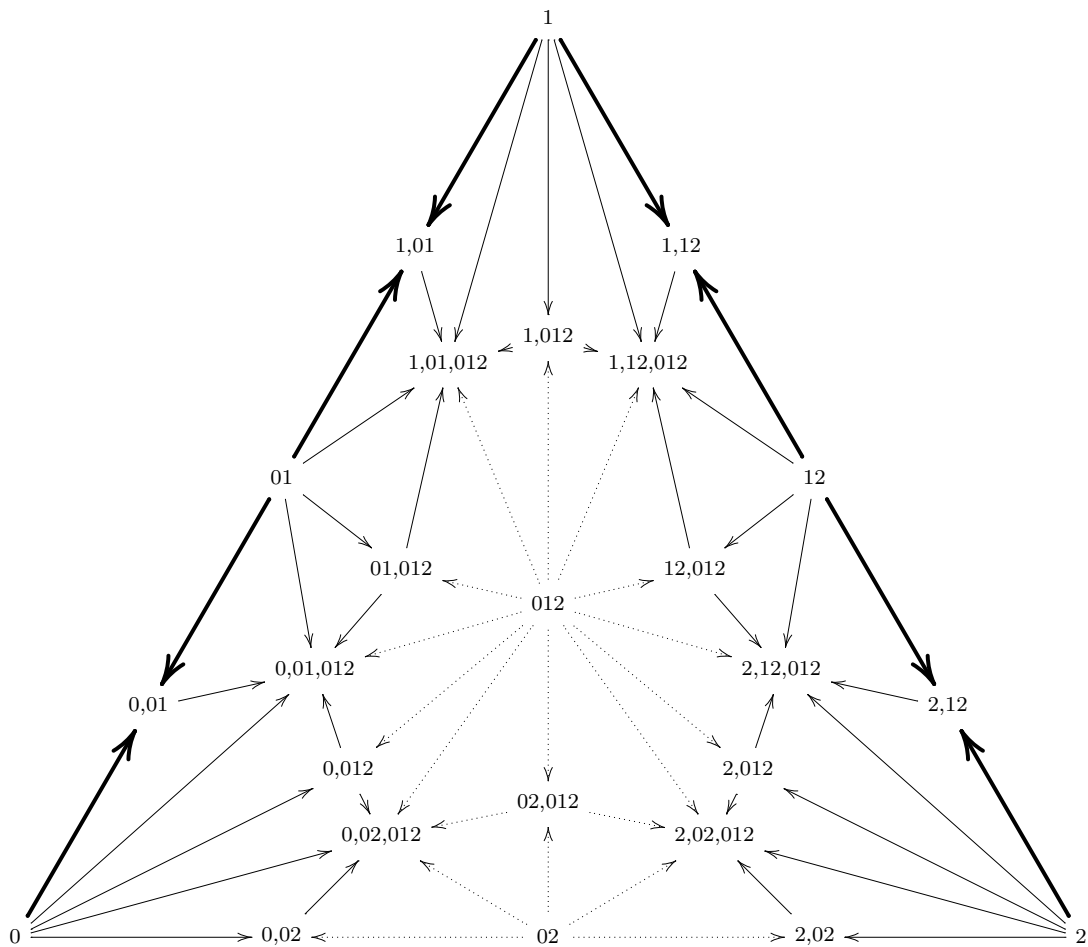

\begin{center}
    \def\objectstyle{\scriptstyle}
    \xy
    0;/r.21pc/:
    (0,0)*{\halo{02}}="1";
    (80,0)*{\halo{2}}="2";
    (-80,0)*{\halo{0}}="3";
    (40,0)*{\halo{2,02}}="4";
    (-40,0)*{\halo{0,02}}="5";
    (0,50)*{\halo{012}}="6";
    (-60,34.5)*{\halo{0,01}}="7";
    (60,34.5)*{\halo{2,12}}="8";
    (-35,40)*{\halo{0,01,012}}="9";
    (35,40)*{\halo{2,12,012}}="10";
    (-30,25)*{\halo{0,012}}="11";
    (30,25)*{\halo{2,012}}="12";
    (-25,15)*{\halo{0,02,012}}="13";
    (25,15)*{\halo{2,02,012}}="14";
    (0,20)*{\halo{02,012}}="15";
    (22,55)*{\halo{12,012}}="16";
    (-22,55)*{\halo{01,012}}="17";
    (40,69)*{\halo{12}}="18";
    (-40,69)*{\halo{01}}="19";
    (15,86)*{\halo{1,12,012}}="20";
    (-15,86)*{\halo{1,01,012}}="21";
    (0,90)*{\halo{1,012}}="22";
    (-20,103.5)*{\halo{1,01}}="23";
    (20,103.5)*{\halo{1,12}}="24";
    (0,138)*{\halo{1}}="25";
    {\ar@{->}"3";"5"};
    {\ar@{->}"3";"13"};
    {\ar@{->}"3";"11"};
    {\ar@{->}"3";"9"};
        {\ar@*{[|1.5pt]}"3";"7"};
    {\ar@{->}"2";"4"};
    {\ar@{->}"2";"14"};
    {\ar@{->}"2";"12"};
    {\ar@{->}"2";"10"};
        {\ar@*{[|1.5pt]}"2";"8"};
    {\ar@{.>}"1";"5"};
    {\ar@{.>}"1";"13"};
    {\ar@{.>}"1";"15"};
    {\ar@{.>}"1";"14"};
    {\ar@{.>}"1";"4"};
    {\ar@{.>}"6";"22"};
    {\ar@{.>}"6";"20"};
    {\ar@{.>}"6";"16"};
    {\ar@{.>}"6";"10"};
    {\ar@{.>}"6";"12"};
    {\ar@{.>}"6";"14"};
    {\ar@{.>}"6";"15"};
    {\ar@{.>}"6";"13"};
    {\ar@{.>}"6";"11"};
    {\ar@{.>}"6";"9"};
    {\ar@{.>}"6";"17"};
    {\ar@{.>}"6";"21"};
    {\ar@{.>}"15";"13"};
    {\ar@{.>}"15";"14"};
    {\ar@{->}"11";"13"};
    {\ar@{->}"11";"9"};
    {\ar@{->}"12";"10"};
    {\ar@{->}"12";"14"};
    {\ar@{->}"7";"9"};
    {\ar@{->}"8";"10"};
        {\ar@*{[|1.5pt]}"19";"7"};
    {\ar@{->}"19";"9"};
    {\ar@{->}"19";"17"};
    {\ar@{->}"19";"21"};
        {\ar@*{[|1.5pt]}"19";"23"};
        {\ar@*{[|1.5pt]}"18";"8"};
    {\ar@{->}"18";"10"};
    {\ar@{->}"18";"16"};
    {\ar@{->}"18";"20"};
        {\ar@*{[|1.5pt]}"18";"24"};
    {\ar@{->}"23";"21"};
    {\ar@{->}"24";"20"};
    {\ar@{->}"17";"9"};
    {\ar@{->}"17";"21"};
    {\ar@{->}"16";"20"};
    {\ar@{->}"16";"10"};
        {\ar@*{[|1.5pt]}"25";"23"};
    {\ar@{->}"25";"21"};
    {\ar@{->}"25";"22"};
    {\ar@{->}"25";"20"};
    {\ar@*{[|1.5pt]}"25";"24"};
    {\ar@{->}"22";"20"};
    {\ar@{->}"22";"21"};
    {\ar@{->}"5";"13"};
    {\ar@{->}"4";"14"};
    \endxy
\end{center}
\caption{Decomposition of the poset $\bfP\Sd\Delta[2]$. The dark
arrows form the poset $\bfP\Sd\Lambda^1[2]$, while its up-closure
$\Out$ consists of all solid arrows. The poset $\Cen$ consists of
all the triangles emanating from 012; these triangles all have two
dotted sides emanating from 012. The poset $\Comp$ consists of the
four triangles at the bottom emanating from 02; these four trinagles
each have two dotted sides emanating from 02. The geometric
realization of all triangles with at least two dotted edges, namely
$|N(\Comp \cup \Cen)|$, is topologically deformation retracted onto
the solid part of its boundary.} \label{subdivisionfigure}
\end{figure}

The adjunction
\begin{equation}
\xymatrix@C=4pc{\mathbf{SSet} \ar@{}[r]|{\perp} \ar@/^1pc/[r]^{\Sd}
& \ar@/^1pc/[l]^{\Ex} \mathbf{SSet}}
\end{equation}
between barycentric subdivision $\Sd$ and Kan's functor $\Ex$ is
crucial to Thomason's transfer from {\bf Cat} to {\bf SSet}. We will
need a good understanding of subdivision for the Thomason structure
on $\mathbf{nFoldCat}$ as well, so we recall it in this section.
Explicit descriptions of certain subsimplices of the double
subdivisions $\Sd^2\Lambda^k[m]$, $\Sd^2\partial\Delta[m]$, and
$\Sd^2\Delta[m]$ will be especially useful later. In Proposition
\ref{upcloseddecomposition}, we present a decomposition of the poset
$\bfP\Sd\Delta[m]$, which is pictured in Figure
\ref{subdivisionfigure} for the case $m=2$ and $k=1$. The nerve of
the poset $\bfP\Sd\Delta[m]$ is of course $\Sd^2\Delta[m]$. This
decomposition allows us to describe a deformation retraction of part
of $|\Sd^2\Delta[m]|$ in a very controlled way (Proposition
\ref{deformationretract}). In particular, each $m$-subsimplex is
deformation retracted onto one of its faces. This allows us to do a
deformation retraction of the $n$-fold categorifications as well in
Corollary \ref{nfolddeformationretract}. These preparations are
essential for verifying the pushout-axiom in Kan's Lemma on Transfer
of Model Structures.

We begin now with our recollection of barycentric subdivision. The
simplicial set $\Sd\Delta[m]$ is the nerve of the poset
$\bfP\Delta[m]$ of non-degenerate simplices of $\Delta[m]$. The
ordering is the face relation. Recall that the poset $\bfP\Delta[m]$
is isomorphic to the poset of nonempty subsets of $[m]$ ordered by
inclusion. Thus a $q$-simplex $v$ of $\Sd\Delta[m]$ is a tuple
$(v_0, \dots, v_q)$ of nonempty subsets of $[m]$  such that $v_i$ is
a subset of $v_{i+1}$ for all $0\leq i\leq q-1$. For example, the
tuple
\begin{equation} \label{vexample}
(\{0\},\{0,2\},\{0,1,2,3\})
\end{equation}
is a 2-simplex of $\Sd\Delta[3]$. A $p$-simplex $u$ is a {\it face}
of a $q$-simplex $v$ in $\Sd\Delta[m]$ if and only if
\begin{equation} \label{facerelation}
\{u_0,\dots, u_p\} \subseteq \{v_0,\dots, v_q\}.
\end{equation}
For example the 1-simplex
\begin{equation} \label{uexample}
(\{0\},\{0,1,2,3\})
\end{equation}
is a face of the 2-simplex in equation \eqref{vexample}. A face that
is a 0-simplex is called a {\it vertex}. The vertices of $v$ are
written simply as $v_0, \dots, v_q$. A $q$-simplex $v$ of
$\Sd\Delta[m]$ is non-degenerate if and only if all $v_i$ are
distinct. The simplices in equations (\ref{vexample}) and
(\ref{uexample}) are both non-degenerate.

The barycentric subdivision of a general simplicial set $K$ is
defined in terms of the barycentric subdivisions $\Sd\Delta[m]$ that
we have just recalled.
\begin{defn}
The {\it barycentric subdivision} of a simplicial set $K$ is
$$\underset{\Delta[n] \rightarrow K}{\colim}
\Sd\Delta[n]$$ where the colimit is indexed over the category of
simplices of $K$.
\end{defn}
The right adjoint to $\Sd$ is the $\Ex$ functor of Kan, and is
defined in level $m$ by
$$(\Ex X)_m=\mathbf{SSet}(\Sd \Delta[m], X).$$

As pointed out on page 311 of \cite{thomasonCat}, there is a
particularly simple description of $\Sd K$ whenever $K$ is a
classical simplicial complex each of whose simplices has a linearly
ordered vertex set compatible with face inclusion. In this case,
$\Sd K$ is the nerve of the poset $\bfP K$ of non-degenerate
simplices of $K$. The cases $K=\Sd \Delta[m], \Lambda^k[m], \Sd
\Lambda^k[m], \partial \Delta[m],$ and $\Sd \partial \Delta[m]$ are
of particular interest to us.

We first consider the case $K=\Sd \Delta[m]$ in order to describe
the simplicial set $\Sd^2\Delta[m]$. This is the nerve of the poset
$\bfP\Sd\Delta[m]$ of non-degenerate simplices of $\Sd\Delta[m]$. A
$q$-simplex of $\Sd^2\Delta[m]$ is a sequence $V=(V_0,\dots, V_q)$
where each $V_i=(v_0^i, \dots, v^i_{r_i})$ is a non-degenerate
simplex of $\Sd\Delta[m]$ and $V_{i-1} \subseteq V_i$. For example,
\begin{equation} \label{Vexample}
\left(\;\;(\{01\}),\;\;(\{0\},\{01\}),\;\;(\{0\},\{01\},\{012\})\;\;
\right)
\end{equation}
is a 2-simplex in $\Sd^2 \Delta[2]$. A $p$-simplex $U$ is a face of
a $q$-simplex $V$ in $\Sd^2\Delta[m]$ if and only if
\begin{equation} \label{OrderOnDoubleSubdivision}
\{U_0,\dots, U_p\} \subseteq \{V_0,\dots, V_q\}.
\end{equation}
For example, the 1-simplex
\begin{equation} \label{Uexample}
\left(\;\;(\{01\}),\;\;(\{0\},\{01\},\{012\})\;\; \right)
\end{equation}
is a subsimplex of the 2-simplex in equation \eqref{Vexample}. The
vertices of $V$ are $V_0,\dots, V_q$. A $q$-simplex $V$ of
$\Sd^2\Delta[m]$ is non-degenerate if and only if all $V_i$ are
distinct. The simplices in equations (\ref{Vexample}) and
(\ref{Uexample}) are both non-degenerate. Figure
\ref{subdivisionfigure} displays the poset $\bfP\Sd\Delta[m]$, the
nerve of which is $\Sd^2\Delta[m]$.

Next we consider $K=\Lambda^k[m]$ in order to describe $\Sd
\Lambda^k[m]$ as the nerve of the poset $\bfP\Lambda^k[m]$ of
non-degenerate simplices of $\Lambda^k[m]$. The simplicial set
$\Lambda^k[m]$ is the smallest simplicial subset of $\Delta[m]$
which contains all non-degenerate simplices of $\Delta[m]$ except the
sole $m$-simplex $1_{[m]}$ and the $(m-1)$-face opposite the vertex
$\{k\}$. The $n$-simplices of $\Lambda^k[m]$ are
\begin{equation} \label{simplicesofhorn}
(\Lambda^k[m])_n=\{\xymatrix@1{f:[n] \ar[r] & [m]}| \, \im f \nsupseteq [m]\backslash \{k\}  \}.
\end{equation}
A $q$-simplex $(v_0, \dots, v_q)$ of $\Sd \Delta[m]$ is in $\Sd
\Lambda^k[m]$ if and only if each $v_i$ is a face of $\Lambda^k[m]$.
More explicitly, $(v_0, \dots, v_q)$ is in $\Sd \Lambda^k[m]$ if and
only if $|v_q| \leq m$ and in case of equality $k \in v_q$. This
follows from equation \eqref{simplicesofhorn}. Similarly, a
$q$-simplex $V$ in $\Sd^2 \Delta[m]$ is in $\Sd^2 \Lambda^k[m]$ if
and only if all $v^i_j$ are faces of $\Lambda^k[m]$. This is the
case if and only if for all $0 \leq i \leq q$, $|v^i_{r_i}|\leq m$
and in case of equality $k \in v^i_{r_i}$. This, in turn, is the
case if and only if $|v^q_{r_q}|\leq m$ and in case of equality $k
\in v^q_{r_q}$. See again Figure \ref{subdivisionfigure}.

Lastly, we similarly describe $\Sd \partial \Delta[m]$ and $\Sd^2
\partial \Delta[m]$. The simplicial set $\partial \Delta[m]$ is the
simplicial subset of $\Delta[m]$ obtained by removing the sole
$m$-simplex $1_{[m]}$. A $q$-simplex $(v_0, \dots, v_q)$ of $\Sd
\Delta[m]$ is in $\Sd \partial \Delta[m]$ if and only if $v_q\neq
\{0,1, \dots, m\}$. A $q$-simplex $V$ of $\Sd^2 \Delta[m]$ is in
$\Sd^2 \partial \Delta[m]$ if and only if $v^i_{r_i}\neq\{0,1,
\dots, m\}$ for all $0 \leq i \leq q$, which is the case if and only
if $v^q_{r_q}\neq\{0,1, \dots, m\}$. See again Figure
\ref{subdivisionfigure}.

\begin{rmk} \label{gluingofsimplices}
Also of interest to us is the way that the non-degenerate
$m$-simplices of $\Sd^2 \Delta[m]$ are glued together along their
$(m-1)$-subsimplices. In the following, let $V=(V_0, \dots, V_m)$ be
a non-degenerate $m$-simplex of $\Sd^2 \Delta[m]$. Each $V_i=(v_0^i,
\dots, v^i_{r_i})$ is then a distinct non-degenerate simplex of
$\Sd\Delta[m]$. See Figure \ref{subdivisionfigure} for intuition.
\begin{enumerate}
\item \label{gluingofsimplicesi}
Then $r_i=i$, $|V_i|=i+1$, and hence also $v^m_m=\{0,1, \dots, m\}$.
\item \label{gluingofsimplicesii}
If $v^{m-1}_{m-1} \neq\{0,1, \dots, m\}$, then the $m$-th face
$(V_0, \dots, V_{m-1})$ of $V$ is not shared with
any other non-degenerate $m$-simplex $V'$ of $\Sd^2 \Delta[m]$.  \\
{\it Proof:} If $v^{m-1}_{m-1} \neq\{0,1, \dots, m\}$, then the
$(m-1)$-simplex $(V_0, \dots, V_{m-1})$ lies in $\Sd^2
\partial \Delta[m]$ by the description of $\Sd^2
\partial \Delta[m]$ above, and hence does not lie in any
other non-degenerate $m$-simplex $V'$ of $\Sd^2\Delta[m]$.
\item \label{gluingofsimplicesiii}
If $v^{m-1}_{m-1}=\{0,1, \dots, m\}$, then the $m$-th face $(V_0,
\dots, V_{m-1})$ of $V$ is shared with one other non-degenerate
$m$-simplex $V'$ of $\Sd^2 \Delta[m]$.
\\ {\it Proof:}
If $v^{m-1}_{m-1}=\{0,1, \dots, m\}$, then there exists a unique
$0\leq i \leq m-1$ with $v_i^{m-1} \backslash
v_{i-1}^{m-1}=\{a,a'\}$ with $a \neq a'$ (since the sequence
$v_0^{m-1},v_1^{m-1}, \dots, v_{m-1}^{m-1}=\{0,1, \dots, m\}$ is
strictly ascending). Here we define $v^{m-1}_{i-1}=\emptyset$
whenever $i=0$. Thus, the $(m-1)$-simplex $(V_0, \dots, V_{m-1})$ is
also a face of the non-degenerate $m$-simplex $V'$ where
$$V'_\ell=V_\ell \text{\hspace{.5in} for $0 \leq \ell \leq m-1$}$$
$$V'_m=(v_0^{m-1}, \dots, v^{m-1}_{i-1}, v^{m-1}_{i-1} \cup \{a'\}, v^{m-1}_i, \dots ,
v^{m-1}_{m-1}),$$ where we also have
$$V_m=(v_0^{m-1}, \dots, v^{m-1}_{i-1}, v^{m-1}_{i-1} \cup \{a\}, v^{m-1}_i, \dots ,
v^{m-1}_{m-1}).$$
\item \label{gluingofsimplicesiv}
If $0\leq j \leq m-1$, then $V$ shares its $j$-th face $(\dots,
\hat{V}_j,\dots, V_m)$ with one other non-degenerate $m$-simplex
$V'$ of $\Sd^2 \Delta[m]$. \\ {\it Proof:} Since $|V_i|=i+1$, we
have $V_{j+1} \backslash V_{j-1}=\{v,v'\}$ with $v\neq v'$ (we
define $V_{j-1}=\emptyset$ whenever $j=0$). Then $(\dots,
\hat{V}_j,\dots, V_m)$ is shared by the two non-degenerate
$m$-simplices
$$V=(V_0, \dots, V_{j-1}, V_{j-1} \cup \{v\}, V_{j+1}, \dots, V_m)$$
$$V'=(V_0, \dots, V_{j-1}, V_{j-1} \cup \{v'\}, V_{j+1}, \dots, V_m)$$
and no others.
\end{enumerate}
\end{rmk}

After this brief discussion of how the non-degenerate $m$-simplices
of $\Sd^2\Delta[m]$ are glued together, we turn to some comments
about the relationships between the second subdivisions of
$\Lambda^k[m]$, $\partial\Delta[m]$, and $\Delta[m]$. Since the
counit $\xymatrix@1{cN \ar@{=>}[r] & 1_{\mathbf{Cat}} }$ is a
natural isomorphism\footnote{The nerve functor is fully faithful, so
the counit is a natural isomorphism by IV.3.1 of
\cite{maclaneworking}}, the categories $c\Sd^2\Lambda^k[m]$,
$c\Sd^2\partial\Delta[m]$, and $c\Sd^2\Delta[m]$ are respectively
the posets $\bfP\Sd\Lambda^k[m]$, $\bfP\Sd\partial\Delta[m]$, and
$\bfP\Sd\Delta[m]$ of non-degenerate simplices. Moreover, the induced
functors
$$\xymatrix@1{c\Sd^2\Lambda^k[m] \ar[r] &
c\Sd^2\Delta[m]} \hspace{.5in} \xymatrix@1{c\Sd^2\partial\Delta[m]
\ar[r] & c\Sd^2\Delta[m] }$$ are simply the poset inclusions
$$\xymatrix@1{\bfP\Sd\Lambda^k[m] \ar[r] &
\bfP\Sd\Delta[m]} \hspace{.5in} \xymatrix@1{\bfP\Sd\partial\Delta[m]
\ar[r] & \bfP\Sd\Delta[m] }.$$

The down-closure of $\bfP\Sd\Lambda^k[m]$ in $\bfP\Sd\Delta[m]$ is
easily described.

\begin{prop} \label{Lambdadownclosed}
The subposet $\bfP\Sd\Lambda^k[m]$ of $\bfP\Sd\Delta[m]$ is
down-closed.
\end{prop}
\begin{pf}
A $q$-simplex $(v_0, \dots, v_q)$ of $\Sd \Delta[m]$ is in $\Sd
\Lambda^k[m]$ if and only if $|v_q| \leq m$ and in case of equality
$k \in v_q$. If $(v_0, \dots, v_q)$ has this property, then so do
all of its subsimplices.
\end{pf}

The rest of this section is dedicated to a decomposition of
$\bfP\Sd\Delta[m]$ into the union of three up-closed subposets:
$\Comp$, $\Cen$, and $\Out$. This culminates in Proposition
\ref{upcloseddecomposition}, and will be used in the construction of
the retraction in Section \ref{retractionsection} as well as the
transfer proofs in Sections \ref{Thomasonsection} and
\ref{Thomasonnfoldsection}. The reader is encouraged to compare with
Figure \ref{subdivisionfigure} throughout. We begin by describing
these posets. The poset $\Out$ is the up-closure of
$\bfP\Sd\Lambda^k[m]$ in $\bfP\Sd\Delta[m]$. Although $\Out$ depends
on $k$ and $m$, we omit these letters from the notation for
readability.
\begin{prop} \label{outer}
Let $\Out$ denote the smallest up-closed subposet of
$\bfP\Sd\Delta[m]$ which contains $\bfP\Sd\Lambda^k[m]$.
\begin{enumerate}
\item \label{outeri}
The subposet $\Out$ consists of those $(v_0, \dots, v_q)\in \bfP
\Sd\Delta[m]$ such that there exists a $(u_0, \dots, u_p) \in
\bfP\Sd \Lambda^k[m]$ with
$$\{u_0,\dots, u_p\} \subseteq \{v_0,\dots, v_q\}.$$ In particular,
$(v_0, \dots, v_q)\in \bfP \Sd\Delta[m]$ is in $\Out$ if and only if
some $v_i$ satisfies $|v_i|\leq m$ and in case of equality $k \in
v_i$.
\item \label{outerii}
Define a functor $\xymatrix@1{r\co \Out \ar[r] &
\bfP\Sd\Lambda^k[m]}$ by $r(v_0, \dots, v_q):=(u_0, \dots, u_p)$
where $(u_0, \dots, u_p)$ is the maximal subset
$$\{u_0,\dots, u_p\} \subseteq \{v_0,\dots, v_q\}$$
that is in $\bfP\Sd \Lambda^k[m]$. Let $\xymatrix@1{\text{\rm inc}
\co \bfP\Sd \Lambda^k[m] \ar[r] & \Out}$ be the inclusion. Then
$r\circ  \text{\rm inc} = 1_{\bfP\Sd \Lambda^k[m]}$ and there is a
natural transformation $\xymatrix@1{\alpha \co \text{\rm inc}\circ r
\ar@{=>}[r] & 1_{\Out}}$ which is the identity morphism on objects
of $\bfP\Sd \Lambda^k[m]$. Consequently, $\vert \bfP\Sd \Lambda^k[m]
\vert$ is a deformation retract of $\vert \Out \vert$. See Figure
\ref{subdivisionfigure} for a geometric picture.
\end{enumerate}
\end{prop}
\begin{pf}
\begin{enumerate}
\item
An element of $\bfP\Sd\Delta[m]$ is in the up-closure of
$\bfP\Sd\Lambda^k[m]$ if and only if it lies above some element of
$\bfP\Sd\Lambda^k[m]$, and the order is the face relation as in
equation \eqref{facerelation}. For the last part, we use the
observation that $(u_0, \dots, u_p) \in \bfP\Sd \Lambda^k[m]$ if and
only if $|u_p|\leq m$ and in the case of equality $k \in u_p$, as in
the discussion after \eqref{simplicesofhorn}, and also the fact that
$(u_j)\leq(u_0, \dots, u_p)$.
\item
For $(v_0, \dots, v_q) \in \Out$, we define $\alpha(v_0, \dots,
v_q)$ to be the unique arrow in $\Out$ from $r(v_0, \dots, v_q)$ to
$(v_0, \dots, v_q)$. Naturality diagrams must commute, since $\Out$
is a poset. The rest is clear.
\end{enumerate}
\end{pf}

The following trivial remark will be of use later.

\begin{rmk} \label{Deltafreecomposites}
Since $\bfP\Sd\Lambda^k[m]$ is down-closed by Proposition
\ref{Lambdadownclosed}, any morphism of $\bfP\Sd\Delta[m]$ that ends
in $\bfP\Sd\Lambda^k[m]$ must also be contained in
$\bfP\Sd\Lambda^k[m]$. Since $\Out$ is the up-closure of the poset
$\bfP\Sd\Lambda^k[m]$ in $\bfP\Sd\Delta[m]$, any morphism that
begins in $\bfP\Sd\Lambda^k[m]$ ends in $\Out$.
\end{rmk}

We can similarly characterize the up-closure $\Cen$ of
$(\{0,1,\dots,m\})$ in $\bfP\Sd\Delta[m]$. We call a non-degenerate
$m$-simplex of $\Sd^2\Delta[m]$ a {\it central $m$-simplex} if it
has $(\{0,1, \dots, m\})$ as its $0$-th vertex.
\begin{prop} \label{center}
The smallest up-closed subposet $\Cen$ of $\bfP\Sd\Delta[m]$ which
contains $(\{0,1,\dots,m\})$ consists of those $(v_0, \dots, v_q)\in
\bfP \Sd\Delta[m]$ such that $v_q=\{0,1,\dots,m\}$. The nerve
$N\Cen$ consists of all central $m$-simplices of $\Sd^2\Delta[m]$
and all their faces.  A $q$-simplex $(V_0,\dots, V_q)$ of
$\Sd^2\Delta[m]$ is in $N\Cen$ if and only if $v^i_{r_i}=
\{0,1,\dots,m\}$ for all $0\leq i \leq q$.
\end{prop}
For example, the 2-simplex
\begin{equation} \label{central2simplex}
\left(\;\; (\{012\}), \;\; (\{01\},\{012\}), \;\;
(\{0\},\{01\},\{012\}) \;\; \right)
\end{equation}
is a central 2-simplex of $\Sd^2 \Delta[2]$ and the 1-simplex
\begin{equation} \label{faceofcentral2simplex}
\left( \;\; (\{01\},\{012\}), \;\; (\{0\},\{01\},\{012\}) \;\;
\right)
\end{equation}
is in $N\Cen$, as it is a face of the 2-simplex in equation
\eqref{central2simplex}. A glance at Figure \ref{subdivisionfigure}
makes all of this apparent.

\begin{rmk}
We need to understand more thoroughly the way the central
$m$-simplices are glued together in $N\Cen$. Suppose $V$ is a
central $m$-simplex, so that $v^i_i=\{0,1, \dots, m\}$ for all $0
\leq i \leq m$ by Proposition \ref{center}. From the description of
$V'$ in Remark \ref{gluingofsimplices} \ref{gluingofsimplicesiii}
and \ref{gluingofsimplicesiv}, and also Proposition \ref{center}
again, we see for $j=1, \dots, m$ that the neighboring
non-degenerate $m$-simplex $V'$ containing the $(m-1)$-face $(V_0,
\dots, \hat{V}_j, \dots)$ of $V$ is also central. The face $(V_1,
\dots, V_m)$ of $V$ opposite $V_0=(\{0,1,\dots,m\})$, is not shared
with any other central $m$-simplex as every central $m$-simplex has
$\{0, \dots, m\}$ as its $0$-th vertex. Thus, each central
$m$-simplex $V$ shares exactly $m$ of its $(m-1)$-faces with other
central $m$-simplices. A glance at Figure \ref{subdivisionfigure}
shows that the central simplices fit together to form a $2$-ball.
More generally, the central $m$-simplices of $\Sd^2\Delta[m]$ fit
together to form an $m$-ball with center vertex $\{0, \dots, m\}$.
\end{rmk}




There is still one last piece of $\bfP\Sd\Delta[m]$ that we discuss,
namely $\Comp$.

\begin{prop} \label{upclosedcomp}
Let $0 \leq k \leq m$. The smallest up-closed subposet $\Comp$  of
$\bfP\Sd\Delta[m]$ that contains the object $(\{0,1, \dots, \hat{k},
\dots, m\})$ consists of those $(v_0, \dots, v_q)\in
\bfP\Sd\Delta[m]$ with
$$\{0,1, \dots, \hat{k},
\dots, m\}\in \{v_0, \dots, v_q\}.$$
\end{prop}




We describe how the non-degenerate $m$-simplices of $N\Comp$ are
glued together in terms of collections $C^\ell$ of non-degenerate
$m$-simplices. A non-degenerate $m$-simplex $V \in
N_m\bfP\Sd\Delta[m]$ is in $N_m\Comp$ if and only if each $V_0,
\dots, V_m$ is in $\Comp$, and this is the case if and only if
$V_0=(\{0, \dots, \hat{k}, \dots, m\})$ (recall $|V_i|=i+1$ and
Proposition \ref{upclosedcomp}). For $1 \leq \ell \leq m$, we let
$C^\ell$ denote the set of those non-degenerate $m$-simplices $V$ in
$N_m\Comp$ which have their first $\ell$ vertices $V_0, \dots,
V_{\ell-1}$ on the $k$-th face of $|\Delta[m]|$. A non-degenerate
$m$-simplex $V \in N_m \Comp$ is in $C^\ell$ if and only if
$v^i_i=\{0, \dots, \hat{k}, \dots, m\}$ for all $0\leq i \leq
\ell-1$ and $v^i_i=\{0, \dots, m\}$ for all $\ell \leq i \leq m$.

\begin{prop} \label{compgluing}
Let $V \in C^\ell$. Then the $j$-th face of $V$ is shared with
some other $V' \in C^\ell$ if and only if $j \neq 0, \ell -1, \ell$.
\end{prop}
\begin{pf}
By Remark \ref{gluingofsimplices} we know exactly which other
non-degenerate $m$-simplex $V'$ shares the $j$-th face of $V$. So,
for each $\ell$ and $j$ we only need to check whether or not $V'$ is
in $C^\ell$. Let $V \in C^\ell$.

{\bf Cases} $1 \leq \ell \leq m$ and $j=0$. \\
For all $U \in C^\ell$, we have $U_0=(\{0, \dots, \hat{k}, \dots,
m\})=V_0,$ so we conclude from the description of $V'$ in Remark
\ref{gluingofsimplices} \ref{gluingofsimplicesiv} that $V'$ is not
in $C^\ell$.

{\bf Case} $\ell=m$ and $j=m-1$. \\
In this case, $v^{m-1}_{m-1}=\{0, \dots, \hat{k}, \dots, m\}$ and
$v^m_m=\{0,1,\dots,m\}$. By Remark \ref{gluingofsimplices}
\ref{gluingofsimplicesiv}, the $m-1$st-face of $V$ is shared with
the $V'$ which agrees with $V$ everywhere except in $V_{m-1}$, where
we have $(v')^{m-1}_{m-1}=\{0,\dots,m\}$ instead of
$v^{m-1}_{m-1}=\{0,\dots,\hat{k},\dots,m\}$. But this $V'$ is not an
element of $C^m$.

{\bf Case} $\ell=m$ and $j=m$. \\
In this case, $v^{m-1}_{m-1}=\{0, \dots, \hat{k}, \dots,
m\}\neq\{0,1, \dots, m\}$, so we are in the situation of Remark
\ref{gluingofsimplices} \ref{gluingofsimplicesii}. The $m$-th face
$(V_0, \dots, V_{m-1})$ does not lie in any other non-degenerate
$m$-simplex $V'$, let alone in a $V'$ in $C^m$.

{\bf Case} $\ell=m$ and $j\neq0,m-1,m$. \\
By Remark \ref{gluingofsimplices} \ref{gluingofsimplicesiv}, the
$j$-th face is shared with the $V'$ that agrees with $V$ in $V_0$,
$V_{m-1}$, and $V_m$, so that $V' \in C^m$.

At this point we conclude from the above cases that if $\ell=m$, the
$j$-th face of $V\in C^m$ is shared with another $V'\in C^m$ if and
only if $j\neq 0,m-1,m$.

{\bf Cases} $1 \leq \ell \leq m-1$ and $j=\ell-1$. \\
The $\ell-1$st face of $V$ is shared with that $V'$ which agrees
with $V$ everywhere except in $V_{\ell-1}$, where we have
$(v')^{\ell-1}_{\ell-1}=\{0,\dots,m\}$ instead of
$v^{\ell-1}_{\ell-1}=\{0, \dots, \hat{k},\dots,m\}$. Hence $V'$ is
not in $C^\ell$.

{\bf Cases} $1 \leq \ell \leq m-1$ and $j=\ell$. \\
Similarly, the $\ell$-th face of $V$ is shared with that $V'$ which
agrees with $V$ everywhere except in $V_\ell$, where we have
$(v')^{\ell}_{\ell}=\{0,\dots,\hat{k},\dots,m\}$ instead of
$v^{\ell}_{\ell}=\{0,\dots,m\}$. Hence $V'$ is not in $C^\ell$.

{\bf Cases} $1 \leq \ell \leq m-1$ and $j\neq0,\ell-1,\ell$. \\
Then the $j$-th face is shared with a $V'$ that agrees with $V$ in
$V_0$, $V_{\ell-1}$, and $V_\ell$, so that $V' \in C^\ell$.

We conclude that the $j$-th face of $V\in C^\ell$ is shared with
some other $V'\in C^\ell$ if and only if $j\neq 0,\ell-1,\ell$.
\end{pf}

\begin{prop} \label{upcloseddecomposition}
Let $0 \leq k \leq m$. Recall that $\Comp$, $\Cen$, and $\Out$ denote the
up-closure in $\bfP\Sd\Delta[m]$ of $(\{0,1, \dots,
\hat{k},\dots,m\})$, $(\{0,1,\dots,m\})$, and $\bfP \Sd
\Lambda^k[m]$ respectively.

Then the poset $\bfP\Sd\Delta[m]=c\Sd^2\Delta[m]$ is the union of
these three up-closed subposets:
$$\bfP\Sd\Delta[m]=\Comp\cup\Cen\cup\Out.$$
The partial order on $\bfP\Sd\Delta[m]$ is given in (\ref{OrderOnDoubleSubdivision}).
\end{prop}

\section{Deformation Retraction of $|N(\Comp \cup \Cen)|$}\label{retractionsection}

In this section we construct a retraction of $|N(\Comp \cup \Cen)|$
to that part of its boundary which lies in $\Out$. As stated in
Proposition \ref{deformationretract}, each stage of the retraction
is part of a deformation retraction, and is thus a homotopy
equivalence. The retraction is done in such a way that we can adapt
it later to the $n$-fold case. We first treat the retraction of
$|N\Comp|$ in detail.

\begin{prop} \label{deformationretract1}
Let $C^m, C^{m-1}, \dots, C^1$ be the collections of non-degenerate
$m$-simplices of $N\Comp$ defined in Section \ref{barycentric}. Then
there is an $m$ stage retraction of $|N\Comp|$ onto
$|N(\Comp\cap(\Cen\cup\Out))|$ which retracts the individual
simplices of $C^m, C^{m-1}, \dots, C^1$ to subcomplexes of their
boundaries. Further, each retraction of each simplex is part of a
deformation retraction.
\end{prop}
\begin{pf}
As an illustration, we first prove the case $m=1$ and $k=0$. The
poset $\bfP\Sd\Delta[1]$ is
$$\xymatrix{\mathbf{(\{0\})} \ar[r] & (\{0\},\{01\}) & (\{01\}) \ar@{.>}[l] \ar@{.>}[r] &  (\{1\},\{01\})
&  (\{1\}) \ar@{.>}[l]_-f}$$ and $\bfP\Sd\Lambda^0[1]$ consists only
of the object $(\{0\})$. Of the nontrivial morphisms in
$\bfP\Sd\Delta[1]$, the only one in $\Out$ is the solid one on the
far left. The poset $\Cen$ consists of the two middle morphisms,
emanating from $(\{01\})$. The only morphism in $\Comp$ is the one
labelled $f$. The intersection $\Comp\cap(\Cen\cup\Out)$ is the
vertex $(\{1\},\{01\})$, which is the target of $f$.

Clearly, after geometrically realizing, the interval $|f|$ can be
deformation retracted to the vertex $(\{1\},\{01\})$. The case $m=1$
with $k=1$ is exactly the same. In fact, $k$ does not matter, since
the simplices no longer have a direction after geometric
realization.

The case $m=2$ and $k=1$ can be similarly observed in Figure
\ref{subdivisionfigure}.

For general $m \in \mathbb{N}$, we construct a {\it topological}
retraction in $m$ steps, starting with Step 0. In Step 0 we retract
those non-degenerate $m$-simplices of $N_m\Comp$ which have an entire
$m-1$-face on the $k$-th face of $\Delta[m]$, \ie in Step 0 we
retract the elements of $C^m$. Generally, in Step $\ell$ we retract
those non-degenerate $m$-simplices of $N_m\Comp$ which have exactly
$\ell$ vertices on the $k$-th face of $\Delta[m]$, \ie in Step
$\ell$ we retract the elements of $C^{m-\ell}$.

We describe Step $m-\ell$ in detail for $2 \leq \ell \leq m$. We
retract each $V \in C^\ell$ to $$(V_0, \dots,
\hat{V}_{\ell-1},V_\ell, \dots) \cup (V_1, \dots, V_m)$$ in such a
way that for each $j\neq 0, \ell-1, \ell$ the $j$-th  face $$(V_0,
\dots, \hat{V}_j, \dots,V_{\ell-1},V_\ell, \dots)$$ is retracted
{\it within itself} to its subcomplex
$$(V_0, \dots, \hat{V}_j, \dots, \hat{V}_{\ell-1},V_\ell, \dots)
\cup(\hat{V}_0, \dots, \hat{V}_j, \dots, V_{\ell-1},V_\ell,
\dots).$$ We can do this to all $V\in C^\ell$ {\it simultaneously}
because the prescription agrees on the overlaps: $V$ shares the face
$(V_0, \dots, \hat{V}_j, \dots, V_{\ell-1},V_\ell, \dots)$ with only
one other non-degenerate $m$-simplex $V'\in C^\ell$, and $V'$ differs
from $V$ only in $V_j'$ by Proposition \ref{compgluing}.

This procedure is done for Step 0 up to and including Step $m-2$.
After Step $m-2$, the only remaining non-degenerate $m$-simplices in
$N_m\Comp$ are those which have only the first vertex (\ie only
$V_0$) on the $k$-th face of $\Delta[m]$. This is the set $C^1$.

Every $V\in C^1$ has
$$V_0=(\{0, \dots, \hat{k}, \dots, m\})$$
$$V_1=(\{0, \dots, \hat{k}, \dots, m\}, \{0, \dots, m \}),$$
so all $V\in C^1$ intersect in this edge. In Step $m-1$, we retract
each $V\in C^1$ to $(V_1, \dots, V_m)$ in such a way that for $j\neq
0,1$ we retract the $j$-th face $V$ to $(V_1, \dots, \hat{V}_j,
\dots)$, and further we retract the 1-simplex $(V_0,V_1)$ to the
vertex $V_1$. We can do this simultaneously to all $V \in C^1$, as
the procedure agrees in overlaps by Proposition \ref{compgluing},
and the observation about $(V_0,V_1)$ we made above. For each $V\in
C^1$, the 0th face $(V_1, \dots, V_m)$ is also the 0th face of a
non-degenerate $m$-simplex $U$ not in $N_m\Comp$, namely
$$U_0=(\{0, \dots, m\})$$
$$U_j=V_j \text{ for } j \geq 1$$
by Remark \ref{gluingofsimplices} \ref{gluingofsimplicesiv}. The
simplex $U$ is even central. Thus, $(V_1, \dots, V_m)$ is in the
intersection $|N(\Comp\cap(\Cen\cup\Out))|$ and we have succeeded in
retracting $|N\Comp|$ to $|N(\Comp\cap(\Cen\cup\Out))|$ in such a
way that each non-degenerate $m$-simplex is retracted within itself.
Further, each retraction is part of a deformation retraction.
\end{pf}

\begin{prop} \label{deformationretract2}
There is a multi-stage retraction of $|N\Cen|$ onto
$|N(\Cen\cap\Out)|$ which retracts each non-degenerate $m$-simplex to
a subcomplex of its boundary. Further, this retraction is part of a
deformation retraction.
\end{prop}
\begin{pf}
We describe how this works for the case $m=2$ pictured in Figure
\ref{subdivisionfigure}. The poset $\Cen$ consists of all the
central triangles emanating from 012. These have two dotted sides
emanating from 012. The intersection $\Cen\cap\Out$ consists of the
indicated solid lines on those triangles and their vertices (the two
triangles at the bottom have no solid lines). To topologically
deformation retract $|N\Cen|$ onto $|N(\Cen\cap\Out)|$, we first
deformation retract the vertical, downward pointing edge 012 -
02,012 by pulling the vertex 02,012 up to 012 while at the same time
deforming the left bottom triangle to the edge 012 - 0,02,012 and
the right bottom triangle to the edge 012 - 2,02,012.

Then we consecutively deform each of the left triangles emanating
from 012 to the its solid edge and the edge of the next one, holding
the vertex 012 fixed. We deform the left triangles in this manner
all the way until we reach the vertically pointing edge 012 - 1,012.

Similarly, we consecutively deform each of the right triangles
emanating from 012 to the its solid edge and the edge of the next
one, holding the vertex 012 fixed. We deform the right triangles in
this manner all the way until we reach the vertically pointing edge
012 - 1,012.

Finally, we deformation retract the last remaining edge 012 - 1,012
up to the vertex 1, 012, and we are finished.

It is possible to describe this in arbitrary dimensions, although it
gets rather technical, as we already have seen in Proposition
\ref{deformationretract1}.
\end{pf}

\begin{prop} \label{deformationretract}
There is a multi-stage retraction of $|N(\Comp \cup \Cen)|$ to
$|N((\Comp \cup \Cen)\cap\Out)|$ which retracts each non-degenerate
$m$-simplex to a subcomplex of its boundary. Further, each
retraction of each simplex is part of a deformation retraction. See
Figure \ref{subdivisionfigure}.
\end{prop}
\begin{pf}
This follows from Proposition \ref{deformationretract1} and
Proposition \ref{deformationretract2}.
\end{pf}

\section{Nerve, Pushouts, and Colimit Decompositions of Subposets of
$\bfP\Sd\Delta[m]$}\label{pushoutsection}

In this section we prove that the nerve is compatible with certain
colimits and express posets satisfying a chain condition as a
colimit of two finite ordinals, in a way compatible with nerve. The
somewhat technical results of this section are crucial for the
verification of the pushout axiom in the proof of the Thomason model
structure on $\mathbf{Cat}$ and $\mathbf{nFoldCat}$ in Sections
\ref{Thomasonsection} and \ref{Thomasonnfoldsection}. The results of
this section will have $n$-fold versions in Section
\ref{sectionnfolddecompositions}.

We begin by proving that the nerve preserves certain pushouts in
Proposition \ref{nervecommuteswithpushout}. The question of
commutation of nerve with certain pushouts is an old one, and has
been studied in Section 5 of \cite{fritschlatch2}.

The next task is to express posets satisfying a chain condition as a
colimit of two finite ordinals $[m-1]$ and $[m]$ in Proposition
\ref{colimitdecomposition}, and similarly express their nerves as a
colimit of $\Delta[m-1]$ and $\Delta[m]$ in Proposition
\ref{simplicial_colimitdecomposition}. As a consequence, the nerve
functor preserves these colimits in Corollary
\ref{nervecommuteswithcolimitdecomposition}. The combinatorial proof
that our posets of interest, namely $\bfP\Sd \Delta[m]$, $\Cen$,
$\Out$, $\Comp$, $\Comp \cup \Cen$, $\bfP\Sd\Lambda^k[m]$, and $\Out
\cap (\Comp \cup \Cen)$, satisfy the chain conditions, is found in
Remark \ref{remarkonpaths} and Proposition \ref{satisfyhypothesis}.
Corollary \ref{cor:specific_colimit_decompositions} summarizes the
nerve commutation for the decompositions of the posets of interest.

\begin{prop} \label{nervecommuteswithpushout}
Suppose $\bfQ$, $\bfR$, and $\bfS$ are categories, and $\bfS$ is a full subcategory of $\bfQ$ and $\bfR$ such that
\begin{enumerate}
\item \label{nervecommuteswithpushouti}
If $\xymatrix@1{f:x \ar[r] & y }$ is a morphism in $\bfQ$ and $x \in
\bfS$, then $y \in \bfS$,
\item \label{nervecommuteswithpushoutii}
If $\xymatrix@1{f:x \ar[r] & y }$ is a morphism in $\bfR$ and $x \in
\bfS$, then $y \in \bfS$.
\end{enumerate}
Then the nerve of the pushout is the pushout of the nerves.
\begin{equation}
N(\bfQ \coprod_\bfS \bfR) \cong N\bfQ \coprod_{N\bfS} N\bfR
\end{equation}
\end{prop}
\begin{pf}
First we claim that there are no free composites in $\bfQ \coprod_\bfS \bfR$.
Suppose $f$ is a morphism in $\bfQ$ and $g$
is a morphism in $\bfR$ and that these are composable in the pushout $\bfQ \coprod_\bfS \bfR$.
$$\xymatrix{w \ar[r]^f & x \ar[r]^g & y }$$
Then $x \in \Obj \bfQ \cap \Obj \bfR = \bfS$, so $y \in \bfS$ by
hypothesis \ref{nervecommuteswithpushoutii}. Since $\bfS$ is full,
$g$ is a morphism of $\bfS$. Then $g \circ f$ is a morphism in
$\bfQ$ and is not free. The other case $f$ in $\bfR$ and $g$ in
$\bfQ$ is exactly the same. Thus the pushout $\bfQ \coprod_\bfS
\bfR$ has no free composites.

Let $(f_1, \dots, f_p)$ be a $p$-simplex in $N(\bfQ \coprod_\bfS
\bfR)$. Then each $f_j$ is a morphism in $\bfQ$ or $\bfR$, as there
are no free composites. Further, by repeated application of the
argument above, if $f_1$ is in $\bfQ$ then every $f_j$ is in $\bfQ$.
Similarly, if $f_1$ is in $\bfR$ then every $f_j$ is in $\bfR$. Thus
we have a morphism $\xymatrix{N(\bfQ \coprod_\bfS \bfR) \ar[r] &
N\bfQ \coprod_{N\bfS} N\bfR}$. Its inverse is the canonical morphism
$\xymatrix{N\bfQ \coprod_{N\bfS} N\bfR \ar[r] & N(\bfQ \coprod_\bfS
\bfR) }$.
\end{pf}

\begin{prop} \label{nervecommuteswithpushouthypothesisverification}
The full subcategory $(\Comp \cup \Cen)\cap\Out$ of the categories
$\Comp \cup \Cen$ and $\Out$ satisfies
\ref{nervecommuteswithpushouti} and \ref{nervecommuteswithpushoutii}
of Proposition \ref{nervecommuteswithpushout}.
\end{prop}
\begin{pf}
Since $\Comp$ and $\Cen$ are up-closed, the union $\Comp \cup \Cen$
is up-closed, as is its intersection with up-closed poset $\Out$.
Hence conditions \ref{nervecommuteswithpushouti} and
\ref{nervecommuteswithpushoutii} of Proposition
\ref{nervecommuteswithpushout} follow.
\end{pf}

\begin{prop} \label{colimitdecomposition}
Let $\bfT$ be a poset and $m \geq 1$ a positive integer such that
the following hold.
\begin{enumerate}
\item \label{colimitdecompositioni}
Any linearly ordered subposet $U=\{U_0 < U_1 < \cdots < U_p \}$ of
$\bfT$ with $|U|\leq m+1$ is contained in a linearly ordered
subposet $V$ of $\bfT$ with $m+1$ distinct elements.
\item \label{colimitdecompositionii}
Suppose $x$ and $y$ are in $\bfT$ and $x \leq y$. If $V$ and $V'$
are linearly ordered subposets of $\mathbf{T}$ with exactly $m+1$
elements, and both $V$ and $V'$ contain $x$ and $y$, then there
exist linearly ordered subposets $W^0,W^1, \dots, W^k$ of $\bfT$
such that
\begin{enumerate}
\item
$W^0=V$
\item
$W^k=V'$
\item
For all $0 \leq j \leq k$, the linearly ordered poset $W^j$ has
exactly $m+1$ elements
\item
For all $0 \leq j \leq k$, we have $x \in W^j$ and $y \in W^j$
\item
For all $0 \leq j \leq k-1$, the poset $W^j \cap W^{j+1}$ has $m$
elements.
\end{enumerate}
\item \label{colimitdecompositioniii}
If $m=1$, we further assume that there are no linearly ordered
subposets with 3 or more elements, that is, there are no nontrivial
composites $x < y < z$. Whenever $m=1$, hypothesis
\ref{colimitdecompositionii} is vacuous.
\end{enumerate}
Let $\bfJ$ denote the poset of linearly ordered subposets $U$ of
$\bfT$ with exactly $m$ or $m+1$ elements. Then $\bfT$ is the
colimit of the functor
$$\xymatrix{F:\bfJ \ar[r] & \mathbf{Cat}}$$
$$\xymatrix{U \ar@{|->}[r] & U.}$$
The components of the universal cocone $\xymatrix@1{\pi:F
\ar@{=>}[r] & \Delta_\bfT}$ are the inclusions $\xymatrix@1{F(U)
\ar[r] & \bfT}$.
\end{prop}
\begin{pf}
Suppose $\bfS\in\mathbf{Cat}$ and $\xymatrix@1{\alpha:F \ar@{=>}[r]
& \Delta_\bfS}$ is a natural transformation. We define a functor
$\xymatrix@1{G:\bfT \ar[r] & \bfS}$ as follows. Let $x$ and $y$ be
elements of $\bfT$ and suppose $x \leq y$. By hypothesis
\ref{colimitdecompositioni}, there is a linearly ordered subposet
$V$ of $\bfT$ which contains $x$ and $y$ and has exactly $m+1$
elements. We define $G(x\leq y):=\alpha_V(x\leq y)$.

We claim $G$ is well defined. If $V'$ is another linearly ordered
subposet of $\mathbf{T}$ which contains $x$ and $y$ and has exactly
$m+1$ elements, then we have a sequence $W^0, \dots, W^k$ as in
hypothesis \ref{colimitdecompositionii}, and the naturality diagrams
below.
$$\xymatrix@C=5pc{W^i \ar[r]^-{\alpha_{W^i}} & \bfS \ar@{=}[d] \\
W^i \cap W^{i+1} \ar[r]^-{\alpha_{W^i \cap W^{i+1}}} \ar[u] \ar[d] & \bfS \ar@{=}[d] \\
W^{i+1} \ar[r]_-{\alpha_{W^{i+1}}} & \bfS }$$ Thus we have a string
of equalities
$$\alpha_{W^0}(x\leq y)=\alpha_{W^1}(x\leq y)=\cdots =\alpha_{W^k}(x\leq y),$$
and we conclude $\alpha_V(x\leq y)=\alpha_{V'}(x\leq y)$ so that
$G(x\leq y)$ is well defined.

The assignment $G$ is a functor, as follows. It preserves identities
because each $\alpha_V$ does. If $m=1$, then there are no nontrivial
composites by hypothesis \ref{colimitdecompositioniii}, so $G$
vacuously preserves all compositions. If $m \geq 2$, and the
elements $x < y < z$ in are in $\mathbf{T}$, then there exists a $V$
containing all three of $x$, $y$, and $z$. The functor $\alpha_V$
preserves this composition, so $G$ does also.

By construction, for each linearly ordered subposet $V$ of $\bfT$
with $m+1$ elements we have $\alpha_V=G\circ \pi_V$. Further, $G$ is
the unique such functor, since such posets $V$ cover $\bfT$ by
hypothesis \ref{colimitdecompositioni}.

Lastly we claim that $\alpha_U=G \circ \pi_U$ for any linearly
ordered subposet $U$ of $\bfT$ with $m$ elements. By hypothesis
\ref{colimitdecompositioni} there exists a linearly ordered subposet
$V$ of $\bfT$ with $m+1$ elements such that $U\subseteq V$. If $i$
denotes the inclusion of $U$ into $V$, by naturality of $\alpha$ and
$\pi$ we have
$$\alpha_U=\alpha_V \circ i = G\circ \pi_V \circ i =G \circ \pi_U.$$
\end{pf}

\begin{prop} \label{simplicial_colimitdecomposition}
Let $\bfT$ be a poset and $m \geq 1$ a positive integer such that
the following hold.
\begin{enumerate}
\item \label{simplicial_colimitdecompositioni}
Any linearly ordered subposet $U=\{U_0 < U_1 < \cdots < U_p \}$ of
$\bfT$ is contained in a linearly ordered subposet $V$ of $\bfT$
with $m+1$ distinct elements, in particular, any linearly ordered
subposet of $\bfT$ has at most $m+1$ elements.
\item \label{simplicial_colimitdecompositionii}
Suppose $x_0<x_1<\cdots<x_{\ell}$ are in $\bfT$ and $\ell \leq m$.
If $V$ and $V'$ are linearly ordered subposets of $\mathbf{T}$ with
exactly $m+1$ elements, and both $V$ and $V'$ contain
$x_0<x_1<\cdots<x_{\ell}$, then there exist linearly ordered
subposets $W^0,W^1, \dots, W^k$ of $\bfT$ such that
\begin{enumerate}
\item
$W^0=V$
\item
$W^k=V'$
\item
For all $0 \leq j \leq k$, the linearly ordered poset $W^j$ has
exactly $m+1$ elements
\item
For all $0 \leq j \leq k$, the elements $x_0<x_1<\cdots<x_{\ell}$
are all in $W^j$
\item
For all $0 \leq j \leq k-1$, the poset $W^j \cap W^{j+1}$ has
exactly $m$ distinct elements.
\end{enumerate}
\end{enumerate}
As in Proposition \ref{colimitdecomposition}, let $\bfJ$ denote the
poset of linearly ordered subposets $U$ of $\bfT$ with exactly $m$
or $m+1$ elements, let $F$ be the functor
$$\xymatrix{F:\bfJ \ar[r] & \mathbf{Cat}}$$
$$\xymatrix{U \ar@{|->}[r] & U,}$$
and $\pi$ the universal cocone $\xymatrix@1{\pi:F \ar@{=>}[r] &
\Delta_\bfT}.$  The components of $\pi$ are the inclusions
$\xymatrix@1{F(U) \ar[r] & \bfT}$. Then $N\bfT$ is the colimit of
the functor
$$\xymatrix{NF:\bfJ \ar[r] & \mathbf{SSet}}$$
$$\xymatrix{U \ar@{|->}[r] & NFU}$$
and $\xymatrix@1{N\pi:NF \ar@{=>}[r] & \Delta_{N\bfT}}$ is its
universal cocone.
\end{prop}
\begin{pf}
The principle of the proof is similar to the direct proof of
Proposition \ref{colimitdecomposition}. Suppose $S\in\mathbf{SSet}$
and $\xymatrix@1{\alpha:NF \ar@{=>}[r] & \Delta_S}$ is a natural
transformation. We induce a morphism of simplicial sets
$\xymatrix@1{G:N\bfT \ar[r] & S}$ by defining $G$ on the
$m$-skeleton as follows.

Let $\Delta_{m}$ denote the full subcategory of $\Delta$ on the
objects $[0], [1], \dots, [m]$ and let $\xymatrix@1{\text{tr}_{m}
\co \mathbf{SSet} \ar[r] & \mathbf{Set}^{\Delta_{m}^\text{op}}}$
denote the $m$-th truncation functor. The truncation $\text{tr}_{m}
N\bfT$ is a union of the truncated simplicial subsets
$\text{tr}_{m}NV$ for $V \in \bfJ$ with $\vert V \vert=m+1$, since
$\bfT$ is a union of such $V$. We define
$$\xymatrix{G_{m}\vert_{\text{tr}_{m}NV} \co \text{tr}_{m}NV \ar[r] & \text{tr}_{m}S}$$
simply as $\text{tr}_{m}\alpha_{V}$.

The morphism $G_m$ is well-defined, for if $0\leq\ell \leq m$ and $x
\in (\text{tr}_{m}NV)_\ell$ and $x \in (\text{tr}_{m}NV')_\ell$ with
$\vert V \vert=m+1=\vert V \vert$, then $V$ and $V'$ can be
connected by a sequence $W^0,W^1, \dots, W^k$ of $(m+1)$-element
linearly ordered subsets of $\bfT$ that all contain the linearly
ordered subposet $x$ and satisfy the properties in hypothesis
\ref{simplicial_colimitdecompositionii}. By a naturality argument as
in the proof of Proposition \ref{colimitdecomposition}, we have have
a string of equalities
$$\alpha_{W^0}(x)=\alpha_{W^1}(x)=\cdots =\alpha_{W^k}(x),$$
and we conclude $\alpha_{V}(x)=\alpha_{V'}(x)$ so that $G_m(x)$ is
well defined.

By definition $\Delta_{G_m} \circ \text{tr}_m N\pi=\text{tr}_m
\alpha$. We may extend this to non-truncated simplicial sets using
the following observation: if $\bfC$ is a category in which
composable chains of morphisms have at most $m$-morphisms, and
$\text{sk}_m$ is the left adjoint to $\text{tr}_m$, then the counit
inclusion $$\xymatrix@1{\text{sk}_m\text{tr}_m(N\bfC) \ar[r] & N
\bfC}$$ is the identity. Thus $G_m$ extends to $\xymatrix@1{G\co
N\bfT \ar[r] & S}$ and $\Delta_{G} \circ N\pi=\alpha$.

Lastly, the morphism $G$ is unique, since the simplicial subsets
$NV$ for $\vert V \vert=m+1$ cover $N\bfT$ by hypothesis
\ref{simplicial_colimitdecompositioni}.
\end{pf}

\begin{cor} \label{nervecommuteswithcolimitdecomposition}
Under the hypotheses of Proposition
\ref{simplicial_colimitdecomposition}, the nerve functor commutes
with the colimit of $F$.
\end{cor}

Since $\Sd^2 \Delta[m]$ geometrically realizes to a {\it connected}
simplicial complex that is a union of non-degenerate $m$-simplices,
it is clear that we can move from any non-degenerate $m$-simplex $V$
of $\Sd^2 \Delta[m]$ to any other $V'$ by a chain of non-degenerate
$m$-simplices in which consecutive ones share an $(m-1)$-subsimplex.
However, if $x$ and $y$ are two vertices contained in both $V$ and
$V'$, it is not clear that a chain can be chosen from $V$ to $V'$ in
which all non-degenerate $m$-simplices contain both $x$ and $y$. The
following extended remark explains how to choose such a chain.

\begin{rmk} \label{remarkonpaths}
Our next task is to prepare for the proof of Proposition
\ref{satisfyhypothesis}, which says that the posets $\bfP\Sd
\Delta[m]$, $\Cen$, $\Out$, $\Comp$, and $\Comp \cup \Cen$ satisfy
the hypotheses of Proposition \ref{simplicial_colimitdecomposition}
for $m$, and the posets $\bfP\Sd\Lambda^k[m]$ and $\Out \cap (\Comp
\cup \Cen)$ satisfy the hypotheses of Proposition
\ref{simplicial_colimitdecomposition} for $m-1$. Building on Remark
\ref{gluingofsimplices}, we describe a way of moving from a
non-degenerate $m$-simplex $V$ of $\Sd^2 \Delta[m]$ to another
non-degenerate $m$-simplex $V'$ of $\Sd^2 \Delta[m]$ via a chain of
non-degenerate $m$-simplices, in which consecutive $m$-simplices
overlap in an $(m-1)$-simplex, and each non-degenerate $m$-simplex
in the chain contains specified vertices $x_0<x_1<\cdots<x_{\ell}$
contained in both $V$ and $V'$. Observe that the respective elements
$x_0, x_1, \dots, x_\ell$ are in the same respective positions in
$V$ and $V'$, for if they were in different respective positions, we
would arrive a linearly ordered subposet of length greater than
$m+1$, a contradiction.

We first prove the analogous statement about moving from $V$ to $V'$
for $\Sd\Delta[m]$. The non-degenerate $m$-simplices of
$\Sd\Delta[m]$ are in bijective correspondence with the permutations
of $\{0,1, \dots, m\}$. Namely, the simplex $v=(v_0, \dots, v_m)$
corresponds to $a_0, \dots, a_m$ where $a_i=v_i \backslash v_{i-1}$.
For example, $(\{1\},\{1,2\}, \{0,1,2\})$ corresponds to $1,2,0$.
Swapping $a_i$ and $a_{i+1}$ gives rise to a non-degenerate
$m$-simplex $w$ which shares an $(m-1)$-subsimplex with $v$, that
is, $v$ and $w$ differ only in the $i$-th spot: $v_i \neq w_i$.
Since transpositions generate the symmetric group, we can move from
any non-degenerate $m$-simplex of $\Sd\Delta[m]$ to any other by a
sequence of moves in which we only change one vertex at a time.
Suppose $v$ and $v'$ are the same at spots $s_0 < s_1 < \cdots<
s_\ell$, that is $v_{s_i}=v_{s_i}'$ for $0\leq i \leq \ell$. Then,
using transpositions, we can traverse from $v$ to $v'$ through a
chain $w^1, \dots, w^k$ of non-degenerate $m$-simplices of
$\Sd\Delta[m]$, each of which is equal to $v_{s_1}, v_{s_2}, \dots,
v_{s_\ell}$ in spots $s_1, s_2, \dots, s_\ell$. Indeed, this
corresponds to the embedding of symmetric groups
\begin{equation*}
\Sym(v_{s_1})\times \left( \prod_{i=2}^\ell \Sym(v_{s_i} \backslash
v_{s_{i-1}})  \right) \times \Sym(\{0,\dots,n\}\backslash
v_{s_\ell}) \xymatrix@1{ \ar[r] & } \Sym(\{0,\dots,n\})
\end{equation*} and
generation by the relevant transpositions.

Similar, but more involved, arguments allow us to navigate the
non-degenerate $m$-simplices of $\Sd^2\Delta[m]$. For a {\it fixed}
non-degenerate $m$-simplex $V_m=(v_0^{m},\dots,v_m^m)$ of
$\Sd\Delta[m]$, the non-degenerate $m$-simplices $V=(V_0, \dots,
V_m)$ of $\Sd^2\Delta[m]$ ending in the {\it fixed} $V_m$ correspond
to permutations $A_0, \dots, A_m$ of the vertices of $V_m$. For
example, the 2-simplex in \eqref{Vexample} corresponds to the
permutation
$$\{01\},\;\; \{0\},\;\; \{012\}.$$
Again, arguing by transpositions, we can move from any
non-degenerate $m$-simplex of $\Sd^2\Delta[m]$ ending in $V_m$ to
any other ending in $V_m$ by a sequence of moves in which we only
change one vertex at a time, and at every step, we preserve the
specified vertices $x_0<x_1<\cdots<x_{\ell}$. Holding $V_m$ fixed
corresponds to moving (in $\Sd^2\Delta[m]$) within the subdivision
of one of the non-degenerate $m$-simplices of $\Sd\Delta[m]$ (the
subdivision is isomorphic to $\Sd\Delta[m]$, the case treated
above). See for example Figure \ref{subdivisionfigure} for a
convincing picture.

But how do we move between non-degenerate $m$-simplices that do not
agree in the $m$-th spot, in other words, how do we move from
non-degenerate $m$-simplices of one subdivided non-degenerate
$m$-simplex of $\Sd\Delta[m]$ to non-degenerate $m$-simplices in
another subdivided non-degenerate $m$-simplex of $\Sd\Delta[m]$?
First, we say how to move without requiring containment of the
specified vertices $x_0<x_1<\cdots<x_{\ell}$. Note that if $V$ and
$W$ in $\Sd^2\Delta[m]$ only differ in the last spot $m$, then $V_m$
and $W_m$ agree in all but one spot, say $v_i^m \neq w_i^m$, and the
permutations corresponding to $V$ and $W$ are respectively
$$A_0, \dots, A_{m-1},v_i^m$$
$$A_0, \dots, A_{m-1},w_i^m.$$
Given arbitrary non-degenerate $m$-simplices $V$ and $V'$ of
$\Sd^2\Delta[m]$, we construct a chain connecting $V$ and $V'$ as
follows. First we choose a chain of $m$-simplices
$\{\overline{W}^p\}_{p=0}^q$ in $\Sd\Delta[m]$
$$\overline{W}_m^p=(w_0^p, \dots, w_m^p)$$ $0\leq p \leq q$
from $V_m$ to $V_m'$ which corresponds to transpositions. This we
can do by the first paragraph of this Remark. We define an
$m$-simplex $\overline{W}^p$ in $\Sd^2\Delta[m]$ by
$$\overline{W}^p:=( \dots, \overline{W}^p_m \backslash w_{i_p}^p,\overline{W}^p_m)$$ where
$w_{i_p}^p$ is the vertex of $\overline{W}^p_m$ which distinguishes
it from $\overline{W}^{p-1}_m$ for $1 \leq p \leq q$. The last
letter in the permutation corresponding to $\overline{W}^p$ is
$w_{i_p}^p$. The other vertices of $\overline{W}^p$ indicated by
$\dots$ are any subsimplices of $\overline{W}^p_m$ written in
increasing order. Now, our chain $\{W^j\}_j$ in $\Sd^2\Delta[m]$
from $V$ to $V'$ begins at $V$ and traverses to $\overline{W}^1$:
starting from $V$, we pairwise transpose $v^m_{i_1}$ to the end of
the permutation corresponding to $V$, then we replace $v^m_{i_1}$ by
$w^1_{i_1}$, and then we pairwise transpose the first $m$ letters of
the resulting permutation to arrive at the permutation corresponding
to $\overline{W}^1$. Similarly, starting from $\overline{W}^1$ we
move $w_{i_2}^1$ to the end, replace it by $w_{i_2}^2$, and pairwise
transpose the first $m$ letters to arrive at $\overline{W}^2$.
Continuing in this fashion, we arrive at $V'$ through a chain
$\{W^j\}_j$ of non-degenerate $m$-simplices $W^j$ in
$\Sd^2\Delta[m]$ in which $W^j$ and $W^{j+1}$ share an
$(m-1)$-subsimplex.

Lastly, we must prove that if $V$ and $V'$ both contain specified
vertices $x_0<x_1<\cdots<x_{\ell}$, then the chain $\{W^j\}_j$ of
non-degenerate $m$-simplices can be chosen so that each $W^j$
contains all of the specified vertices $x_0<x_1<\cdots<x_{\ell}$.
Suppose
$$V_{s_i}=x_i=V_{s_i}'$$
for all $0 \leq  i \leq \ell$ and $s_0 < s_1 < \cdots < s_\ell$.
Then $V_m$ and $V_m'$ both contain all of the vertices of $x_0,x_1,
\dots, x_\ell$ since
$$V_m\supseteq V_{s_\ell}=x_\ell \supseteq x_{\ell-1} \supseteq \cdots \supseteq x_0=V_{s_0}$$
$$V_m' \supseteq V_{s_\ell}'=x_\ell \supseteq x_{\ell-1} \supseteq \cdots \supseteq x_0=V_{s_0}'.$$
We first choose the chain $\{\overline{W}^p_m\}_p$ in $\Sd\Delta[m]$
so that each $\overline{W}^p_m$ contains all of the vertices of
$x_0, x_1, \dots, x_\ell$ (this can be done by the discussion of
$\Sd\Delta[m]$ above). Since we have $\overline{W}^p_m \supseteq
x_\ell$, all $w^p_{i_p}$ must satisfy $i_p > s_\ell$. The first
vertices of the non-degenerate $m$-simplex $\overline{W}^p$ in
$\Sd^2\Delta[m]$ indicated by $\dots$ are chosen so that in spots
$s_0, s_1, \dots, s_\ell$ we have $x_0, x_1, \dots, x_\ell$. For
fixed $W^p_m$ we can transpose as we wish, without perturbing
$x_0,x_1, \dots, x_\ell$ (again by the discussion of $\Sd\Delta[m]$
above, but this time applied to the $\Sd\Delta[m]$ isomorphic to the
collection of $m$-simplices of $\Sd\Delta[m]$ ending in
$\overline{W}^p_m$.) On the other hand, the part of $\{W^j\}_j$ in
which we move $w^{p-1}_{i_p}$ to the right does not perturb any of
$x_0,x_1, \dots, x_\ell$ because $i_p > s_\ell$. Thus, each $W^j$
has $x_0,x_1, \dots, x_\ell$ in spots $s_0,s_1, \dots, s_\ell$
respectively.
\end{rmk}

\begin{prop} \label{satisfyhypothesis}
Let $m \geq 1$ be a positive integer. The posets $\bfP\Sd
\Delta[m]$, $\Cen$, $\Out$, $\Comp$, and $\Comp \cup \Cen$ satisfy
\ref{simplicial_colimitdecompositioni} and
\ref{simplicial_colimitdecompositionii} of Proposition
\ref{simplicial_colimitdecomposition} for $m$. Similarly,
$\bfP\Sd\Lambda^k[m]$ and $\Out \cap (\Comp \cup \Cen)$ satisfy
\ref{colimitdecompositioni} and
\ref{simplicial_colimitdecompositionii} of Proposition
\ref{simplicial_colimitdecomposition} for $m-1$. The hypotheses of
Proposition \ref{simplicial_colimitdecomposition} imply those of
Proposition \ref{colimitdecomposition}, so Proposition
\ref{colimitdecomposition} also applies to these posets.
\end{prop}
\begin{pf}
We first consider $m=1$ and the various subposets of $\bfP\Sd
\Delta[1]$. Let $k=0$ (the case $k=1$ is symmetric). The poset
$\bfP\Sd\Delta[1]$ is
$$\xymatrix{\mathbf{(\{0\})} \ar[r] & (\{0\},\{01\}) & (\{01\}) \ar@{.>}[l] \ar@{.>}[r] &  (\{1\},\{01\})
&  (\{1\}) \ar@{.>}[l]_-f}$$ and $\bfP\Sd\Lambda^0[1]$ consists only
of the object $(\{0\})$ (the typography is chosen to match with
Figure \ref{subdivisionfigure}). Of the nontrivial morphisms in
$\bfP\Sd\Delta[1]$, the only one in $\Out$ is the solid one on the
far left. The poset $\Cen$ consists of the two middle morphisms,
emanating from $(\{01\})$. The only morphism in $\Comp$ is the one
labelled $f$. The union $\Comp \cup \Cen$ consist of all the dotted
arrows and their sources and targets. The intersection $\Out \cap
(\Comp \cup \Cen)$ consists only of the vertex $(\{0\},\{0,1\})$.
The hypotheses \ref{simplicial_colimitdecompositioni} and
\ref{simplicial_colimitdecompositionii} of Proposition
\ref{simplicial_colimitdecomposition} are clearly true by inspection
for $\bfP\Sd \Delta[1]$, $\Cen$, $\Out$, $\Comp$, and $\Comp \cup
\Cen$ and also $\bfP\Sd\Lambda^k[1]$ and $\Out \cap (\Comp \cup
\Cen)$.

We next prove that $\bfP\Sd \Delta[m]$ satisfies hypothesis
\ref{simplicial_colimitdecompositioni} of Proposition
\ref{simplicial_colimitdecomposition} for $m \geq 2$, and also its
various subposets satisfy hypothesis
\ref{simplicial_colimitdecompositioni}. Suppose $U=\{U_0 < U_1 <
\cdots < U_p\}$ is a linearly ordered subposet of $\bfP\Sd
\Delta[m]$. As before, we write $U_i=(u_0^i, \dots, u^i_{r_i})$. We
extend $U$ to a linearly ordered subposet $V$ with $m+1$ elements so
that $U_i$ occupies the $r_i$-th place (the lowest element is in the
0-th place). For $j \leq r_0$, let $V_j=(u_0^0,\dots,u_j^0)$. For
$j=r_i$, $V_j:=U_i$. For $r_i \leq j < r_{i+1}-1$, we define
$V_{j+1}$ as $V_j$ with one additional element of $U_{i+1}
\backslash U_i$. If $|U_p|=m+1$, then we are now finished. If
$|U_p|=r_p+1< m+1$, then extend $U_p$ to a strictly increasing chain
of subsets of $\{0,\dots,m\}$ of length $m+1$, where the new subsets
are $v_1,\dots,v_{m+1-(r_p+1)}$ and define for $j=1, \dots, m-r_p$
$$V_{r_p+j}:=V_{r_p} \cup \{v_1, \dots, v_j\}.$$ Then we have $U$
contained in $V=\{V_0< \dots < V_m\}$.

Easy adjustments show that the poset $\Cen$ satisfies hypothesis
\ref{simplicial_colimitdecompositioni} for $m\geq2$. If $U$ is a
linearly ordered subposet of $\Cen$, then each $u^i_{r_i}$ is
$\{0,1,\dots,m\}$ by Proposition \ref{center}. We take
$V_0=(\{0,1,\dots,m\})$ and then successively throw in $u_0^0,
\dots, u^0_{r_0-1}$ to obtain $V_1,\dots,V_{r_0}$. The higher
$V_j$'s are as above. By Proposition \ref{center}, the extension $V$
lies in $\Cen$. A similar argument works for $\Comp$, since it is
also the up-closure of a single point, namely $(\{0,1, \dots,
\hat{k}, \dots, m\})$. The union $\Comp \cup \Cen$ also satisfies
hypothesis \ref{simplicial_colimitdecompositioni} for $m\geq 2$: if
$U$ is a subposet of the union, then $U_0$ is in at least one of
$\Comp$ or $\Cen$, and all the other $U_i$'s are also contained in
that one, so the proof for $\Comp$ or $\Cen$ then finishes the job.

The poset $\Out$ satisfies hypothesis
\ref{simplicial_colimitdecompositioni} for $m\geq 2$, for if $U$ is
a subposet of $\Out$, then $U_0$ must contain some $u^0_i$ in
$\Lambda^k[m]$ by Proposition \ref{outer}. We extend to the left of
$U_0$ by taking $V_0=(u^0_i)$ and then successively throwing in the
remaining elements of $U_0$. The rest of the extension proceeds as
above, since everything above $U_0$ also contains $u^0_i \in
\Lambda^k[m]$. The poset $\Out \cap \Comp$ satisfies hypothesis
\ref{simplicial_colimitdecompositioni} for $m-1$ rather than $m$
because any element in the intersection must have at least 2
vertices, namely a vertex in $\Lambda^k[m]$ and $\{0, \dots,
\hat{k}, \dots, m\}$. Similarly, the poset $\Out \cap \Cen$
satisfies hypothesis \ref{simplicial_colimitdecompositioni} for
$m-1$ rather than $m$ because any element in the intersection must
have at least 2 vertices, namely a vertex in $\Lambda^k[m]$ and
$\{0, \dots, m\}$. The proofs that $\Out \cap \Comp$ and $\Out \cap
\Cen$ satisfy hypothesis \ref{simplicial_colimitdecompositioni} are
similar to the above. Since unions of subposets of
$\bfP\Sd\Delta[m]$ that satisfy hypothesis
\ref{simplicial_colimitdecompositioni} for $m-1$ also satisfy
hypothesis \ref{simplicial_colimitdecompositioni} for $m-1$, we see
that
\begin{equation}\label{equ:union_intersection}
(\Out \cap \Comp)\cup(\Out \cap \Cen)=\Out \cap (\Comp\cup\Cen )
\end{equation}
also satisfies hypothesis \ref{simplicial_colimitdecompositioni} for
$m-1$.

Lastly $\bfP\Sd\Lambda^k[m]$ satisfies hypothesis
\ref{simplicial_colimitdecompositioni} for $m-1$. It is down closed
by Proposition \ref{Lambdadownclosed}, so for a subposet $U$, the
extension of $U$ to the left in $\bfP\Sd \Delta[m]$ described above
is also in $\bfP\Sd\Lambda^k[m]$. Any extension to the right which
includes $k$ in the final $m$-element set is also in
$\bfP\Sd\Lambda^k[m]$ by the discussion after equation
\eqref{simplicesofhorn}.

Next we turn to hypothesis \ref{simplicial_colimitdecompositionii}
of Proposition \ref{simplicial_colimitdecomposition} for the
subposets of $\bfP\Sd \Delta[m]$ in question, where $m\geq 2$. The
poset $\bfP\Sd \Delta[m]$ satisfies hypothesis
\ref{simplicial_colimitdecompositionii} by Remark
\ref{remarkonpaths}.

The poset $\Cen$ is the up-closure of $(\{0,1, \dots, m\})$ in
$\bfP\Sd \Delta[m]$. Every linearly ordered subposet of $\Cen$ with
$m+1$ elements must begin with $(\{0,1, \dots, m\})$.  Given
$(m+1)$-element, linearly ordered subposets $V$ and $V'$ of $\Cen$
with specified elements $x_0<x_1<\cdots<x_{\ell}$ in common, we can
select the chain $\{W^j\}_j$ in Remark \ref{remarkonpaths} so that
each $W^j$ has $(\{0,1, \dots, m\})$ as its $0$-vertex. Thus $\Cen$
satisfies hypothesis \ref{simplicial_colimitdecompositionii}. The
poset $\Comp$ similarly satisfies hypothesis
\ref{simplicial_colimitdecompositionii}, as it is also the
up-closure of an element in $\bfP\Sd \Delta[m]$.

The union $\Comp \cup \Cen$ satisfies hypothesis
\ref{simplicial_colimitdecompositionii} as follows. If $V$ and $V'$
(of cardinality $m+1$) are both linearly ordered subposets of
$\Comp$ or are both linearly ordered subposets of $\Cen$
respectively with the specified elements in common, then we may
simply take the chain in $\Comp$ or $\Cen$ respectively. If $V$ is
in $\Cen$ and $V'$ is in $\Comp$, then $V_0=(\{0,1, \dots, m\})$ and
$V_0'=(\{0, \dots, \hat{k},\dots, m\})$. Suppose
$$V_{s_i}=x_i=V_{s_i}'$$
for all $0 \leq  i \leq \ell$ and $s_0 < s_1 < \cdots < s_\ell$.
Then $x_0$ contains both $\{0,1, \dots, m\}$ and $\{0, \dots,
\hat{k},\dots, m\}$. Then we move from $V'$ to $V''$ by transposing
$\{0,1, \dots, m\}$ down to vertex 0, leaving everything else
unchanged. This chain from $V'$ to $V''$ is in $\Comp$ until it
finally reaches $V''$, which is in $\Cen$. From $V$ we can reach
$V''$ via a chain in $\Cen$ as above. Putting these two chains
together, we move from $V$ to $V'$ as desired.

To show $\Out$ satisfies hypothesis
\ref{simplicial_colimitdecompositionii}, suppose $V$ and $V'$ are
linearly ordered subposets of cardinality $m+1$ with
$V_{s_i}=x_i=V_{s_i}'$ for all $0 \leq  i \leq \ell$ and $s_0 < s_1
< \cdots < s_\ell$. If $V_0=V_0'$, then we can make certain that the
chain $\{W^j\}_j$ in Remark \ref{remarkonpaths} satisfies
$W^j_0=V_0=V_0' \in \bfP\Sd\Lambda^k[m]$. Then each $W^j$ lies in
$\Out$, and we are finished. If $V_0\neq V_0'$, then we move from
$V'$ to $V''$ with $V_0''=V_0$ as follows. The elements $V_0$ and
$V_0'$ are both in $V_{s_0}=x_0=V_{s_0}'$, so we can transpose $V_0$
in $V'$ down to the 0-vertex and interchange $V_0$ and $V_0'$. Each
step of the way is in $\Out$. The result is $V''$, to which we can
move from $V$ on a chain in $\Out$.

We claim that the subposet $\Out \cap \Comp$ of $\bfP\Sd\Delta[m]$
satisfies hypothesis \ref{simplicial_colimitdecompositionii} for
$m-1$. Suppose  $V$ and $V'$ are linearly ordered subposets of
cardinality $m$ with $V_{s_i}=x_i=V_{s_i}'$ for all $0 \leq  i \leq
\ell$ and $s_0 < s_1 < \cdots < s_\ell$, where $1\leq \ell \leq
m-1$. Then $V_0=(v , \{0, \dots, \hat{k}, \dots, m\})$ and
$V_0'=(v', \{0, \dots, \hat{k}, \dots, m\})$ where $v$ and $v'$ are
elements of $\bfP\Sd\Lambda^k[m]$. We extend the $m$-element
linearly ordered posets $V$ and $V'$ to $(m+1)$-element linearly
ordered posets $\bar{V}$ and $\bar{V}'$ in $\Comp$ by putting $(\{0,
\dots, \hat{k}, \dots, m\})$ in the 0-th spot of $\bar{V}$ and
$\bar{V}'$. If $v=v'$, then we can find a chain $\{W^j\}_j$ from
$\bar{V}$ to $\bar{V}'$ in $\Comp$ which preserves $x_0,x_1, \dots,
x_\ell,$ and $v$ using the above result that $\Comp$ satisfies
\ref{simplicial_colimitdecompositionii} for $m$. Truncating the 0-th
spot of each $W^j$, we obtain the desired chain in $\Out \cap
\Comp$. If $v \neq v'$, then we find a chain in $\Comp$ from
$\bar{V}'$ to a $\bar{V}''$ with $v''=v$, like above, and then find
a chain in $\Comp$ from $\bar{V}$ to $\bar{V}''$. Combining chains,
and truncating the 0-th spot again gives us the desired path from
$V$ to $V'$.

By a similar argument, with the role of $\{0, \dots, \hat{k}, \dots,
m\}$ played by $\{0,1, \dots, m\}$, the poset $\Out \cap \Cen$
satisfies hypothesis \ref{simplicial_colimitdecompositionii} for
$m-1$. Next we claim that the union of $\Out \cap \Comp$ with $\Out
\cap \Cen$ also satisfies hypothesis
\ref{simplicial_colimitdecompositionii} for $m-1$. Suppose $V
\subseteq \Out \cap \Comp$ and $V' \subseteq \Out \cap \Cen$ are
$m$-element linearly ordered subposets with $V_{s_i}=x_i=V_{s_i}'$
for all $0 \leq  i \leq \ell$ and $s_0 < s_1 < \cdots < s_\ell$,
where $1\leq \ell \leq m-1$. Then $v$, $v'$, $\{0, \dots, \hat{k},
\dots, m\}$, and $\{0,1, \dots, m\}$ are in $x_0$, so we can
transpose $v$ and $\{0, \dots, \hat{k}, \dots, m\}$ down in $V'$ to
take the place of $v'$ and $\{0,1, \dots, m\}$, without perturbing
$x_0, x_1, \dots, x_\ell$. The resulting poset $V''$ is in $\Out
\cap \Comp$, and was reached from $V'$ by a chain in $\Out \cap
\Cen$. By the above, we can reach $V''$ from $V$ by a chain in $\Out
\cap \Comp$. Thus we have connected $V$ and $V'$ by a chain in
\eqref{equ:union_intersection}, always preserving $x_0, x_1, \dots,
x_\ell$, and therefore $\Out \cap (\Comp \cup \Cen)$ satisfies
hypothesis \ref{simplicial_colimitdecompositionii} for $m-1$.
\end{pf}

\begin{rmk}
The posets $\bfC^\ell$ do not satisfy the hypotheses of Proposition
\ref{simplicial_colimitdecomposition}, nor those of Proposition
\ref{colimitdecomposition}.
\end{rmk}

\begin{cor} \label{cor:specific_colimit_decompositions}
Let $m \geq 1$ be a positive integer.
\begin{enumerate}
\item \label{cor:specific_colimit_decompositionsi}
The posets $\bfP\Sd \Delta[m]$, $\Cen$, $\Out$, $\Comp$, and $\Comp
\cup \Cen$ are each a colimit of finite ordinals $[m-1]$ and $[m]$.
Similarly, the posets $\bfP\Sd\Lambda^k[m]$ and $\Out \cap (\Comp
\cup \Cen)$ are each a colimit of finite ordinals $[m-2]$ and
$[m-1]$. (By definition $[-1]=\emptyset$.)
\item \label{cor:specific_colimit_decompositionsii}
The simplicial sets $N(\bfP\Sd \Delta[m])$, $N(\Cen)$, $N(\Out)$,
$N(\Comp)$, and $N(\Comp \cup \Cen)$ are each a colimit of
simplicial sets of the form $\Delta[m-1]$ and $\Delta[m]$.
Similarly, the simplicial sets $N(\bfP\Sd\Lambda^k[m])$ and $N(\Out
\cap (\Comp \cup \Cen))$ are each a colimit of simplicial sets of
the form $\Delta[m-2]$ and $\Delta[m-1]$. (By definition
$[-1]=\emptyset$.)
\item \label{cor:specific_colimit_decompositionsiii}
The nerve of the colimit decomposition in $\mathbf{Cat}$ in
\ref{cor:specific_colimit_decompositionsi} is the colimit
decomposition in $\mathbf{SSet}$ in
\ref{cor:specific_colimit_decompositionsii}.
\end{enumerate}
\end{cor}
\begin{pf}
\begin{enumerate}
\item
By Proposition \ref{satisfyhypothesis}, the posets $\bfP\Sd
\Delta[m]$, $\Cen$, $\Out$, $\Comp$, and $\Comp \cup \Cen$ satisfy
hypotheses \ref{simplicial_colimitdecompositioni} and
\ref{simplicial_colimitdecompositionii} of Proposition
\ref{simplicial_colimitdecomposition} for $m$, as do the posets
$\bfP\Sd\Lambda^k[m]$ and $\Out \cap (\Comp \cup \Cen)$ for $m-1$.
The hypotheses of Proposition \ref{simplicial_colimitdecomposition}
imply the hypotheses of Proposition \ref{colimitdecomposition}, so
part \ref{cor:specific_colimit_decompositionsi} of the current
corollary follows from Proposition \ref{colimitdecomposition}.
\item
By Proposition \ref{satisfyhypothesis}, the posets $\bfP\Sd
\Delta[m]$, $\Cen$, $\Out$, $\Comp$, and $\Comp \cup \Cen$ satisfy
hypotheses \ref{simplicial_colimitdecompositioni} and
\ref{simplicial_colimitdecompositionii} of Proposition
\ref{simplicial_colimitdecomposition} for $m$, as do the posets
$\bfP\Sd\Lambda^k[m]$ and $\Out \cap (\Comp \cup \Cen)$ for $m-1$.
So Proposition \ref{simplicial_colimitdecomposition} applies and we
immediately obtain part \ref{cor:specific_colimit_decompositionsii}
of the current corollary.
\item
This follows from Corollary
\ref{nervecommuteswithcolimitdecomposition} and Proposition
\ref{satisfyhypothesis}.
\end{enumerate}
\end{pf}

\section{Thomason Structure on {\bf Cat}} \label{Thomasonsection}

The Thomason structure on {\bf Cat} is transferred from the standard
model structure on {\bf SSet} by transferring across the adjunction
\begin{equation}
\xymatrix@C=4pc{\mathbf{SSet} \ar@{}[r]|{\perp}
\ar@/^1pc/[r]^{\Sd^2} &  \ar@/^1pc/[l]^{\Ex^2} \mathbf{SSet}
\ar@{}[r]|{\perp} \ar@/^1pc/[r]^{c} & \ar@/^1pc/[l]^{N}
\mathbf{Cat}}
\end{equation}
as in \cite{thomasonCat}. In other words, a functor $F$ in
$\mathbf{Cat}$ is a weak equivalence or fibration if and only if
$\Ex^2 NF$ is. We present a quick proof that this defines a model
structure using a corollary to Kan's Lemma on Transfer. Although
Thomason did not do it exactly this way, it is practically the same,
in spirit. Our proof relies on the results in the previous sections:
the decomposition of $\Sd^2 \Delta[m]$, the commutation of nerve
with certain colimits, and the deformation retraction.

This proof of the Thomason structure on ${\bf Cat}$ will be the
basis for our proof of the Thomason structure on ${\bf nFoldCat}$.
The key corollary to Kan's Lemma on Transfer is the following
Corollary, inspired by Proposition 3.4.1 in
\cite{worytkiewicz2Cat}.\footnote{The difference between
\cite{worytkiewicz2Cat} and the present paper is that in hypothesis
\ref{KanCorollaryi} we require $Fi$ and $Fj$ to be small with
respect to the entire category $\bfD$, rather than merely small with
respect to $FI$ and $FJ$.}

\begin{cor} \label{KanCorollary}
Let $\mathbf{C}$ be a cofibrantly generated model category with
generating cofibrations $I$ and generating acyclic cofibrations $J$.
Suppose {\bf D} is complete and cocomplete, and that $F \dashv G$ is
an adjunction as in {\rm (\ref{adjunction})}.
\begin{equation} \label{adjunction}
\xymatrix@C=4pc{{\bf C} \ar@{}[r]|{\perp} \ar@/^1pc/[r]^{F} &
\ar@/^1pc/[l]^{G} {\bf D}}
\end{equation}
Assume the following.
\begin{enumerate}
\item \label{KanCorollaryi}
For every $i\in I$ and $j \in J$, the objects $\dom Fi$ and $\dom
Fj$ are small with respect to the entire category $\bfD$.
\item \label{KanCorollaryii}
For any ordinal $\lambda$ and any colimit preserving functor
$\xymatrix@1{X\co\lambda \ar[r] & {\bf C}}$ such that
$\xymatrix@1{X_{\beta} \ar[r] & X_{\beta+1}}$ is a weak equivalence,
the transfinite composition
$$\xymatrix{X_0 \ar[r] &}\underset{\lambda}{\colim}X$$
is a weak equivalence.
\item \label{KanCorollaryiii}
For any ordinal $\lambda$ and any colimit preserving functor
$\xymatrix@1{Y\co\lambda \ar[r] & {\bf D}}$, the functor $G$
preserves the colimit of $Y$.
\item \label{KanCorollaryiv}
If $j'$ is a pushout of $F(j)$ in {\bf D} for $j \in J$, then
$G(j')$ is a weak equivalence in {\bf C}.
\end{enumerate}
Then there exists a cofibrantly generated model structure on {\bf D}
with generating cofibrations $FI$ and generating acyclic
cofibrations $FJ$. Further, $f$ is a weak equivalence in {\bf D} if
and only $G(f)$ is a weak equivalence in {\bf C}, and $f$ is a
fibration in {\bf D} if and only $G(f)$ is a fibration in {\bf C}.
\end{cor}
\begin{pf}
For a proof of a similar statement, see \cite{fiorepaolipronk1}. The
only difference between the statement here and the one proved in
\cite{fiorepaolipronk1} is that here we only require in hypothesis
\ref{KanCorollaryiii} that $G$ preserves colimits indexed by an
ordinal $\lambda$, rather than more general filtered colimits. The
proof of the statement here is the same as in
\cite{fiorepaolipronk1}: it is a straightforward application of
Kan's Lemma on Transfer.
\end{pf}

\begin{lem} \label{ExPreservesAndReflectsWEs}
The functor $\Ex$ preserves and reflects weak equivalences. That is,
a morphism $f$ of simplicial sets is a weak equivalence if and only
if $\Ex f$ is a weak equivalence.
\end{lem}
\begin{pf}
There is a natural weak equivalence $\xymatrix@1{1_{\mathbf{SSet}}
\ar@{=>}[r] & \Ex}$ by Lemma 3.7 of \cite{kancss}, or more recently
Theorem 6.2.4 of \cite{joyaltierneysimplicial}, or Theorem 4.6 of
\cite{goerssjardine}. Then the Proposition follows from the
naturality diagram below.
$$\xymatrix@C=4pc{ X \ar[r]^-{\text{w.e.}} \ar[d]_f & \Ex X \ar[d]^{\Ex f} \\ Y \ar[r]_-{\text{w.e.}} & \Ex Y }$$
\end{pf}

We may now prove Thomason's Theorem.

\begin{thm} \label{CatCase}
There is a model structure on $\mathbf{Cat}$ in which a functor $F$
is a weak equivalence respectively fibration if and only if
$\Ex^2NF$ is a weak equivalence respectively fibration in
$\mathbf{SSet}$. This model structure is cofibrantly generated with
generating cofibrations
$$\{\xymatrix{c\Sd^2\partial\Delta[m] \ar[r] & c\Sd^2\Delta[m]}|\; m \geq 0\}$$
and generating acyclic cofibrations
$$\{\xymatrix{c\Sd^2\Lambda^k[m] \ar[r] & c\Sd^2\Delta[m]}|\; 0 \leq k \leq m \text{ and } m \geq 1\}.$$
These functors were explicitly described in Section
\ref{barycentric}.
\end{thm}
\begin{pf}
\begin{enumerate}
\item \label{CatCasei}
The categories $c\Sd^2\partial\Delta[m]$ and $c\Sd^2\Lambda^k[m]$
each have a finite number of morphisms, hence they are finite, and
are small with respect to $\mathbf{Cat}$. For a proof, see
Proposition 7.6 of \cite{fiorepaolipronk1}.
\item \label{CatCaseii}
The model category $\mathbf{SSet}$ is cofibrantly generated, and the
domains and codomains of the generating cofibrations and generating
acyclic cofibrations are finite. By Corollary 7.4.2 in \cite{hovey},
this implies that transfinite compositions of weak equivalences in
$\mathbf{SSet}$ are weak equivalences.
\item \label{CatCaseiii}
The nerve functor preserves filtered colimits. Every ordinal is
filtered, so the nerve functor preserves $\lambda$-sequences.

The $\Ex$ functor preserves colimits of $\lambda$-sequences as well.
We use the idea in the proof of Theorem 4.5.1 of
\cite{worytkiewicz2Cat}. First recall that for each $m$, the
simplicial set $\Sd\Delta[m]$ is finite, so that $\mathbf{SSet}(\Sd
\Delta[m], -)$ preserve colimits of all $\lambda$-sequences. If
$\xymatrix@1{Y\co\lambda \ar[r] & \mathbf{SSet}}$ is a
$\lambda$-sequence, then
$$\aligned
(\Ex\; \underset{\lambda}{\colim} Y)_m &= \mathbf{SSet}(\Sd
\Delta[m],\underset{\lambda}{\colim} Y) \\ &\cong
\underset{\lambda}{\colim}
\mathbf{SSet}(\Sd \Delta[m], Y)\\
&\cong (\underset{\lambda}{\colim} \Ex Y)_m.
\endaligned$$
Colimits in $\mathbf{SSet}$ are formed pointwise, we see that $\Ex$
preserves $\lambda$-sequences.

Thus $\Ex^2N$ preserves $\lambda$-sequences.
\item \label{CatCaseiv}
Let $\xymatrix@1{j\co\Lambda^k[m] \ar[r] & \Delta[m]}$ be a
generating acyclic cofibration for $\mathbf{SSet}$. Let the functor
$j'$ be the pushout along $L$ as in the following diagram with $m
\geq 1$.
$$\xymatrix@C=3pc@R=3pc{c\Sd^2\Lambda^k[m] \ar[d]_{c\Sd^2j} \ar[r]^-L & \bfB \ar[d]^{j'}
\\ c\Sd^2\Delta[m] \ar[r] & \bfP}$$
We factor $j'$ into two inclusions
$$\xymatrix{\bfB \ar[r]^i & \bfQ \ar[r] & \bfP}$$
and show that the nerve of each is a weak equivalence.

By Remark \ref{Deltafreecomposites} the only free composites that
occur in the pushout $\bfP$ are of the form $(f_1,f_2)$
$$\xymatrix{\ar[r]^{f_1} & \ar[r]^{f_2} & }$$
where $f_1$ is a morphism in $\bfB$ and $f_2$ is a morphism of
$\Out$ with source in $c\Sd^2\Lambda^k[m]$ and target outside of
$c\Sd^2\Lambda^k[m]$ (see for example the drawing of
$c\Sd^2\Delta[m]$ in Figure \ref{subdivisionfigure}). Hence, $\bfP$ is the union
\begin{equation} \label{decomposition}
\bfP=\overbrace{(\bfB \coprod_{c\Sd^2\Lambda^k[m]} \Out)}^\bfQ  \cup
\overbrace{(\Comp \cup \Cen)}^{\bfR}
\end{equation}
by Proposition \ref{upcloseddecomposition}, all free composites are
in $\bfQ$, and they have the form $(f_1,f_2)$.

We claim that the nerve of the inclusion $\xymatrix@1{i\co\bfB
\ar[r] & \bfQ}$ is a weak equivalence. Let
$\xymatrix@1{\overline{r}\co \bfQ \ar[r] & \bfB}$ be the identity on
$\bfB$, and for any $(v_0, \dots, v_q)\in \Out$ we define
$\overline{r}(v_0, \dots, v_q)=(u_0, \dots, u_p)$ where $(u_0,
\dots, u_p)$ is the maximal subset
$$\{u_0,\dots, u_p\} \subseteq \{v_0,\dots, v_q\}$$
that is in $\bfP\Sd \Lambda^k[m]$ (recall Proposition \ref{outer}
\ref{outerii}). On free composites in $\bfQ$ we then have
$\overline{r}(f_1,f_2)=(f_1,\overline{r}(f_2))$. More conceptually,
we define $\xymatrix@1{\overline{r}\co \bfQ \ar[r] & \bfB}$ using
the universal property of the pushout $\bfQ$ and the maps $1_{\bfB}$
and $Lr$ (the functor $r$ is as in Proposition \ref{outer}
\ref{outerii}).

Then $\overline{r}i=1_\bfB$, and there is a unique natural
transformation $\xymatrix@1{i\overline{r} \ar@{=>}[r] & 1_\bfQ}$
which is the identity morphism on the objects of $\bfB$. Thus
$\xymatrix@1{|Ni|\co|N\bfB|\ar[r] & |N\bfQ|}$ includes $|N\bfB|$ as
a deformation retract of $|N\bfQ|$.

Next we show that the nerve of the inclusion $\xymatrix@1{\bfQ
\ar[r] & \bfP}$ is also a weak equivalence. The intersection of
$\bfQ$ and $\bfR$ in (\ref{decomposition}) is equal to
$$\bfS=\Out\cap(\Comp\cup\Cen).$$ Proposition
\ref{nervecommuteswithpushouthypothesisverification} then implies
that $\bfQ$, $\bfR$, and $\bfS$ satisfy the hypotheses of
Proposition \ref{nervecommuteswithpushout}. Then
\begin{equation} \label{QPinclusion}
\aligned
|N\bfQ| & \cong |N\bfQ| \coprod_{|N\bfS|} |N\bfS| \text{ (pushout along identity) }\\
& \simeq |N\bfQ| \coprod_{|N\bfS|} |N\bfR| \text{ (Prop. \ref{deformationretract} and Gluing Lemma)}\\
& \cong |N\bfQ \coprod_{N\bfS} N\bfR| \text{ (realization is a left adjoint})  \\
& \cong |N(\bfQ \coprod_\bfS \bfR)| \text{ (Prop. \ref{nervecommuteswithpushout} and Prop.
\ref{nervecommuteswithpushouthypothesisverification})} \\
& = |N\bfP|.
\endaligned
\end{equation}
In the second line, for the application of the Gluing Lemma (Lemma
8.12 in \cite{goerssjardine} or Proposition 13.5.4 in
\cite{hirschhorn}), we use two identities and the inclusion
$\xymatrix@1{|N\bfS| \ar[r] & |N\bfR|}$. It is a homotopy
equivalence whose inverse is the retraction in Proposition
\ref{deformationretract}. We conclude that the inclusion
$\xymatrix@1{|N\bfQ| \ar[r] & |N\bfP|}$ is a weak equivalence, as it
is the composite of the morphisms in equation \eqref{QPinclusion}.
It is even a homotopy equivalence by Whitehead's Theorem.

We conclude that $|Nj'|$ is the composite of two weak equivalences
$$\xymatrix@C=3pc{|N\bfB| \ar[r]^{|Ni|} & |N\bfQ| \ar[r] & |N \bfP|}$$
and is therefore a weak equivalence.  By Lemma
\ref{ExPreservesAndReflectsWEs}, the functor $\Ex$ preserves weak
equivalences, so that $\Ex^2Nj'$ is also a weak equivalence of
simplicial sets. Part \ref{KanCorollaryiv} of Corollary
\ref{KanCorollary} then holds, and we have the Thomason model
structure on $\mathbf{Cat}$.
\end{enumerate}
\end{pf}

\section{Pushouts and Colimit Decompositions of $c^n\delta_! \Sd^2
\Delta[m]$} \label{sectionnfolddecompositions}

Next we enhance the proof of the $\mathbf{Cat}$-case to obtain the
$\mathbf{nFoldCat}$-case. The preparations of Section
\ref{barycentric}, \ref{retractionsection}, and \ref{pushoutsection}
are adapted in this section to $n$-fold categorification.

\begin{prop} \label{standardgluingsofdoublecats}
Let $\xymatrix@1{d^i:[m-1] \ar[r] & [m]}$ be the injective order
preserving map which skips $i$. Then the pushout in
$\mathbf{nFoldCat}$
\begin{equation} \label{standardgluingsofdoublecatspushout}
\xymatrix@R=3pc@C=4pc{[m-1] \boxtimes \cdots \boxtimes [m-1] \ar[r]^-{d^i \boxtimes \cdots \boxtimes
d^i}
\ar[d]_-{d^i \boxtimes \cdots \boxtimes d^i} & [m] \boxtimes \cdots \boxtimes [m] \ar[d] \\
[m] \boxtimes \cdots \boxtimes [m] \ar[r]   & \bbP}
\end{equation}
does not have any free composites, and is an $n$-fold poset.
\end{prop}
\begin{pf}
We do the proof for $n=2$.

We consider horizontal morphisms, the proof for vertical morphisms
and more generally squares is similar. We denote the two copies of
$[m] \boxtimes [m]$ by $\bbN_1$ and $\bbN_2$ for convenience. A free
composite occurs whenever there are
$$\xymatrix{f_1\co A_1 \ar[r] & B_1}$$
$$\xymatrix{g_2\co B_2 \ar[r] & C_2}$$
in $\bbN_1$ and $\bbN_2$ respectively such that $B_1$ and $B_2$ are
identified in the pushout, and further, the images of $[m-1]
\boxtimes [m-1]$ contain neither $f_1$ nor $g_2$. Inspection of $d^i
\boxtimes d^i$ shows that this does not occur.
\end{pf}

\begin{rmk}
The gluings of Proposition \ref{standardgluingsofdoublecats} are the
only kinds of gluings that occur in $c^n \delta_! \Sd^2 \Delta[m]$
and $c^n \delta_! \Sd^2 \Lambda^k[m]$ because of the description of
glued simplices in Remark \ref{gluingofsimplices} and the fact that
$c^n \delta_!$ is a left adjoint.
\end{rmk}

\begin{cor}
Consider the pushout $\bbP$ in Proposition
\ref{standardgluingsofdoublecats}. The application of $\delta^\ast
N^n$ to Diagram (\ref{standardgluingsofdoublecatspushout}) is a
pushout and is drawn in Diagram (\ref{gluingtogluing}).
\begin{equation} \label{gluingtogluing}
\xymatrix@R=3pc@C=6pc{\Delta[m-1] \times \cdots \times \Delta[m-1]
\ar[r]^-{\delta^\ast N^n(d^i \boxtimes \cdots \boxtimes d^i)}
\ar[d]_-{\delta^\ast N^n(d^i \boxtimes \cdots \boxtimes d^i)} & \Delta[m] \times \cdots \times \Delta[m] \ar[d] \\
\Delta[m] \times \cdots \times \Delta[m] \ar[r] & \delta^\ast
N^n\bbP}
\end{equation}
\end{cor}
\begin{pf}
The functor $N^n$ preserves a pushout whenever there are no free
composites in that pushout, which is the case here by Proposition
\ref{standardgluingsofdoublecats}. Also, $\delta_!$ is a left
adjoint, so it preserves any pushout.
\end{pf}

The $n$-fold version of Proposition \ref{colimitdecomposition} is as
follows.

\begin{prop} \label{colimitdecompositionnfold}
Let $\bfT$ and $F$ be as in Proposition \ref{colimitdecomposition}.
In particular, $\bfT$ could be $\bfP\Sd \Delta[m], \Cen, \Out,
\Comp$ or $\Comp \cup \Cen$ by Proposition \ref{satisfyhypothesis}.
Then $c^n \delta_!N\bfT$ is the union inside of $\bfT \boxtimes \bfT
\boxtimes \cdots \boxtimes \bfT$ given by
\begin{equation} \label{equ:cndelta!NT}
c^n \delta_!N\bfT=\underset{|U|=m+1}{\underset{U \subseteq \bfT \;
\; \text{\rm lin. ord.}}{\bigcup}} U \boxtimes U \boxtimes \cdots
\boxtimes U.
\end{equation}
Similarly, if $\bfS=\bfP\Sd\Lambda^k[m]$ or $\bfS=\Out \cap (\Comp
\cup \Cen)$, then by Proposition \ref{satisfyhypothesis}, $c^n
\delta_!N\bfS$ is the union inside of $\bfS \boxtimes \bfS \boxtimes
\cdots \boxtimes \bfS$ given by
\begin{equation} \label{equ:cndelta!NTLambda}
c^n \delta_!N\bfS=\underset{|U|=m}{\underset{U \subseteq \bfS \; \;
\text{\rm lin. ord.}}{\bigcup}} U \boxtimes U \boxtimes \cdots
\boxtimes U.
\end{equation}
If $\bfT$ or $\bfS$ is any of the respective posets above, then
$$c^n \delta_! N \bfT \subseteq \bfP\Sd \Delta[m]\boxtimes \bfP\Sd \Delta[m] \boxtimes \cdots \boxtimes \bfP\Sd \Delta[m]$$
$$c^n \delta_! N \bfS \subseteq \bfP\Sd \Delta[m]\boxtimes \bfP\Sd \Delta[m] \boxtimes \cdots \boxtimes \bfP\Sd \Delta[m].$$
\end{prop}
\begin{pf}
For any linearly ordered subposet $U$ of $\bfT$ we have
 $$\aligned c^n\delta_! NU &=  c^n(NU
\boxtimes NU \boxtimes \cdots \boxtimes NU) \\
&= cNU \boxtimes cNU \boxtimes \cdots \boxtimes cNU \\
&=U \boxtimes U \boxtimes \cdots \boxtimes U.
\endaligned$$
Thus we have
$$\aligned c^n \delta_! N\bfT &= c^n \delta_! N(\underset{\bfJ}{\colim} F)
\text{  by Proposition \ref{colimitdecomposition}}\\
&= c^n \delta_!( \underset{\bfJ}{\colim} NF)
\text{ by Corollary \ref{nervecommuteswithcolimitdecomposition}}\\
&=\underset{\bfJ}{\colim} c^n \delta_! NF \text{  because $c^n\delta_!$ is a left adjoint} \\
&=\underset{U \in \bfJ}{\colim} U \boxtimes U \boxtimes \cdots \boxtimes U \\
&=\underset{|U|=m+1}{\underset{U \subseteq \bfT \; \; \text{\rm lin. ord.}}{\bigcup}} U \boxtimes U \boxtimes \cdots  \boxtimes U.
\endaligned$$
This last equality follows for the same reason that $\bfT$ (=colimit of $F$) is the union of the
linearly ordered subposets $U$ of $\bfT$ with exactly $m+1$ elements.  See also Proposition
\ref{standardgluingsofdoublecats}.
\end{pf}

\begin{rmk}
Note that
\begin{equation*}
\bfT \boxtimes \bfT \boxtimes \cdots  \boxtimes \bfT \supsetneq
\underset{|U|=m+1}{\underset{U \subseteq \bfT \; \; \text{\rm lin. ord.}}{\bigcup}}
U \boxtimes U \boxtimes \cdots  \boxtimes U.
\end{equation*}
\end{rmk}

\begin{defn}
An $n$-fold category is an {\it $n$-fold preorder} if for any two
objects $A$ and $B$, there is at most one $n$-cube with $A$ in the
$(0,0, \dots, 0)$-corner and $B$ in the $(1,1,\dots,1)$-corner. If
$\bbD$ is an $n$-fold preorder, we define an ordinary preorder on
$\Obj \bbD$ by $A \leq B :\Longleftrightarrow$ there exists an $n$-cube with $A$ in the
$(0,0, \dots, 0)$-corner and $B$ in the $(1,1,\dots,1)$-corner. We call an $n$-fold preorder
an {\it $n$-fold poset} if $\leq$ is additionally antisymmetric as a preorder on $\Obj \bbD$, that is,
$(\Obj \bbD, \leq)$ is a poset. If $\bbT$ is an $n$-fold preorder and $\bbS$ is a sub-$n$-fold preorder, then $\bbS$ is
{\it down-closed in $\bbT$} if $A \leq B$ and $B \in \bbS$ implies $A \in \bbS$.
 If $\bbT$ is an $n$-fold preorder and $\bbS$ is a sub-$n$-fold preorder, then the
{\it up-closure} of $\bbS$ in $\bbT$ is the full sub-$n$-category of $\bbT$ on the objects $B$ in $\bbT$ such that
$B \geq A$ for some object $A \in \bbS$.
\end{defn}

\begin{examp}
If $\bfT$ is a poset, then $\bfT \boxtimes \bfT \boxtimes \cdots
\boxtimes \bfT$ is an $n$-fold poset, and $(a_1, \dots, a_n) \leq
(b_1, \dots, b_n)$ if and only if $a_i \leq b_i$ in $\bfT$ for all
$1 \leq i \leq n$. If $\bfT$ is as in Proposition
\ref{colimitdecomposition}, then the $n$-fold category
$c^n\delta_!N\bfT$ is also an $n$-fold poset, as it is contained in
the $n$-fold poset $\bfT \boxtimes \bfT \boxtimes \cdots  \boxtimes
\bfT$ by equation \eqref{equ:cndelta!NT}.
\end{examp}

\begin{prop} \label{Lambdadownclosednfold}
The $n$-fold poset $c^n\delta_!N\bfP\Sd\Lambda^k[m]$ is down-closed
in $c^n\delta_!N\bfP\Sd \Delta[m]$.
\end{prop}
\begin{pf}
Suppose $(a_1, \dots, a_n) \leq (b_1, \dots, b_n)$ in
$c^n\delta_!N\bfP\Sd \Delta[m]$ and $(b_1, \dots, b_n) \in
c^n\delta_!N\bfP\Sd\Lambda^k[m]$. We make use of equations \eqref{equ:cndelta!NT} and \eqref{equ:cndelta!NTLambda}
in Proposition \ref{colimitdecompositionnfold}. There exists a linearly
ordered subposet $V$ of $\bfP\Sd\Lambda^k[m]$ such that $|V|=m$ and
$b_1, \dots, b_n \in V$.  There also exists a linearly ordered subposet $U$ of $\bfP\Sd \Delta[m]$ such that $|U|=m+1$
and $a_1, \dots, a_n \in U$. In particular, $\{a_1, \dots, a_n\}$ is linearly ordered.

The preorder on $\Obj c^n\delta_!N\bfP\Sd \Delta[m]$ then implies
that $a_i \leq b_i$ in $\bfP\Sd \Delta[m]$ for all $i$, so that $a_i
\in \bfP\Sd\Lambda^k[m]$ by Proposition \ref{Lambdadownclosed}.
Since the length of a maximal chain in $\bfP\Sd\Lambda^k[m]$ is $m$,
the linearly ordered poset $\{a_1, \dots, a_n\}$ has at most $m$
elements. By Proposition \ref{satisfyhypothesis}, there exists a
linearly ordered subposet $U'$ of $\bfP\Sd\Lambda^k[m]$ such that
$|U'|=m$ and $a_1, \dots, a_n \in U'$. Consequently, $(a_1, a_2,
\ldots, a_n)\in c^n\delta_!N\bfP\Sd\Lambda^k[m]$, again by equation
\eqref{equ:cndelta!NTLambda}.
\end{pf}

\begin{prop} \label{upclosurenfold}
The up-closure of $c^n\delta_!N\bfP\Sd\Lambda^k[m]$ in
$c^n\delta_!N\bfP\Sd \Delta[m]$ is contained in $c^n \delta_!N
\Out$.
\end{prop}
\begin{pf}
An explicit description of all three $n$-fold posets is given in
equations \eqref{equ:cndelta!NT} and \eqref{equ:cndelta!NTLambda} of
Proposition \ref{colimitdecompositionnfold}. Recall that
$\bfP\Sd\Lambda^k[m]$ and $\Out$ satisfy hypothesis
\ref{colimitdecompositioni} of Proposition
\ref{colimitdecomposition} for $m$ and $m+1$ respectively (by
Proposition \ref{satisfyhypothesis}).

Suppose $$A=(a_1, a_2, \ldots, a_n)\leq (b_1, b_2, \ldots, b_n)=B$$
in $c^n\delta_!N\bfP\Sd \Delta[m]$, $A \in
c^n\delta_!N\bfP\Sd\Lambda^k[m]$, and $B \in c^n\delta_!N\bfP\Sd
\Delta[m]$. Then $\{a_1, a_2, \ldots, a_n \} \subseteq U$ for some
linearly ordered subposet $U \subseteq \bfP\Sd\Lambda^k[m]$ with
$\vert U \vert =m$, and $\{b_1, b_2, \ldots, b_n \} \subseteq V$ for
some linearly ordered subposet $V \subseteq \bfP\Sd\Delta[m]$ with
$\vert V \vert =m+1$. We also have $a_i \leq b_i$ in
$\bfP\Sd\Delta[m]$ for all $i$, so that each $b_i$ is in the
up-closure of $\bfP\Sd\Lambda^k[m]$ in $\bfP\Sd\Delta[m]$, namely in
$\Out$. Since equation \eqref{equ:cndelta!NT} holds for $\Out$, we
see $B \in c^n \delta_!N \Out$, and therefore the up-closure of
$c^n\delta_!N\bfP\Sd\Lambda^k[m]$ is contained in $c^n \delta_!N
\Out$.
\end{pf}

\begin{rmk} \label{Deltafreecompositesnfold}
\begin{enumerate}
\item
If $\alpha$ is an $n$-cube in $c^n\delta_!N\bfP\Sd \Delta[m]$ whose
$i$-th target is in $c^n\delta_!N\bfP\Sd\Lambda^k[m]$, then $\alpha$
is in $c^n\delta_!N\bfP\Sd\Lambda^k[m]$.
\item
If $\alpha$ is an $n$-cube in $c^n\delta_!N\bfP\Sd \Delta[m]$ whose
$i$-th source is in $c^n\delta_!N\bfP\Sd\Lambda^k[m]$, then $\alpha$
is in $c^n \delta_!N \Out$.
\end{enumerate}
\end{rmk}
\begin{pf}
\begin{enumerate}
\item
If $\alpha$ is an $n$-cube in $c^n\delta_!N\bfP\Sd \Delta[m]$ whose
$i$-th target is in $c^n\delta_!N\bfP\Sd\Lambda^k[m]$, then its
$(1,1,\ldots,1)$-corner is in $c^n\delta_!N\bfP\Sd\Lambda^k[m]$, as
this corner lies on the $i$-th target. By Proposition
\ref{Lambdadownclosednfold}, we then have $\alpha$ is in
$c^n\delta_!N\bfP\Sd\Lambda^k[m]$.
\item
If $\alpha$ is an $n$-cube in $c^n\delta_!N\bfP\Sd \Delta[m]$ whose
$i$-th source is in $c^n\delta_!N\bfP\Sd\Lambda^k[m]$, then the
$(0,0,\ldots,0)$-corner is in $c^n\delta_!N\bfP\Sd\Lambda^k[m]$, as
this corner lies on the $i$-th source. By Proposition
\ref{upclosurenfold}, we then have $\alpha$ is in $c^n \delta_!N
\Out$.
\end{enumerate}
\end{pf}

Next we describe the diagonal of the nerve of certain $n$-fold
categories as a union of $n$-fold products of standard simplices in
Proposition \ref{diagonaldecomposition}. This proposition is also an
analogue of Corollary \ref{nervecommuteswithcolimitdecomposition}
since it says the composite functor $\delta^*N^nc^n\delta_!N$
preserves colimits of certain posets.

\begin{lem} \label{lem:n-fold_categorification_of_lin_ordered_poset}
For any finite, linearly ordered poset $U$ we have
$$\delta^*N^n c^n\delta_! NU=NU \times NU \times \cdots \times NU.$$
\end{lem}
\begin{pf}
Since $U$ is a finite, linearly ordered poset, $NU$ is isomorphic to
$\Delta[m]$ for some non-negative integer $m$, and we have
$$\aligned \delta^*N^n c^n\delta_! NU &= \delta^*N^n c^n(NU
\boxtimes NU \boxtimes \cdots \boxtimes NU) \\
&=\delta^*N^n (cNU \boxtimes cNU \boxtimes \cdots \boxtimes cNU) \\
&=\delta^*N^n (U \boxtimes U \boxtimes \cdots \boxtimes U) \\
&=\delta^*(NU \boxtimes NU \boxtimes \cdots \boxtimes NU) \\
&=NU \times NU \times \cdots \times NU.
\endaligned$$
\end{pf}

\begin{lem} \label{lem:n-fold_skeletality}
For any finite, linearly ordered poset $U$, the simplicial set
$$\delta^*N^n c^n\delta_! NU=NU \times NU \times \cdots \times NU$$ is
$M$-skeletal for a large enough $M$ depending on $n$ and the
cardinality of $U$.
\end{lem}
\begin{pf}
We prove that there is an $M$ such that all simplices in degrees
greater than $M$ are degenerate.

Without loss of generality, we may assume $U$ is $[m]$. We have
$$\aligned
c^n\delta_! N[m] &= c^n \delta_! \Delta[m] \\
&= c^n \left( \Delta[m] \boxtimes \Delta[m] \boxtimes \cdots
\boxtimes
\Delta[m] \right) \\
&= \left( c\Delta[m] \right) \boxtimes \left( c\Delta[m]\right)
\boxtimes \cdots \boxtimes
\left(c\Delta[m]\right) \\
&= [m] \boxtimes [m] \boxtimes \cdots \boxtimes [m]
\endaligned$$
by Example \ref{examp:categorification_and_external_products}. An
$\ell$-simplex in $\delta^*N^n\left([m] \boxtimes [m] \boxtimes
\cdots \boxtimes [m] \right)$ is an $\ell \times \ell \times \cdots
\times \ell$ array of composable $n$-cubes in $[m] \boxtimes [m]
\boxtimes \cdots \boxtimes [m]$, that is, a collection of $n$
sequences of $\ell$ composable morphisms in $[m]$, namely $\left(
(f^1_j)_j, (f^2_j)_j, \dots, (f^n_j)_j \right)$ where $1\leq j \leq
\ell$ and $f^i_{j+1} \circ f^i_j$ is defined for $j+1\leq\ell$. An
$\ell$-simplex is degenerate if and only if there is a $j_0$ such
that $f^1_{j_0}, f^2_{j_0}, \dots,f^n_{j_0}$ are all identities. An
$\ell$-simplex has $\ell$-many $n$-cubes along its diagonal, namely
$$(f^1_{j}, f^2_{j}, \dots,f^n_{j})$$
for $1 \leq j \leq \ell$. Since $[m]$ is finite, there is an integer
$M$ such that for any $\ell\geq 0$ and any $\ell$-simplex $y$, there
are at most $M$-many nontrivial $n$-cubes in $y$, that is, there are
at most $M$-many tuples
$$(f^1_{j_1}, f^2_{j_2}, \dots,f^n_{j_n})$$ which have at least one
$f^i_{j_i}$ nontrivial.

If $\ell > M$ then at least one of the $\ell$-many $n$-cubes on the
diagonal must be trivial, by the pigeon-hole principle. Hence, for
$\ell > M$, every $\ell$-simplex of $\delta^*N^n c^n\delta_! N[m]$
is degenerate. Finally, $\delta^*N^n c^n\delta_! N[m]$ is
$M$-skeletal.
\end{pf}

\begin{prop} \label{diagonaldecomposition}
Let $m \geq 1$ be a positive integer and $\bfT$ a poset satisfying
the hypotheses \ref{simplicial_colimitdecompositioni} and
\ref{simplicial_colimitdecompositionii} of Proposition
\ref{simplicial_colimitdecomposition}. In particular, $\bfT$ could
be $\bfP\Sd \Delta[m]$, $\Cen$, $\Out$, $\Comp$, or $\Comp \cup
\Cen$ by Proposition \ref{satisfyhypothesis}. Let the functor
$\xymatrix@1{F\co\bfJ \ar[r] & \mathbf{Cat}}$ and the universal
cocone $\xymatrix@1{\pi\co F \ar@{=>}[r] & \Delta_\bfT}$ be as
indicated in Proposition \ref{colimitdecomposition}. Then
$$\aligned \delta^* N^n c^n \delta_! N\bfT &=\underset{\bfJ}{\colim}
\delta^* N^n c^n \delta_! NF \\ &= \underset{\bfJ}{\colim} (NF
\times \cdots \times NF)\endaligned$$ where $NF(U)$ is isomorphic to
$\Delta[m-1]$ or $\Delta[m]$ for all $U \in \bfJ$. Similarly, the
simplicial sets $\delta^* N^n c^n \delta_! N(\bfP\Sd\Lambda^k[m])$
and
$$\delta^* N^n c^n \delta_! N(\Out \cap (\Comp \cup \Cen))$$ are
each a colimit of simplicial sets of the form $\Delta[m-2]\times
\cdots \times \Delta[m-2]$ and $\Delta[m-1] \times \cdots \times
\Delta[m-1]$. (By definition $[-1]=\emptyset$.)
\end{prop}
\begin{pf}
We first prove directly that $\delta^* N^n c^n \delta_! N\bfT$ is a
colimit of $\xymatrix@1{\delta^* N^n c^n \delta_! NF \co \bfJ \ar[r]
& \mathbf{SSet}}$ along the lines of the proof of Proposition
\ref{simplicial_colimitdecomposition}.

Let $M>m$ be a large enough integer such that the simplicial set
$\delta^*N^n c^n\delta_! N[m]$ is $M$-skeletal. Such an $M$ is
guaranteed by Lemma \ref{lem:n-fold_skeletality}.

Suppose $S\in\mathbf{SSet}$ and $\xymatrix@1{\alpha:\delta^* N^n c^n
\delta_! NF \ar@{=>}[r] & \Delta_S}$ is a natural transformation. We
induce a morphism of simplicial sets $$\xymatrix@1{G:\delta^* N^n
c^n \delta_! N\bfT \ar[r] & S}$$ by defining $G$ on the $M$-skeleton
as follows.

As in the proof of Proposition
\ref{simplicial_colimitdecomposition}, $\Delta_{M}$ denotes the full
subcategory of $\Delta$ on the objects $[0], [1], \dots, [M]$ and
$\xymatrix@1{\text{tr}_{M} \co \mathbf{SSet} \ar[r] &
\mathbf{Set}^{\Delta_{M}^\text{op}}}$ denotes the $M$-th truncation
functor. The truncation $\text{tr}_{M} (\delta^* N^n (c^n \delta_!
N\bfT))$ is a union of the truncated simplicial subsets
$\text{tr}_{M}(\delta^* N^n (c^n \delta_! N\bfV))$ for $V \in \bfJ$
with $\vert V \vert=m+1$, since
\begin{itemize}
\item
$c^n\delta_!N\bfT$ is a union of such $c^n \delta_! N\bfV$ by
Proposition \ref{colimitdecompositionnfold},
\item
any maximal linearly ordered subset of $\bfT$ has $m+1$ elements,
and
\item
$\delta^*N^n$ preserves unions.
\end{itemize}
We define
$$\xymatrix{G_{M}\vert_{\text{tr}_{M}(\delta^* N^n (c^n \delta_! N\bfV))}
\co \text{tr}_{M}(\delta^* N^n (c^n \delta_! N\bfV)) \ar[r] &
\text{tr}_{M}S}$$ simply as $\text{tr}_{M}\alpha_{V}$.

The morphism $G_M$ is well-defined, for if $0\leq\ell \leq M$ and $x
\in (\text{tr}_{M}(\delta^* N^n c^n \delta_! N\bfV)_\ell$ and $x \in
(\text{tr}_{M}(\delta^* N^n c^n \delta_! N\bfV)_\ell$ with $\vert V
\vert=m+1=\vert V \vert$, then $V$ and $V'$ can be connected by a
sequence $W^0,W^1, \dots, W^k$ of $(m+1)$-element linearly ordered
subsets of $\bfT$ that all contain the linearly ordered subposet $x$
and satisfy the properties in hypothesis
\ref{simplicial_colimitdecompositionii}. By a naturality argument as
in the proof of Proposition \ref{colimitdecomposition}, we have have
a string of equalities
$$\alpha_{W^0}(x)=\alpha_{W^1}(x)=\cdots =\alpha_{W^k}(x),$$
and we conclude $\alpha_{V}(x)=\alpha_{V'}(x)$ so that $G_M(x)$ is
well defined.

By definition $\Delta_{G_M} \circ \text{tr}_M N\pi=\text{tr}_M
\alpha$. We may extend this to non-truncated simplicial sets by
recalling from above that the simplicial set $\delta^* N^n c^n
\delta_! N\bfT$ is {\it $M$-skeletal}, that is, the counit inclusion
$$\xymatrix@1{\text{sk}_M\text{tr}_M(\delta^* N^n c^n \delta_! N\bfT) \ar[r] & \delta^* N^n c^n \delta_! N\bfT}$$ is
the identity.

Thus $G_M$ extends to $\xymatrix@1{G\co N\bfT \ar[r] & S}$ and
$\Delta_{G} \circ N\pi=\alpha$.

Lastly, the morphism $G$ is unique, since the simplicial subsets
$\delta^* N^n c^n \delta_! N\bfV$ for $\vert V \vert=m+1$ in $\bfJ$
cover $\delta^* N^n c^n \delta_! N\bfT$ by hypothesis
\ref{simplicial_colimitdecompositioni}.

So far we have proved $\delta^* N^n c^n \delta_! N\bfT
=\underset{\bfJ}{\colim} \delta^* N^n c^n \delta_! NF$. It only
remains to show $\underset{\bfJ}{\colim} \delta^* N^n c^n \delta_!
NF=\underset{\bfJ}{\colim} (NF \times \cdots \times NF)$. But this
follows from Lemma
\ref{lem:n-fold_categorification_of_lin_ordered_poset} and that fact
that $FV=V$ for all $V \in \bfJ$.
\end{pf}

The $n$-fold version of Proposition \ref{deformationretract} is the
following.

\begin{cor} \label{nfolddeformationretract}
The space $|\delta^*N^n c^n \delta_! N (\Out \cap (\Comp \cup \Cen))|$
includes into the space $|\delta^* N^n c^n \delta_!N(\Comp \cup
\Cen)|$ as a deformation retract.
\end{cor}
\begin{pf}
Recall that realization $|\cdot|$ commutes with colimits, since it
is a left adjoint, and that $|\cdot|$ also commutes with products.
We do the multi-stage deformation retraction of Proposition
\ref{deformationretract} to each factor $|\Delta[m]|$ of
$|\Delta[m]|\times \cdots \times |\Delta[m]|$ in the colimit of
Proposition \ref{diagonaldecomposition}. This is the desired
deformation retraction of $|\delta^* N^n c^n \delta_! N(\Comp \cup
\Cen)|$ to $|\delta^*N^n c^n \delta_! N(\Out \cap (\Comp \cup
\Cen))|$.
\end{pf}

\begin{prop} \label{PushoutDescription}
Consider $n=2$. Let $\xymatrix@1{j\co\Lambda^k[m] \ar[r] &
\Delta[m]}$ be a generating acyclic cofibration for $\mathbf{SSet}$,
$\bbB$ a double category, and $L$ a double functor as below. Then
the pushout $\bbQ$ in the diagram
\begin{equation} \label{PushoutDiagram}
\xymatrix@C=3pc@R=3pc{c^2\delta_!\Sd^2\Lambda^k[m]
\ar[d]_{c^2\delta_!\Sd^2j} \ar[r]^-L  & \bbB \ar[d]
\\ c^2\delta_!N \Out \ar[r] & \bbQ}
\end{equation}
has the following form.
\begin{enumerate}
\item \label{PushoutDescription_object_set}
The object set of $\bbQ$ is the pushout of the object sets.
\item \label{PushoutDescription_horizontal}
The set of horizontal morphisms of $\bbQ$ consists of the set of
horizontal morphisms of $\bbB$, the set of horizontal morphisms of
$c^2\delta_!N\Out$, and the set of formal composites of the form
$$\xymatrix@C=3pc{\ar[r]^{f_1} & \ar[r]^{(1,f_2)} & }$$
where $f_1$ is a horizontal morphism in $\bbB$, $f_2$ is a morphism
in $\Out$, and the target of $f_1$ is the source of $(1,f_2)$ in
$\Obj \bbQ$.
\item \label{PushoutDescription_vertical}
The set of vertical morphisms of $\bbQ$ consists of the set of
vertical morphisms of $\bbB$, the set of vertical morphisms of
$c^2\delta_!N\Out$, and the set of formal composites of the form
$$\xymatrix@R=3pc{\ar[d]^{g_1} \\ \ar[d]^{(g_2,1)} \\ \\}$$
where $g_1$ is a vertical morphism in $\bbB$, $g_2$ is a morphism in
$\Out$, and the target of $g_1$ is the source of $(g_2,1)$ in $\Obj
\bbQ$.
\item \label{PushoutDescription_squares}
The set of squares of $\bbQ$ consists of the set of squares of
$\bbB$, the set of squares of $c^2\delta_!N\Out$, and the set of
formal composites of the following three forms.
\begin{enumerate}
\item \label{PushoutDescription_squares_a}
$$\xymatrix@R=4pc@C=4pc{ \ar[r]^{f_1} \ar[d]_{g_1}
\ar@{}[dr]|{\alpha_1} & (W,A') \ar[r]^{(1_W,f_2)}
\ar[d]|{\tb{(g,1_{A'})}} & (W,B') \ar[d]^{(g,1_{B'})} \\
 \ar[r]_{p_1} & (A,A') \ar[r]_{(1_A,f_2)}
 & (A,B')}$$
\item \label{PushoutDescription_squares_b}
$$\xymatrix@R=4pc@C=4pc{\ar[r]^{f_1} \ar[d]_{g_1}
\ar@{}[dr]|{\beta_1} &
\ar[d]^{q_1}  \\
(A,W') \ar[r]|{\lr{(1_A,f)}} \ar[d]_{(g_2,1_{W'})} & (A,A')
\ar[d]^{(g_2,1_{A'})} \\
(B,W') \ar[r]_{(1_B,f)} & (B,A')  }$$
\item \label{PushoutDescription_squares_c}
$$\xymatrix@R=4pc@C=4pc{\ar[r]^{f_1} \ar[d]_{g_1}
\ar@{}[dr]|{\gamma_1} & (W,A') \ar[r]^{(1_W,f_2)}
\ar[d]|{\tb{(g,1_{A'})}} & (W,B') \ar[d]^{(g,1_{B'})} \\
(A,W') \ar[r]|{\lr{(1_A,f)}} \ar[d]_{(g_2,1_{W'})} & (A,A')
\ar[r]|{\lr{(1_A,f_2)}} \ar[d]|{\tb{(g_2,1_{A'})}} &
(A,B') \ar[d]^{(g_2,1_{B'})} \\
(B,W') \ar[r]_{(1_B,f)} & (B,A') \ar[r]_{(1_B,f_2)} & (B,B') }$$
\end{enumerate}
where $\alpha_1,\beta_1,\gamma_1$ are squares in $\bbB$, the
horizontal morphisms $f_1,p_1$ are in $\bbB$, the vertical morphisms
$g_1,q_1$ are in $\bbB$, and the morphisms $f$, $f_2$, $g$, $g_2$
are in $\Out$. Further, each boundary of each square in
$c^2\delta_!N\Out$ must belong to a linearly ordered subset of
$\Out$ of cardinality $m+1$ (see Proposition
\ref{colimitdecompositionnfold}). So for example, $f$ and $g_2$ must
belong to a linearly ordered subset of $\Out$ of cardinality $m+1$,
and $f_2$ and $g$ must belong to another linearly ordered subset of
$\Out$ of cardinality $m+1$. Of course, the sources and targets in
each of \ref{PushoutDescription_squares_a},
\ref{PushoutDescription_squares_b}, and
\ref{PushoutDescription_squares_c} must match appropriately.
\end{enumerate}
\end{prop}
\begin{pf}
All of this follows from the colimit formula in $\mathbf{DblCat}$,
which is Theorem 4.6 of \cite{fiorepaolipronk1}, and is also a
special case of Proposition
\ref{prop:forgetful_admits_right_adjoint} in the present paper. The
horizontal and vertical 1-categories of $\bbQ$ are the pushouts of
the horizontal and vertical 1-categories, so
\ref{PushoutDescription_object_set} follows, and then
\ref{PushoutDescription_horizontal} and
\ref{PushoutDescription_vertical} follow from Remark
\ref{Deltafreecomposites}. To see \ref{PushoutDescription_squares},
one observes that the only free composite pairs of squares that can
occur are of the first two forms, again from Remark
\ref{Deltafreecomposites}. Certain of these can be composed with a
square in $c^2 \delta_! N \Out$ to obtain the third form. No further
free composites can be obtained from these ones because of Remark
\ref{Deltafreecomposites} and the special form of $c^2 \delta_! N
\Out$.
\end{pf}

\begin{prop} \label{pushoutsimplexdescription}
Consider $n=2$ and the pushout $\bbQ$ in diagram
(\ref{PushoutDiagram}). Then any $q$-simplex in $\delta^*N^2\bbQ$ is
a $q\times q$-matrix of composable squares of $\bbQ$ which has the
form in Figure \ref{pushoutsimplex}.
\begin{figure}
\setlength{\unitlength}{1mm}
\begin{picture}(50,50)
\put(0,0){\line(1,0){50}} \put(0,0){\line(0,1){50}}
\put(50,0){\line(0,1){50}} \put(0,50){\line(1,0){50}}
\put(30,30){\line(0,1){20}} \put(35,30){\line(0,1){20}}
\put(0,30){\line(1,0){35}} \put(0,35){\line(1,0){35}}
\put(15,40){\makebox(0,0)[b]{$\bbB$}}
\put(32.5,31){\makebox(0,0)[b]{$c$}}
\put(32.5,40){\makebox(0,0)[b]{$a$}}
\put(15,31){\makebox(0,0)[b]{$b$}}
\put(25,15){\makebox(0,0)[b]{$c^2\delta_!N\Out$}}
\end{picture}
\caption{A $q$-simplex in $\delta^*N^2\bbQ$.}\label{pushoutsimplex}
\end{figure}
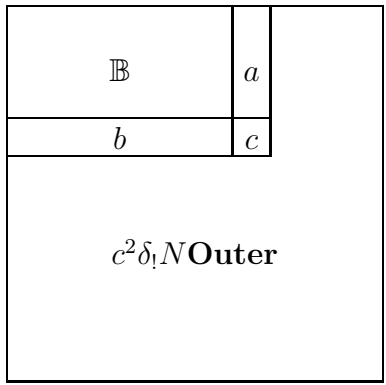
The submatrix labelled $\bbB$ is a matrix of squares in $\bbB$. The
submatrix labelled $a$ is a single column of squares of the form
\ref{PushoutDescription_squares_a} in Proposition
\ref{PushoutDescription} \ref{PushoutDescription_squares} (the
$\alpha_1$'s may be trivial). The submatrix labelled $b$ is a single
row of squares of the form \ref{PushoutDescription_squares_b} in
Proposition \ref{PushoutDescription}
\ref{PushoutDescription_squares} (the $\beta_1$'s may be trivial).
The submatrix labelled $c$ is a single square of the form
\ref{PushoutDescription_squares_c} in Proposition
\ref{PushoutDescription} \ref{PushoutDescription_squares} (part of
the square may be trivial). The remaining squares in the $q$-simplex
are squares of $c^2\delta_!N\Out$.
\end{prop}
\begin{pf}
These are the only composable $q\times q$-matrices of squares
because of the special form of the horizontal and vertical
1-categories.
\end{pf}

\begin{rmk}
The analogues of Propositions \ref{PushoutDescription} and
\ref{pushoutsimplexdescription} clearly hold in higher dimensions as
well, only the notation gets more complicated. Proposition
\ref{prop:forgetful_admits_right_adjoint} provides the key to
proving the higher dimensional versions, namely, it allows us to
calculate the pushout in $\mathbf{nFoldCat}$ in steps: first the
object set of the pushout, then sub-1-categories of the pushout in
all $n$-directions, then the squares in the sub-double-categories of
the pushout in each direction $ij$, then the cubes in the
sub-3-fold-categories of the pushout in each direction $ijk$, and so
on. Since we do not need the explicit formulations of Propositions
\ref{PushoutDescription} and \ref{pushoutsimplexdescription} for
$n>2$ in this paper, we refrain from stating and proving them. In
fact, we do not even need the case $n=2$ for this paper; we only
presented Propositions \ref{PushoutDescription} and
\ref{pushoutsimplexdescription} as an illustration of how the
pushout in $\mathbf{nFoldCat}$ works in a specific case.
\end{rmk}

The $n$-fold version of \ref{nervecommuteswithpushout} is the
following.

\begin{prop} \label{nfoldnervecommuteswithpushout}
Suppose $\bbQ$, $\bbR$, and $\bbS$ are $n$-fold categories, and
$\bbS$ is an $n$-foldly full $n$-fold subcategory of $\bbQ$ and
$\bbR$ such that
\begin{enumerate}
\item \label{nfoldnervecommuteswithpushouti}
If $\xymatrix@1{f:x \ar[r] & y }$ is a 1-morphism in $\bbQ$ (in any
direction) and $x \in \bbS$, then $y \in \bbS$,
\item \label{nfoldnervecommuteswithpushoutii}
If $\xymatrix@1{f:x \ar[r] & y }$ is a 1-morphism in $\bbR$ (in any
direction) and $x \in \bbS$, then $y \in \bbS$.
\end{enumerate}
Then the nerve of the pushout of $n$-fold categories is the pushout
of the nerves.
\begin{equation}
N^n(\bbQ \coprod_\bbS \bbR) \cong N^n\bbQ \coprod_{N^n\bbS} N^n\bbR
\end{equation}
\end{prop}
\begin{pf}

We claim that there are no free composite $n$-cubes in the pushout
$\bbQ \coprod_\bbS \bbR$. Suppose that $\alpha$ is an $n$-cube in
$\bbQ$ and $\beta$ is an $n$-cube in $\bbR$ and that these are
composable in the $i$-th direction. In other words, the $i$-th
target of $\alpha$ is the $i$-th source of $\beta$, which we will
denote by $\gamma$. Then $\gamma$ must be an $(n-1)$-cube in $\bbS$,
as it lies in both $\bbQ$ and $\bbR$. Since the corners of $\gamma$
are in $\bbS$, we can use hypothesis
\ref{nfoldnervecommuteswithpushoutii} to conclude that all corners
of $\beta$ are in $\bbS$ by travelling along edges that emanate from
$\gamma$. By the fullness of $\bbS$, the cube $\beta$ is in $\bbS$,
and also $\bbQ$. Then $\beta \circ_i \alpha$ is in $\bbQ$ and is not
free.

If $\alpha$ is in $\bbR$ and $\beta$ is in $\bbQ$, we can similarly
conclude that $\beta$ is in $\bbS$, $\beta \circ_i \alpha$ is in
$\bbR$, and $\beta \circ_i \alpha$ is not a free composite.

Thus, the pushout $\bbQ \coprod_\bbS \bbR$ has no free composite
$n$-cubes, and hence no free composites of any cells at all.

Let $(\alpha_{\overline{j}})_{\overline{j}}$ be a $p$-simplex in
$N^n(\bbQ \coprod_\bbS \bbR)$. Then each $\alpha_{\overline{j}}$ is
an $n$-cube in $\bbQ$ or $\bbR$, since there are no free composites.
By repeated application of the argument above, if
$\alpha_{(0,\ldots,0)}$ is in $\bbQ$ then every
$\alpha_{\overline{j}}$ is in $\bbQ$. Similarly, if
$\alpha_{(0,\ldots,0)}$ is in $\bbR$ then every
$\alpha_{\overline{j}}$ is in $\bbR$. Thus we have a morphism
$\xymatrix{N^n(\bbQ \coprod_\bbS \bbR) \ar[r] & N^n\bbQ
\coprod_{N^n\bbS} N^n\bbR}$. Its inverse is the canonical morphism
$\xymatrix{N^n\bbQ \coprod_{N^n\bbS} N^n\bbR \ar[r] & N^n(\bbQ
\coprod_\bbS \bbR) }$.

Note that we have not used the higher dimensional versions of
Propositions \ref{PushoutDescription} and
\ref{pushoutsimplexdescription} anywhere in this proof.
\end{pf}

\section{Thomason Structure on $\mathbf{nFoldCat}$}\label{Thomasonnfoldsection}

We apply Corollary \ref{KanCorollary} to transfer across the
adjunction below.

\begin{equation}\label{nfoldcatadjunction}
\xymatrix@C=4pc{\mathbf{SSet} \ar@{}[r]|{\perp}
\ar@/^1pc/[r]^-{\Sd^2} &  \ar@/^1pc/[l]^-{\Ex^2} \mathbf{SSet}
\ar@{}[r]|{\perp} \ar@/^1pc/[r]^-{\delta_!} &
\ar@/^1pc/[l]^-{\delta^\ast} \mathbf{SSet^n} \ar@{}[r]|{\perp}
\ar@/^1pc/[r]^-{c^n} & \ar@/^1pc/[l]^-{N^n} \mathbf{nFoldCat}}
\end{equation}

\begin{prop} \label{Ex_squared}
Let $F$ be an $n$-fold functor. Then the morphism of simplicial sets
$\delta^\ast N^nF$ is a weak equivalence if and only if
$\Ex^2\delta^\ast N^nF$ is a weak equivalence.
\end{prop}
\begin{pf}
This follows from two applications of Lemma
\ref{ExPreservesAndReflectsWEs}.
\end{pf}

\begin{thm} \label{MainModelStructure}
There is a model structure on $\mathbf{nFoldCat}$ in which an
$n$-fold functor $F$ a weak equivalence respectively fibration if
and only if $\Ex^2\delta^\ast N^nF$ is a weak equivalence
respectively fibration in $\mathbf{SSet}$. Moreover, this model
structure on $\mathbf{nFoldCat}$ is cofibrantly generated with
generating cofibrations
$$\{\xymatrix{c^n\delta_!\Sd^2\partial\Delta[m] \ar[r] & c^n\delta_!\Sd^2\Delta[m]}|\; m \geq 0\}$$
and generating acyclic cofibrations
$$\{\xymatrix{c^n\delta_!\Sd^2\Lambda^k[m] \ar[r] & c^n\delta_!\Sd^2\Delta[m]}|\; 0 \leq k
\leq m \text{ and } m \geq 1\}.$$
\end{thm}
\begin{pf}
We apply Corollary \ref{KanCorollary}.
\begin{enumerate}
\item
The $n$-fold categories $c^n\delta_!\Sd^2\partial\Delta[m]$ and
$c^n\delta_!\Sd^2\Lambda^k[m]$ each have a finite number of
$n$-cubes, hence they are finite, and are small with respect to
$\mathbf{nFoldCat}$. For a proof, see Proposition 7.7 of
\cite{fiorepaolipronk1} and the remark immediately afterwards.
\item
This holds as in the proof of \ref{CatCaseii} in Theorem
\ref{CatCase}.
\item
The $n$-fold nerve functor $N^n$ preserves filtered colimits. Every
ordinal is filtered, so $N^n$ preserves $\lambda$-sequences. The
functor $\delta^\ast$ preserves all colimits, as it is a left
adjoint. The functor $\Ex$ preserves $\lambda$-sequences as in the
proof of \ref{CatCaseiii} in Theorem \ref{CatCase}.
\item
Let $\xymatrix@1{j\co\Lambda^k[m] \ar[r] & \Delta[m]}$ be a
generating acyclic cofibration for $\mathbf{SSet}$. Let the functor
$j'$ be the pushout along $L$ as in the following diagram with $m
\geq 1$.
\begin{equation} \label{L}
\xymatrix@C=3pc@R=3pc{c^n\delta_!\Sd^2\Lambda^k[m]
\ar[d]_{c^n\delta_!\Sd^2j} \ar[r]^-L  & \bbB \ar[d]^{j'}
\\ c^n\delta_!\Sd^2\Delta[m] \ar[r] & \bbP}
\end{equation}
We factor $j'$ into two inclusions
\begin{equation} \label{nfoldinclusions}
\xymatrix{\bbB \ar[r]^i & \bbQ \ar[r] & \bbP}
\end{equation}
and show that $\delta^*N^n$ applied to each yields a weak
equivalence. For the first inclusion $i$, we will see in Lemma
\ref{nfoldsimplicialdefretract} that $\delta^* N^n i$ is a weak
equivalence of simplicial sets.

By Remark \ref{Deltafreecompositesnfold}, the only free composites
of an $n$-cube in $c^n \delta_! \Sd^2 \Delta[m]$ with an $n$-cube in
$\bbB$ that can occur in $\bbP$ are of the form $\beta \circ_i
\alpha$ where $\alpha$ is an $n$-cube in $\bbB$ and $\beta$ is an
$n$-cube in $c^n \delta_! N \Out$ with $i$-th source in $c^n
\delta_! N \bfP \Sd \Lambda^k[m]$ and $i$-th target outside of $c^n
\delta_! N \bfP \Sd \Lambda^k[m]$. Of course, there are other free
composites in $\bbP$, most generally of a form analogous to
Proposition \ref{PushoutDescription}
\ref{PushoutDescription_squares_c}, but these are obtained by
composing the free composites of the form $\beta \circ_i \alpha$
above.  Hence $\bbP$ is the union
\begin{equation} \label{decompositionnfold}
\bbP=\overbrace{(\bbB \coprod_{c^n\delta_!N\bfP\Sd\Lambda^k[n]}
c^n\delta_!N\Out)}^\bbQ  \cup \overbrace{(c^n\delta_!N(\Comp \cup
\Cen ))}^{\bbR}.
\end{equation}
Note that we have not used the higher dimensional versions of
Propositions \ref{PushoutDescription} and
\ref{pushoutsimplexdescription} to draw this conclusion.

We show that $\delta^*N^n$ applied to the second inclusion
$\xymatrix@1{\bbQ \ar[r] & \bbP}$ in equation
\eqref{nfoldinclusions} is a weak equivalence. The intersection of
$\bbQ$ and $\bbR$ in (\ref{decompositionnfold}) is equal to
$$\aligned
\bbS&=c^n\delta_!N(\Out)\cap c^n\delta_!N(\Comp\cup\Cen)
\\ &=c^n\delta_!N(\Out\cap(\Comp\cup\Cen)).
\endaligned$$ Propositions
\ref{nervecommuteswithpushouthypothesisverification} and
\ref{colimitdecompositionnfold} then imply that $\bbQ$, $\bbR$, and
$\bbS$ satisfy the hypotheses of Proposition
\ref{nfoldnervecommuteswithpushout}. Then
\begin{equation*}
\aligned
|\delta^*N^n\bbQ| & \cong |\delta^*N^n\bbQ| \coprod_{|\delta^*N^n\bbS|} |\delta^*N^n\bbS| \text{ (pushout along identity) }\\
& \simeq |\delta^*N^n\bbQ| \coprod_{|\delta^*N^n\bbS|} |\delta^*N^n\bbR| \text{ (Cor. \ref{nfolddeformationretract} and Gluing Lemma)}\\
& \cong |\delta^*\left(N^n\bbQ \coprod_{N^n\bbS} N^n\bbR \right)| \text{ (the functors $|\cdot|$ and $\delta^*$ are left adjoints})  \\
& \cong |\delta^*N^n(\bbQ \coprod_\bbS \bbR)| \text{ (Prop. \ref{nfoldnervecommuteswithpushout}) } \\
& = |\delta^*N^n\bbP|.
\endaligned
\end{equation*}
In the second line, for the application of the Gluing Lemma, we use
two identities and the inclusion $\xymatrix@1{|\delta^*N^n\bbS|
\ar[r] & |\delta^*N^n\bbR|}$. It is a homotopy equivalence whose
inverse is the retraction in Corollary
\ref{nfolddeformationretract}. We conclude that the inclusion
$\xymatrix@1{|\delta^*N^n\bbQ| \ar[r] & |\delta^*N^n\bbP|}$ is a
weak equivalence, as it is the composite of the morphisms above. It
is even a homotopy equivalence by Whitehead's Theorem.

We conclude that $|\delta^*N^n j'|$ is the composite of two weak
equivalences
$$\xymatrix@C=4pc{|\delta^*N^n\bbB| \ar[r]^{|\delta^*N^ni|} & |\delta^*N^n\bbQ| \ar[r] & |\delta^*N^n \bbP|}$$
and is therefore a weak equivalence. Thus $\delta^*N^n j'$ is a weak
equivalence of simplicial sets. By Lemma
\ref{ExPreservesAndReflectsWEs}, the functor $\Ex$ preserves weak
equivalences, so that $\Ex^2\delta^*N^n j'$ is also a weak
equivalence of simplicial sets. Part \ref{KanCorollaryiv} of
Corollary \ref{KanCorollary} then holds, and we have the Thomason
model structure on $\mathbf{nFoldCat}$.
\end{enumerate}
\end{pf}

\begin{lem} \label{nfoldsimplicialdefretract}
The inclusion $\xymatrix@1{\delta^\ast N^n i\co \delta^\ast N^n\bbB
\ar[r] & \delta^\ast N^n \bbQ}$ embeds the simplicial set
$\delta^\ast N^n\bbB$ into $\delta^\ast N^n \bbQ$ as a simplicial
deformation retract.
\end{lem}
\begin{pf}
Recall $\xymatrix@1{i\co \bbB \ar[r] & \bbQ}$ is the inclusion in
equation \eqref{nfoldinclusions} and $\bbQ$ is defined as in
equation \eqref{decompositionnfold}. We define an $n$-fold functor
$\xymatrix@1{\overline{r}\co \bbQ \ar[r] & \bbB}$ using the
universal property of the pushout $\bbQ$ and the functor from
Proposition \ref{outer} \ref{outerii} called $\xymatrix@1{r\co \Out
\ar[r] & \bfP\Sd \Lambda^k[m]}$. If $(v_0, \dots, v_q)\in \Out$ then
$r(v_0, \dots, v_q):=(u_0, \dots, u_p)$ where $(u_0, \dots, u_p)$ is
the maximal subset
$$\{u_0,\dots, u_p\} \subseteq \{v_0,\dots, v_q\}$$
that is in $\bfP\Sd \Lambda^k[m]$. We have
$$\aligned
c^n\delta_!N\bfP\Sd \Lambda^k[m]&= \underset{|U|=m}{\underset{U
\subseteq \bfP\Sd \Lambda^k[m] \; \;
\text{\rm lin. ord.}}{\bigcup}} U \boxtimes U \boxtimes \cdots  \boxtimes U \\
&\subseteq \underset{|U|=m+1}{\underset{U \subseteq \Out \; \;
\text{\rm lin. ord.}}{\bigcup}} U \boxtimes U \boxtimes \cdots  \boxtimes U \\
&= c^n\delta_!N\Out.
\endaligned$$
Recall $L$ is the $n$-fold functor in diagram (\ref{L}). We define
$\overline{r}$ on $c^n\delta_!N\Out$ to be
$$\xymatrix@1{L \circ (r\boxtimes r \boxtimes \cdots \boxtimes r) \co c^n\delta_!N\Out \ar[r] & \bbB}$$
and we define $\overline{r}$ to be the identity on $\bbB$. This induces the desired $n$-fold functor
$\xymatrix@1{\overline{r}\co \bbQ \ar[r] & \bbB}$ by the universal property of the pushout $\bbQ$.

By definition we have $\overline{r}i=1_\bbB$. We next define an
$n$-fold natural transformation $\xymatrix@1{\overline{\alpha}\co
i\overline{r} \ar@{=>}[r] & 1_\bbQ}$ (see Definition
\ref{defn_nfold_nat_transf}), which will induce a simplicial
homotopy from $\delta^*N^n(i\overline{r})$ to $1_{\delta^*N^n\bbQ}$
as in Proposition \ref{nfoldnat_gives_simplicial_homotopy}. Let
$$\xymatrix@1{f_1\co \bbB \ar[r] & \bbB^{[1]\boxtimes \cdots
\boxtimes [1]}}$$ $$\xymatrix@1{f_2\co c^n\delta_!N\Out \ar[r] &
\bbB^{[1]\boxtimes \cdots [1]}}$$ be the $n$-fold functors
corresponding to the $n$-fold natural transformations
$$\xymatrix@1{pr_\bbB \co \bbB \times ([1]\boxtimes \cdots \boxtimes
[1]) \ar[r] & \bbB}$$  $$\xymatrix@1{L \circ (\alpha \boxtimes
\cdots \boxtimes \alpha)\co c^n\delta_!N\Out \times ([1]\boxtimes
\cdots \boxtimes [1]) \ar[r] & \bbB}$$ (recall $\mathbf{nFoldCat}$
is Cartesian closed by Ehresmann--Ehresmann \cite{ehresmannthree},
the definition of $\alpha$ in Proposition \ref{outer} \ref{outerii},
and Example
\ref{examp:n_naturaltransfs_yield_an_nfold_naturaltransf}). Then the
necessary square involving $f_1$, $f_2$, $L$ and the inclusion
$$\xymatrix{c^n\delta_!N\bfP\Sd \Lambda^k[m] \ar[r] & c^n\delta_!N\Out}$$ commutes
($\alpha\boxtimes \cdots \boxtimes \alpha$ is trivial on
$c^n\delta_!N\bfP\Sd \Lambda^k[m]$), so we have an $n$-fold functor
$\xymatrix@1{f\co \bbQ \ar[r] & \bbB^{[1]\boxtimes \cdots \boxtimes
[1]}}$, which corresponds to an $n$-fold natural transformation
$$\xymatrix{\overline{\alpha}\co i\overline{r} \ar@{=>}[r] &
1_\bbQ}.$$ Thus $\overline{\alpha}$ induces a simplicial homotopy
from $\delta^*N^n(i)\circ \delta^*N^n(\overline{r})$ to
$1_{\delta^*N^n\bbQ}$ and from above we have
$\delta^*N^n(\overline{r}) \circ
\delta^*N^n(i)=1_{\delta^*N^n\bbB}$. This completes the proof that
the inclusion $\xymatrix@1{\delta^\ast N^n i\co \delta^\ast N^n\bbB
\ar[r] & \delta^\ast N^n \bbQ}$ embeds the simplicial set
$\delta^\ast N^n\bbB$ into $\delta^\ast N^n \bbQ$ as a simplicial
deformation retract.

We next write out what this simplicial homotopy is in the case
$n=2$. We denote by $\sigma$ this simplicial homotopy from
$\delta^*N^n(i\overline{r})$ to $1_{\delta^*N^n\bbQ}$. For each $q$,
we need to define $q+1$ maps
$\xymatrix@1{\sigma_\ell:(\delta^*N^n\bbQ)_q \ar[r] &
(\delta^*N^n\bbQ)_{q+1}}$ compatible with the face and degeneracy
maps, $\delta^*N^n(i\overline{r})$, and $1_{\delta^*N^n\bbQ}$. We
define $\sigma_\ell$ on a $q$-simplex $\alpha$ of the form in
Proposition \ref{pushoutsimplexdescription}. {\it This $q$-simplex
$\alpha$ has nothing to do with the $n$-fold natural transformation
$\alpha$ above.} Suppose that the unique square of type
\ref{PushoutDescription_squares_c} of Proposition
\ref{PushoutDescription} is in entry $(u,v)$ and $u \leq v$.

If $\ell<u$, then $\sigma_\ell(\alpha)$ is obtained from $\alpha$ by
inserting a row of vertical identities between rows $\ell$ and
$\ell+1$ of $\alpha$, as well as a column of horizontal identity
squares between columns $\ell$ and $\ell+1$ of $\alpha$. Thus
$\sigma_\ell(\alpha)$ is vertically trivial in row $\ell+1$ and
horizontally trivial in column $\ell+1$ of $\alpha$.

If $\ell=u$ and $u <v$, then to obtain $\sigma_{\ell}(\alpha)$ from
$\alpha$, we replace row $u$ by the two rows that make row $u$ into
a row of formal vertical composites, and we insert a column of
horizontal identity squares between column $u$ and column $u+1$ of
$\alpha$.

If $\ell=u$ and $u=v$, then to obtain $\sigma_{\ell}(\alpha)$ from
$\alpha$, we replace row $u$ by the two rows that make row $u$ into
a row of formal vertical composites, and we replace column $u$ by
the two columns that make column $u$ into a column of formal
horizontal composites.

If $u<\ell<v$, then to obtain $\sigma_{\ell}(\alpha)$ from $\alpha$,
we replace row $u$ by the row of squares $\beta_1$ in $\bbB$ that
make up the first part of the formal vertical composite row $u$
(consisting partly of region $b$ of Proposition
\ref{pushoutsimplexdescription}), then rows $u+1, u+2, \ldots, \ell$
of $\sigma_{\ell}(\alpha)$ are identity rows, row $\ell+1$ of
$\sigma_{\ell}(\alpha)$ is the composite of the bottom half of row
$u$ of $\alpha$ with rows $u+1, u+2, \ldots, \ell$ of $\alpha$, and
the remaining rows of $\sigma_\ell(\alpha)$ are the remaining rows
of $\alpha$ (shifted down by 1). We also insert a column of
horizontal identity squares between column $\ell$ and column
$\ell+1$ of $\alpha$.

If $u<\ell=v$, then to obtain $\sigma_{\ell}(\alpha)$ from $\alpha$,
we do the row construction as in the case $u<\ell<v$, and we also
replace column $v$ by the two columns that make column $v$ into a
column of formal horizontal composites.

If $u\leq v<\ell$, then to obtain $\sigma_{\ell}(\alpha)$ from
$\alpha$, we do the row construction as in the case $u<\ell<v$, and
we also do the analogous column construction.

The maps $\sigma_\ell$ for $0\leq\ell\leq q$ are compatible with the
boundary operators, $\delta^*N^n(i\overline{r})$, and
$1_{\delta^*N^n\bbQ}$ for the same reason that the analogous maps
associated to a natural transformation of functors are compatible
with the face and degeneracy maps and the functors. Indeed, the
$\sigma_\ell$'s are defined precisely as those for a natural
transformation, we merely take into account the horizontal and
vertical aspects.

In conclusion, we have morphisms of simplicial sets
$$\xymatrix{\delta^*N^n(i)\co \delta^\ast N^n\bbB \ar@{^{(}->}[r] & \delta^\ast N^n \bbQ}$$
$$\xymatrix{\delta^*N^n(\overline{r})\co \delta^\ast N^n \bbQ \ar[r] &  \delta^\ast N^n\bbB}$$
such that $(\delta^*N^n(\overline{r})) \circ
(\delta^*N^n(i))=1_{\delta^\ast N^n\bbB}$ and $(\delta^*N^n(i))
\circ (\delta^*N^n(\overline{r}))$ is simplicially homotopic to
$1_{\delta^*N^n\bbQ}$ via the simplicial homotopy $\sigma$.
\end{pf}

\section{Unit and Counit are Weak Equivalences} \label{unitcounitsection}

In this section we prove that the unit and counit of the adjunction
in (\ref{nfoldcatadjunction}) are weak equivalences. Our main tool
is the $n$-fold Grothendieck construction and the theorem that, in
certain situations, a natural weak equivalence between functors
induces a weak equivalence between the colimits of the functors. We
prove that $N^n$ and the $n$-fold Grothendieck construction are
``homotopy inverses''. From this, we conclude that our Quillen
adjunction (\ref{nfoldcatadjunction}) is actually a Quillen
equivalence. The left and right adjoints of
(\ref{nfoldcatadjunction}) preserve weak equivalences, so the unit
and counit are weak equivalences.

\begin{defn} \label{nfoldGrothendieck}
Let $\xymatrix@1{Y\co(\Delta^{\times n})^{\op} \ar[r] & \mathbf{Set}}$ be a
multisimplicial set. We define the {\it $n$-fold Grothendieck
construction} $\Delta^{\boxtimes n} / Y \in \mathbf{nFoldCat}$ as
follows. The {\it objects} of the $n$-fold category
$\Delta^{\boxtimes n} / Y$ are
$$\Obj \Delta^{\boxtimes n} / Y =\{(y,\overline{k})|\overline{k}=([k_1], \ldots, [k_n]) \in \Delta^{\times n},
y \in Y_{\overline{k}}\}.$$ An {\it $n$-cube} in $\Delta^{\boxtimes
n} / Y$ with $(0,0,\ldots,0)$-vertex $(y,\overline{k})$ and
$(1,1,\ldots,1)$-vertex $(z,\overline{\ell})$ is a morphism
$\xymatrix@1{\overline{f}=(f_1, \ldots, f_n)\co \overline{k} \ar[r]
& \overline{\ell}}$ in $\Delta^{\times n}$ such that
\begin{equation} \label{nfoldGrothendieckmorphism}
\overline{f}^*(z)=y.
\end{equation}
For $\epsilon_\ell\in \{0,1\}$, the
$(\epsilon_1,\epsilon_2,\ldots,\epsilon_n)$-vertex of such an
$n$-cube is
$$(f_1^{1-\epsilon_1},f_2^{1-\epsilon_2}, \ldots, f_n^{1-\epsilon_n})^*(z).$$
For $1\leq i \leq n$, a {\it morphism in direction $i$} is an
$n$-cube that has $f_j$ the identity except at $j=i$. A {\it square
in direction $i i'$} is an $n$-cube such that $f_j$ is the identity
except at $j=i$ and $j=i'$, etc. In this way, the edges, subsquares,
subcubes, etc. of an $n$-cube $\overline{f}$ are determined.
\end{defn}

\begin{examp}
If $n=1$, then the Grothendieck construction of Definition
\ref{nfoldGrothendieck} is the usual Grothendieck construction of a
simplicial set.
\end{examp}

\begin{examp}
The Grothendieck construction $\Delta/\Delta[m]$ of the simplicial
set $\Delta[m]$ is the comma category $\Delta/[m]$.
\end{examp}

\begin{examp}
The Grothendieck construction commutes with external products, that
is, for simplicial sets $X_1,X_2, \ldots, X_n$ we have
$$\Delta^{\boxtimes n}/(X_1\boxtimes X_2 \boxtimes \cdots \boxtimes
X_n)=(\Delta/X_1) \boxtimes (\Delta/X_2) \boxtimes \cdots \boxtimes
(\Delta/X_n).$$
\end{examp}

\begin{rmk} \label{p-multisimplices}
We describe the $n$-fold nerve of the $n$-fold Grothendieck
construction. We learned the $n=1$ case from Chapter 6 of
\cite{joyaltierneysimplicial}. Let $\xymatrix@1{Y\co(\Delta^{\times
n})^{\op} \ar[r] & \mathbf{Set}}$ be a multisimplicial set and
$\overline{p}=([p_1],\ldots,[p_n]) \in \Delta^{\times n}$. Then a
$\overline{p}$-multisimplex of $N^n(\Delta^{\boxtimes n} / Y)$
consists of $n$ composable paths of morphisms in $\Delta$ of lengths
$p_1, p_2, \ldots, p_n$
$$\xymatrix{\langle f^1_1, \dots, f^1_{p_1} \rangle\co [k_0^1] \ar[r]^-{f^1_1} & [k_1^1] \ar[r]^{f^1_2}
& \cdots \ar[r]^{f^1_{p_1}} & [k^1_{p_1}] }$$
$$\xymatrix{\langle f^2_1, \dots, f^2_{p_2} \rangle\co [k_0^2] \ar[r]^-{f^2_1} & [k_1^2] \ar[r]^{f^2_2}
& \cdots \ar[r]^{f^2_{p_2}} & [k^2_{p_2}] }$$
$$\cdots$$
$$\xymatrix{\langle f^n_1, \dots, f^n_{p_n} \rangle\co [k_0^n] \ar[r]^-{f^n_1} & [k_1^n] \ar[r]^{f^n_2}
& \cdots \ar[r]^{f^n_{p_n}} & [k^n_{p_n}] }$$ and a multisimplex
$z$ of $Y$ in degree $$\overline{k_{\overline{p}}}:=(k^1_{p_1},
k^2_{p_2},\ldots,k^n_{p_n}).$$ The last vertex in this
$\overline{p}$-array of $n$-cubes in $\Delta^{\boxtimes n} / Y$ is
$$(z,([k^1_{p_1}], [k^2_{p_2}],\ldots,[k^n_{p_n}])).$$
The other vertices of this array are determined from $z$ by applying
the $f$'s and their composites as in equation
\eqref{nfoldGrothendieckmorphism}. Thus, the set of
$\overline{p}$-multisimplices of $N^n(\Delta^{\boxtimes n} / Y)$ is
\begin{equation} \label{overlinepmultisimplices}
\underset{\langle f^n_1, \dots, f^n_{p_n}
\rangle}{\underset{\cdots}{\underset{\langle f^2_1, \dots, f^2_{p_2}
\rangle}{\underset{\langle f^1_1, \dots, f^1_{p_1}
\rangle}{\coprod}}}}Y_{\overline{k_{\overline{p}}}}.
\end{equation}
\end{rmk}

\begin{prop} \label{NdnGrothendieckpreservescolimits}
The functor $Y \mapsto N^n(\Delta^{\boxtimes n} / Y)$ preserves
colimits.
\end{prop}
\begin{pf}
The set of $\overline{p}$-multisimplices of $N^n(\Delta^{\boxtimes
n} / Y)$ is (\ref{overlinepmultisimplices}). The assignment of $Y$
to the expression in (\ref{overlinepmultisimplices}) preserves
colimits.
\end{pf}

\begin{rmk}
We can also describe the $p$-simplices of
$\delta^*N^n(\Delta^{\boxtimes n} / Y)$. We learned the $n=1$ case
from Joyal and Tierney in Chapter 6 of
\cite{joyaltierneysimplicial}. A $p$-simplex of $\delta^*
N^n(\Delta^{\boxtimes n} / Y )$ is a composable path of $p$
$n$-cubes
$$\xymatrix{\overline{f^i}\co (y^{i-1},\overline{k^{i-1}}) \ar[r] & (y^{i},\overline{k^i})}$$
($i=1,\ldots,p$). Each $y^i$ is determined from $y^p$ by the
$\overline{f^i}$'s, as in equation
\eqref{nfoldGrothendieckmorphism}. The last target, namely
$(y^p,\overline{k^p})$, is the same as a morphism of multisimplicial
sets $\xymatrix@1{\Delta^{\times n}[\overline{k^i}] \ar[r] & Y}$. So
by Yoneda, a $p$-simplex is the same as a composable path of
morphisms of multisimplicial sets
$$\xymatrix{\Delta^{\times n}[\overline{k^0}] \ar[r] & \Delta^{\times n}[\overline{k^1}] \ar[r]
& \cdots \ar[r] & \Delta^{\times n}[\overline{k^p}] \ar[r] & Y }.$$
The set of $p$-simplices of $\delta^* N^n(\Delta^{\boxtimes n} / Y
)$ is
\begin{equation} \label{psimplicesofGrothendieck}
\coprod_{\Delta^{\times n}[\overline{k^0}] \rightarrow
\Delta^{\times n}[\overline{k^1}] \rightarrow \cdots \rightarrow
\Delta^{\times n}[\overline{k^p}]} Y_{\overline{k^p}}.
\end{equation}
\end{rmk}

Let us recall the natural morphism of simplicial sets $\xymatrix@1{N
(\Delta/ X) \ar[r] & X}$ in 6.1 of \cite{joyaltierneysimplicial},
which we shall call $\rho_X$ as in Appendix A of
\cite{moerdijksvenssonOnshapiro}. First note that any path of
morphisms in $\Delta$
\begin{equation} \label{pathinDelta}
\xymatrix{ [k_0] \ar[r] & [k_1] \ar[r] & \cdots \ar[r] & [k_p] }
\end{equation}
determines a morphism
\begin{equation} \label{imagemap}
\begin{array}{c}
\xymatrix{[p] \ar[r] & [k_p]} \\
i \mapsto \im{k_i}
\end{array}
\end{equation}
where $\im{k_i}$ refers to the image of $k_i$ under the composite of
the last $p-i$ morphisms in (\ref{pathinDelta}). Note also that
paths of the form (\ref{pathinDelta}) are in bijective
correspondence with paths of the form
\begin{equation}
\xymatrix{ \Delta[k_0] \ar[r] & \Delta[k_1] \ar[r] & \cdots \ar[r] &
\Delta[k_p] }
\end{equation}
by the Yoneda Lemma. \label{rhoXdefinition} The morphism
$\xymatrix@1{\rho_X\co N (\Delta/ X) \ar[r] & X}$ sends a
$p$-simplex
$$\xymatrix{\Delta[k^0] \ar[r] & \Delta[k^1] \ar[r]
& \cdots \ar[r] & \Delta[k^p] \ar[r] & X }$$ to the composite
$$\xymatrix{\Delta[p] \ar[r] & \Delta[k^p] \ar[r] & X}$$
where the first morphism is the image of (\ref{imagemap}) under the
Yoneda embedding. As is well known, the morphism $\xymatrix@1{N
(\Delta/ X) \ar[r] & X}$ is a natural weak equivalence (see Theorem
6.2.2 of \cite{joyaltierneysimplicial}, page 21 of \cite{illusieII},
page 359 of \cite{waldhausen}).

We analogously define a morphism of multisimplicial sets
$$\xymatrix{\rho_Y\co N^n(\Delta^{\boxtimes n} / Y) \ar[r] & Y}$$
natural in $Y$. Consider a $\overline{p}$-multisimplex of
$N^n(\Delta^{\boxtimes n} / Y)$ as in Remark \ref{p-multisimplices}.
For each $1 \leq j \leq n$, the path $\langle f^j_1, \dots,
f^j_{p_j} \rangle$ gives rise to a morphism in $\Delta$
$$\xymatrix{[p_j] \ar[r] & [k^j_{p_j}]}$$
as in (\ref{pathinDelta}) and (\ref{imagemap}). Together these form
a morphism in $\Delta^{\times n}$, which induces a morphism of
multisimplicial sets
$$\xymatrix{\Delta^{\times n}[\overline{p}] \ar[r] & \Delta^{\times n}[\overline{k_{\overline{p}}}]}.$$
The morphism $\rho_Y$ assigns to the $\overline{p}$-multisimplex we
are considering the $\overline{p}$-multisimplex
$$\xymatrix{\Delta^{\times n}[\overline{p}] \ar[r] & \Delta^{\times n}[\overline{k_{\overline{p}}}] \ar[r]^-z & Y}.$$
This completes the definition of the natural transformation $\rho$.

\begin{rmk} \label{rhowithexternalproducts}
The natural transformation $\rho$ is compatible with external
products. If $X_1,X_2,\dots,X_n$ are simplicial sets and
$Y=X_1\boxtimes X_2\boxtimes \cdots \boxtimes X_n$, then
$$\xymatrix{\rho_{Y}
\co N^n(\Delta^{\boxtimes n}/Y) \ar[r] & Y}$$ is equal to
$$\rho_{X_1} \boxtimes \rho_{X_2} \boxtimes \cdots \boxtimes \rho_{X_n}
\co$$
$$\xymatrix{ N(\Delta /X_1) \boxtimes N(\Delta /X_2) \boxtimes \cdots \boxtimes N(\Delta /X_n) \ar[r]
& X_1 \boxtimes X_2 \boxtimes \cdots \boxtimes X_n.}$$ Thus
$\delta^*\rho_Y=\rho_{X_1}\times \rho_{X_2} \times \cdots \times
\rho_{X_n}$ is a weak equivalence, since in $\mathbf{SSet}$ any
finite product of weak equivalences is a weak equivalence. We
conclude that $\rho_Y$ is a weak equivalence of multisimplicial sets
whenever $Y$ is an external product. (For us, a morphism $f$ of
multisimplicial sets is a {\it weak equivalence} if and only if
$\delta^*f$ is a weak equivalence of simplicial sets.) As we shall
soon see, $\rho_Y$ is a weak equivalence for all $Y$.
\end{rmk}

We quickly recall what we will need regarding Reedy model
structures. The following definition and proposition are part of
Definitions 5.1.2, 5.2.2, and Theorem 5.2.5 of \cite{hovey}, or
Definitions 15.2.3, 15.2.5, and Theorem 15.3.4 of \cite{hirschhorn}
\begin{defn}
Let $(\mathcal{B},\mathcal{B}_+,\mathcal{B}_-)$ be a Reedy category
and $\mathcal{C}$ a category with all small colimits and limits. For
$i \in \mathcal{B}$, the {\it latching category} $\mathcal{B}_i$ is
the full subcategory of $\mathcal{B}_+/i$ on the {\it non-identity}
morphisms $\xymatrix@1{b \ar[r] & i}$. For $F \in
\mathcal{C}^\mathcal{B}$ the {\it latching object of $F$ at $i$} is
the colimit $L_iF$ of the composite functor
\begin{equation} \label{latchingobjectequation}
\xymatrix{\mathcal{B}_i \ar[r] & \mathcal{B} \ar[r]^F & \mathcal{C}
}.
\end{equation}
For $i \in \mathcal{B}$, the {\it matching category} $\mathcal{B}^i$
is the full subcategory of $i/\mathcal{B}_-$ on the {\it
non-identity} morphisms $\xymatrix@1{i \ar[r] & b}$. For $F \in
\mathcal{C}^\mathcal{B}$ the {\it matching object of $F$ at $i$} is
the limit $M_iF$ of the composite functor
\begin{equation} \label{matchingobjectequation}
\xymatrix{\mathcal{B}^i \ar[r] & \mathcal{B} \ar[r]^F & \mathcal{C}
}.
\end{equation}
\end{defn}

\begin{thm}[Kan]
Let $(\mathcal{B},\mathcal{B}_+,\mathcal{B}_-)$ be a Reedy category
and $\mathcal{C}$ a model category. Then the levelwise weak
equivalences, Reedy fibrations, and Reedy cofibrations form a model
structure on the category $\mathcal{C}^\mathcal{B}$ of functors
$\xymatrix@1{\mathcal{B} \ar[r] & \mathcal{C}}$.
\end{thm}

\begin{rmk} \label{remarkconsequence}
A consequence of the definitions is that a functor
$\xymatrix@1{\mathcal{B} \ar[r] & \mathcal{C}}$ is {\it Reedy
cofibrant} if and only if the induced morphism $\xymatrix{L_iF
\ar[r] & Fi}$ is a cofibration in $\mathcal{C}$ for all objects $i$
of $\mathcal{B}$.
\end{rmk}

\begin{prop}[Compare Example 15.1.19 of \cite{hirschhorn}]
The category of multisimplices
$$\Delta^{\times n}Y:=\Delta^{\times n}/Y$$ of a multisimplicial
set $\xymatrix@1{Y\co (\Delta^{\times n})^{\op} \ar[r] & \mathbf{Set}}$
is a Reedy category. The degree of a
$\overline{p}$-multisimplex is $p_1+p_2\cdots+p_n$. The direct
subcategory $(\Delta^{\times n}Y)_+$ consists of those morphisms $(f_1,\dots, f_n)$ that are
iterated coface maps in each coordinate, \ie injective maps in each coordinate.
The inverse subcategory $(\Delta^{\times n}Y)_-$
consists of those morphisms $(f_1,\dots, f_n)$ that are iterated
codegeneracy maps in each coordinate, \ie surjective maps in each coordinate.
\end{prop}

\begin{prop}[Compare Proposition 15.10.4(1) of \cite{hirschhorn}] \label{multisimpliceshavefibrantconstants}
If $\mathcal{B}$ is the category of multisimplices of a
multisimplicial set, then for every $i \in \mathcal{B}$, the
matching category $\mathcal{B}^i$ is either connected or empty.
\end{prop}
\begin{pf}
This follows from the multidimensional Eilenberg-Zilber Lemma, recalled in Proposition \ref{prop:EZmultsimplicial}.
Let $\xymatrix@1{Y\co(\Delta^{\times n})^{\op} \ar[r] & \mathbf{Set}}$ be a multipsimplicial set
and $\mathcal{B}=\Delta^{\times n} Y$ its category of multisimplices.

Let $\xymatrix@1{i\co \Delta^{\times n}\left[ \overline{p} \right] \ar[r] & Y }$ be a degenerate multisimplex. Then there
exists a non-trivial, componentwise surjective map $\overline{\tau}$ and a totally non-degenerate multisimplex $t$ with $i=(\overline{\tau})^*t$.  The pair $(\overline{\tau},t)$ is an object of the matching category $\mathcal{B}^i$. If $(\overline{\eta},b)$
is another object of $\mathcal{B}^i$, there exists a componentwise surjective map $\overline{g}$ and a totally non-degenerate
$b' \in \mathcal{B}$ such that $b=(\overline{g})^*b'$. But $i=(\overline{\eta})^*b=(\overline{\eta})^*(\overline{g})^*b'$ implies
that $b'=t$, $\overline{g} \circ \overline{\eta}=\overline{\tau}$, and $\overline{g}$ is a morphism in $\mathcal{B}^i$ from
$(\overline{\eta},b)$ to $(\overline{\tau},t)$. Thus, whenever $i$ is degenerate, there is a morphism from any object of $\mathcal{B}^i$
to $(\overline{\tau},t)$ and $\mathcal{B}^i$ is connected. One can also show $(\overline{\tau},t)$ is a terminal object of $\mathcal{B}^i$, but we do not need this.

Let $\xymatrix@1{i\co \Delta^{\times n}\left[ \overline{p} \right] \ar[r] & Y }$ be a totally non-degenerate multisimplex.
An object of the matching category $\mathcal{B}^i$ is a non-trivial,
componentwise surjective map $\overline{\eta}$ and a multisimplex $b$ with $i=(\overline{\eta})^*b$.
Such $\overline{\eta}$ and $b$ cannot exist because $i$ is totally non-degenerate.
Thus, whenever $i$ is totally non-degenerate, the matching category $\mathcal{B}^i$ is empty .
\end{pf}

\begin{thm} \label{colimitpreservesweakequivalences}
Suppose $\mathcal{C}$ is a model category and $\mathcal{B}$ is a
Reedy category such that for all $i \in \mathcal{B}$, the matching
category $\mathcal{B}^i$ is either connected or empty. Then the
colimit functor
$$\xymatrix{\text{\rm colim}\co \mathcal{C}^\mathcal{B} \ar[r] &
\mathcal{C}}$$ takes levelwise weak equivalences between Reedy
cofibrant functors to weak equivalences between cofibrant objects of
$\mathcal{C}$.
\end{thm}
\begin{pf}
This is merely a summary of Definition 15.10.1(2), Proposition
15.10.2(2), and Theorem 15.10.9(2) of \cite{hirschhorn}.
\end{pf}

\begin{notation} \label{BC}
Let $\xymatrix@1{Y\co (\Delta^{\times n})^{\op} \ar[r] & \mathbf{Set}}$ be a
multisimplicial set, $\mathcal{B}=\Delta^{\times n}Y$,
$\mathcal{C}=\mathbf{SSet}$, and $\xymatrix@1{i \co \Delta^{\times
n}[\overline{m}] \ar[r] & Y}$ an object of $\mathcal{B}$. Then the
set of nonidentity morphisms in $\mathcal{B}_+$ with target $i$ is
the set of morphisms $(f_1, \dots, f_n)$ in $\Delta^{\times n}$ with
target $[\overline{m}]$ such that each $f_j$ is injective and not
all $f_j$'s are the identity.
\end{notation}
\begin{notation} \label{twofunctors}
Let $F$ and $G$ be the following two functors.
$$\xymatrix{F\co\Delta^{\times n} Y \ar[r] & \mathbf{SSet^n}}$$
$$\left[\Delta^{\times n}[\overline{m}] \rightarrow Y \right] \mapsto
N^n(\Delta^{\boxtimes n}/ \Delta^{\times n}[\overline{m}])$$
$$\xymatrix{G\co\Delta^{\times n} Y \ar[r] & \mathbf{SSet^n}}$$
$$\left[\Delta^{\times n}[\overline{m}] \rightarrow Y \right] \mapsto
\Delta^{\times n}[\overline{m}]$$ Note that $\delta^\ast \circ F$
and $\delta^\ast \circ G$ are in $\mathcal{C}^\mathcal{B}$. The
natural transformation $\rho$ induces a natural transformation we
denote by
$$\xymatrix{\rho^Y\co F \ar@{=>}[r] & G}.$$
\end{notation}

\begin{rmk} \label{levelwiseweakequivalence}
The natural transformation $\rho^Y$ is levelwise a weak equivalence
by Remark \ref{rhowithexternalproducts}.
\end{rmk}

\begin{lem} \label{colimrho^Y=rho_Y}
The morphism in $\mathbf{SSet^n}$
\begin{equation*}
\underset{\Delta^{\times n} Y}{\colim}\rho^Y\co
\underset{\Delta^{\times n} Y}{\colim} F \xymatrix{ \ar[r]
&}\underset{\Delta^{\times n} Y}{\colim} G
\end{equation*}
is equal to $$\xymatrix{\rho_Y\co N^n(\Delta^{\boxtimes n}/Y) \ar[r]
& Y.}$$
\end{lem}
\begin{pf}
By Proposition \ref{NdnGrothendieckpreservescolimits}, we have
$$\aligned \underset{\Delta^{\times n} Y}{\colim} F
&=\underset{\Delta^{\times n}[\overline{m}] \rightarrow
Y}{\colim}N^n(\Delta^{\boxtimes n}/ \Delta^{\times n}[\overline{m}])
\\
&=N^n(\Delta^{\boxtimes n}/(\underset{\Delta^{\times
n}[\overline{m}] \rightarrow Y}{\colim} \Delta^{\times
n}[\overline{m}])) \\
&=N^n(\Delta^{\boxtimes n}/Y).
\endaligned$$
\end{pf}

\begin{lem} \label{delta*Fcofibrant}
The functor
$$\xymatrix{\delta^*\circ F\co\Delta^{\times n} Y \ar[r] & \mathbf{SSet}}$$
$$\left[\Delta^{\times n}[\overline{m}] \rightarrow Y \right] \mapsto
N(\Delta/ \Delta[m_1]) \times N(\Delta/  \Delta[m_2]) \times \cdots
\times N(\Delta/ \Delta[m_n])$$ is Reedy cofibrant.
\end{lem}
\begin{pf}
We use Notations \ref{BC} and \ref{twofunctors}. The colimit of
equation \eqref{latchingobjectequation} is
$$L_i(\delta^\ast \circ F)=\underset{1 \leq j \leq n}{\bigcup}
N(\Delta/\Delta[m_1]) \times \cdots \times N(\Delta/\partial
\Delta[m_j]) \times \cdots \times N(\Delta/\Delta[m_n])$$ and
$\delta^\ast \circ F (i)=N(\Delta/  \Delta[m_2]) \times \cdots
\times N(\Delta/ \Delta[m_n]).$ The map
$$\xymatrix{L_i(\delta^\ast \circ F) \ar[r] & \delta^\ast \circ
F(i)}$$ is injective, or equivalently, a cofibration. Remark
\ref{remarkconsequence} now implies that $\delta^\ast \circ F$ is
Reedy cofibrant.
\end{pf}

\begin{lem} \label{delta*Gcofibrant}
The functor
$$\xymatrix{\delta^* \circ G\co\Delta^{\times n} Y \ar[r] & \mathbf{SSet}}$$
$$\left[\Delta^{\times n}[\overline{m}] \rightarrow Y \right] \mapsto
\Delta[m_1] \times \Delta[m_2] \times \cdots \times \Delta[m_n]$$ is
Reedy cofibrant.
\end{lem}
\begin{pf}
We use Notations \ref{BC} and \ref{twofunctors}. The colimit of
equation \eqref{latchingobjectequation} is
$$L_i(\delta^\ast \circ G)=\underset{1 \leq j \leq n}{\bigcup} \Delta[m_1]
\times \cdots \times \partial \Delta[m_j] \times \cdots \times
\Delta[m_n]$$ and $\delta^\ast \circ G (i)=\Delta[m_1] \times
\Delta[m_2] \times \cdots \times \Delta[m_n].$ The morphism
$$\xymatrix{L_i(\delta^\ast \circ G) \ar[r] & \delta^\ast \circ
G(i)}$$ is injective, or equivalently, a cofibration. Remark
\ref{remarkconsequence} now implies that $\delta^\ast \circ G$ is
Reedy cofibrant.
\end{pf}

\begin{thm} \label{rhowe}
For every multisimplicial set $\xymatrix@1{Y\co (\Delta^{\times n})^{\op}
\ar[r] & \mathbf{Set}}$, the morphism
$$\xymatrix{\rho_Y\co N^n(\Delta^{\times n}/Y) \ar[r] &
Y}$$ is a weak equivalence of multisimplicial sets.
\end{thm}
\begin{pf}
Fix a multisimplicial set $Y$, and let $F$, $G$, and $\rho^Y$ be as
in Notation \ref{twofunctors}. The natural transformation
$\xymatrix@1{\delta^* \rho^Y \co \delta^* F \ar@{=>}[r] & \delta^*
G}$ is levelwise a weak equivalence of simplicial sets by Remark
\ref{levelwiseweakequivalence}, and is a natural transformation
between Reedy cofibrant functors by Lemmas \ref{delta*Fcofibrant}
and \ref{delta*Gcofibrant}. By Proposition
\ref{multisimpliceshavefibrantconstants}, each matching category of
the Reedy category $\Delta^{\times n} Y$ is connected or empty.
Theorem \ref{colimitpreservesweakequivalences} then guarantees that
the morphism
\begin{equation*}
\underset{\Delta^{\times n} Y}{\colim}\delta^*\rho^Y\co
\underset{\Delta^{\times n} Y}{\colim} \delta^*\circ F \xymatrix{
\ar[r] &}\underset{\Delta^{\times n} Y}{\colim} \delta^*\circ G
\end{equation*}
is a weak equivalence of simplicial sets. Since $\delta^*$ is a left
adjoint, it commutes with colimits, and we have
$$\underset{\Delta^{\times n} Y}{\colim}\delta^*\rho^Y
=\delta^*\underset{\Delta^{\times n}
Y}{\colim}\rho^Y=\delta^*\rho_Y$$ by Lemma \ref{colimrho^Y=rho_Y}.
We conclude $\delta^*\rho_Y$ is a weak equivalence, and that
$\rho_Y$ is a weak equivalence of multisimplicial sets.
\end{pf}

We also define an $n$-fold functor
$$\xymatrix{\lambda_\bbD \co \Delta^{\boxtimes n} / N^n(\bbD) \ar[r] &  \bbD}$$
natural in $\bbD$, by analogy to Appendix A of
\cite{moerdijksvenssonOnshapiro}, and many others. If
$(y,\overline{k})$ is an object of $\Delta^{\boxtimes n} /
N^n(\bbD)$, then $\lambda(y, \overline{k})$ is the $n$-fold category in the last
vertex of the array of $n$-cubes $y$, namely
$$\lambda_\bbD(y,\overline{k})=y_{\overline{k}}.$$

\begin{thm}\label{lambdawe}
For any $n$-fold category $\bbD$, we have
$N^n(\lambda_\bbD)=\rho_{N^n(\bbD)}$. In particular,
$\lambda_\bbD$ is a weak equivalence of $n$-fold categories.
\end{thm}

\begin{cor} \label{Ndinducesequivalence}
The functor $\xymatrix@1{N^n\co \mathbf{nFoldCat} \ar[r] &
\mathbf{SSet^n}}$ induces an equivalence of categories
$$\Ho\mathbf{nFoldCat} \simeq \Ho \mathbf{SSet^n}.$$
Here $\text{\rm Ho}$ refers to the category obtained by formally
inverting weak equivalences. There is no reference to any model
structure.
\end{cor}
\begin{pf}
An ``inverse'' to $N^n$ is the $n$-fold Grothendieck construction,
since $\rho$ and $\lambda$ induce natural isomorphisms after passing
to homotopy categories by Theorems \ref{rhowe} and \ref{lambdawe}.
\end{pf}

The following simple proposition, pointed out to us by Denis-Charles
Cisinski, will be of use.
\begin{prop} \label{prop:cisinski}
Let $\xymatrix@C=4pc{{\bf C} \ar@{}[r]|{\perp} \ar@/^1pc/[r]^{F} &
\ar@/^1pc/[l]^{G} {\bf D}}$ be a Quillen equivalence. If both $F$ and $G$ preserve
weak equivalences, then
\begin{enumerate}
\item \label{item:FG_detect}
Both $F$ and $G$ detect weak equivalences,
\item \label{item:unit_counit_wes}
The unit and counit of the adjunction $F \dashv G$ are weak
equivalences.
\end{enumerate}
\end{prop}
\begin{pf}
\noindent \ref{item:FG_detect}
We prove $F$ detects weak equivalences; the proof that $G$ detects weak equivalences is similar.
Let $\xymatrix@1{Q\co \bfC \ar[r] & \bfC}$ be a cofibrant replacement functor on $\bfC$, that is, $QC$ is cofibrant for all objects
$C$ in $\bfC$ and there is a natural acyclic fibration $\xymatrix@1{q\co QC \ar[r] & C}$. Suppose $Ff$ is a weak equivalence. Then
$FQf$ is a weak equivalence (apply $F$ to the naturality diagram for $f$ and $Q$ and use the 3-for-2 property).
The total left derived functor $\bfL F$ is
the composite
$$\xymatrix@C=3pc{\Ho \bfC \ar[r]_{\Ho Q} \ar@/^1pc/[rr]^{\bfL F} & \Ho \bfC_c \ar[r]_{\Ho F\vert_{\bfC_c}} & \Ho \bfD },$$
where $\bfC_c$ is the full subcategory of $\bfC$ on the cofibrant objects of $\bfC$. Then $\bfL F [f]$ is an isomorphism in $\Ho \bfD$,
as $FQf$ is a weak equivalence in $\bfD$.
The functor $\bfL F$ detects isomorphisms, as it is an equivalence of categories, so $[f]$ is an isomorphism in $\Ho \bfC$. Finally,
a morphism in $\bfC$ is a weak equivalence if and only if its image in $\Ho \bfC$ is an isomorphism, so $f$ is a weak equivalence in
$\bfC$, and $F$ detects weak equivalences.

\noindent \ref{item:unit_counit_wes} We prove that the unit of the
adjunction $F\dashv G$ is a natural weak equivalence; the proof that
the counit is a natural weak equivalence is similar. Let
$\xymatrix@1{Q\co \bfC \ar[r] & \bfC}$ be a cofibrant replacement
functor on $\bfC$, that is, $QC$ is cofibrant for every object $C$
in $\bfC$ and there is a natural acyclic fibration
$\xymatrix@1{q_C\co QC \ar[r] & C}$. Let $\xymatrix@1{R \co \bfD
\ar[r] & \bfD}$ be a fibrant replacement functor on $\bfD$, that is,
$RD$ is fibrant for every object $D$ in $\bfD$ and there is a
natural acyclic cofibration $\xymatrix@1{r_D\co D \ar[r] & RD}$.
Since $F\dashv G$ is a Quillen equivalence, the composite
$$\xymatrix@C=4pc{QC \ar[r]^-{\eta_{QC}} & GFQX \ar[r]^-{Gr_{FQX}} & GRFQX}$$
is a weak equivalence by Proposition 1.3.13 of \cite{hovey}. Then $\eta_{QC}$ is a weak equivalence by the 3-for-2 property and the
hypothesis that $G$ preserves weak equivalences. An application of 3-for-2 to the naturality diagram for $\eta$
$$\xymatrix@C=3pc{QC \ar[r]^-{\eta_{QC}} \ar[d]_{q_C} & GFQC \ar[d]^{GFq_C} \\ C \ar[r]_-{\eta_C} & GFC }$$ shows that
$\eta_C$ is a weak equivalence (recall $GF$ preserves weak equivalences).
\end{pf}

\begin{lem} \label{lem:HoG_equiv_implies_right_derived_equiv}
Let $\xymatrix@1{G\co \bfD \ar[r] & \bfC}$ be a right Quillen
functor. Suppose $\xymatrix@1{\Ho G \co \Ho \bfD \ar[r] & \Ho \bfC}$
is an equivalence of categories. Then the total right derived
functor
$$\xymatrix@C=3pc{\Ho \bfD \ar[r]_{\Ho R} \ar@/^1pc/[rr]^{\bfR G}
& \Ho \bfD_f \ar[r]_{\Ho G\vert_{\bfD_f}} & \Ho \bfC }$$ is an
equivalence of categories. Here $R$ is a fibrant replacement functor
on $\bfD$, and $\bfD_f$ is the full subcategory of $\bfD$ on the
fibrant objects.
\end{lem}
\begin{pf}
The functors $\xymatrix@1@C=3pc{\Ho \bfD \ar@<.5ex>[r]^{\Ho R} &
\ar@<.5ex>[l]^{\Ho i} \Ho \bfD_f}$ are equivalences of categories,
``inverse'' to one another. Then $\Ho G\vert_{\bfD_f}=(\Ho G) \circ
(\Ho i)$ is a composite of equivalences.
\end{pf}

\begin{lem} \label{lem:right_derived_equiv_implies_Quillen_equiv}
Suppose $L \dashv R$ is an adjunction and $R$ is an equivalence of
categories. Then the unit $\eta$ and counit $\varepsilon$ of this
adjunction are natural isomorphisms.
\end{lem}
\begin{pf}
By Theorem 1 on page 93 of \cite{maclaneworking}, $R$ is part of an
adjoint equivalence $L' \dashv R$ with unit $\eta'$ and counit
$\varepsilon'$. By the universality of $\eta$ and $\eta'$ there
exists an isomorphism $\xymatrix@1{\theta_X\co LX  \ar[r] & L'X}$
such that $(R\theta_X) \circ \eta_X=\eta_X'$. Since $\eta_X'$ is
also an isomorphism, we see that $\eta_X$ is an isomorphism. A
similar argument shows that the counit $\varepsilon$ is a natural
isomorphism.
\end{pf}

\begin{prop} \label{unitcounitwe}
The unit and counit of (\ref{nfoldcatadjunction})
\begin{equation*}
\xymatrix@C=4pc{\mathbf{SSet} \ar@{}[r]|{\perp}
\ar@/^1pc/[r]^-{\Sd^2} &  \ar@/^1pc/[l]^-{\Ex^2} \mathbf{SSet}
\ar@{}[r]|{\perp} \ar@/^1pc/[r]^-{\delta_!} &
\ar@/^1pc/[l]^-{\delta^\ast} \mathbf{SSet^n} \ar@{}[r]|{\perp}
\ar@/^1pc/[r]^-{c^n} & \ar@/^1pc/[l]^-{N^n} \mathbf{nFoldCat}}
\end{equation*}
are weak equivalences.
\end{prop}
\begin{pf}
Let $F\dashv G$ denote the adjunction in (\ref{nfoldcatadjunction}).
This is a Quillen adjunction by Theorem \ref{MainModelStructure}. We
first prove it is even a Quillen equivalence. The functor $\Ex^2
\delta^*$ is known to induce an equivalence of homotopy categories,
and $N^n$ induces an equivalence of homotopy categories by Corollary
\ref{Ndinducesequivalence}, so $G=\Ex^2 \delta^* N^n$ induces an
equivalence of homotopy categories $\Ho G$. Lemma
\ref{lem:HoG_equiv_implies_right_derived_equiv} then says that the
total right derived functor $\bfR G$ is an equivalence of
categories. The derived adjunction $\bfL F \dashv \bfR G$ is then an
adjoint equivalence by Lemma
\ref{lem:right_derived_equiv_implies_Quillen_equiv}, so $F \dashv G$
is a Quillen equivalence.

By Ken Brown's Lemma, the left Quillen functor $F$ preserves weak
equivalences (every simplicial set is cofibrant). The right Quillen
functor $G$ preserves weak equivalences by definition. Proposition
\ref{prop:cisinski} now guarantees that the unit and counit are weak
equivalences.
\end{pf}

We now summarize our main results of Theorem
\ref{MainModelStructure}, Corollary \ref{Ndinducesequivalence},
Proposition \ref{unitcounitwe}.

\begin{thm} \label{maintheoremsummary}
\begin{enumerate}
\item
There is a cofibrantly generated model structure on
$\mathbf{nFoldCat}$ such that an $n$-fold functor $F$ is a weak
equivalence (respectively fibration) if and only if $\Ex^2 \delta^*
N^n(F)$ is a weak equivalence (respectively fibration). In
particular, an $n$-fold functor is a weak equivalence if and only if
the diagonal of its nerve is a weak equivalence of simplicial sets.
\item
The adjunction
\begin{equation*}
\xymatrix@C=4pc{\mathbf{SSet} \ar@{}[r]|{\perp}
\ar@/^1pc/[r]^-{\Sd^2} &  \ar@/^1pc/[l]^-{\Ex^2} \mathbf{SSet}
\ar@{}[r]|{\perp} \ar@/^1pc/[r]^-{\delta_!} &
\ar@/^1pc/[l]^-{\delta^\ast} \mathbf{SSet^n} \ar@{}[r]|{\perp}
\ar@/^1pc/[r]^-{c^n} & \ar@/^1pc/[l]^-{N^n} \mathbf{nFoldCat}}
\end{equation*}
is a Quillen equivalence.
\item
The unit and counit of this Quillen equivalence are weak
equivalences.
\end{enumerate}
\end{thm}

\begin{cor} \label{maincorollary}
The homotopy category of $n$-fold categories is equivalent to the
homotopy category of topological spaces.
\end{cor}

Another approach to proving that $N^n$ and the $n$-fold Grothendieck
construction are homotopy inverse would be to apply a
multisimplicial version of the following Weak Equivalence Extension
Theorem of Joyal-Tierney. We apply the present Weak Equivalence
Extension Theorem to prove that there is a natural isomorphism
$$\xymatrix@1{\delta^* N^n(\Delta^{\boxtimes n} / \delta_!\text{-}
) \ar@{=>}[r] & 1_{\text{\rm Ho}\;\mathbf{SSet}}}.$$

\begin{thm}[Theorem 6.2.1 of \cite{joyaltierneysimplicial}]
\label{weakequivalenceextension} Let $\xymatrix@1{\phi \co F
\ar@{=>}[r] & G}$ be a natural transformation between functors
$\xymatrix@1{F,G\co \Delta \ar[r] & \mathbf{SSet}}$. We denote by
$\xymatrix@1{\phi^+\co F^+ \ar@{=>}[r] & G^+}$ the left Kan
extension along the Yoneda embedding $\xymatrix@1{Y\co \Delta \ar[r]
& \mathbf{SSet}}$.
$$\xymatrix{\mathbf{SSet} \ar[dr]^{F^+,G^+} & \\ \Delta \ar[r]_-{F,G} \ar[u]^Y & \mathbf{SSet}}$$
Suppose that $G$ satisfies the following condition.
\begin{itemize}
\item
$\im G\epsilon^0 \cap \im G \epsilon^1 = \emptyset,$ where
$\xymatrix@1{\epsilon^i\co [0] \ar[r] & [1]}$ is the injection which
misses $i$.
\end{itemize}
If $\xymatrix@1{\phi[m]\co F[m] \ar[r] & G[m]}$ is a weak
equivalence for all $m \geq 0$, then $$\xymatrix@1{\phi^+X\co F^+ X
\ar[r] & G^+X }$$ is a weak equivalence for every simplicial set
$X$.
\end{thm}

\begin{lem} \label{Grothendieckpreservescolimits}
The functor
$$\xymatrix{\mathbf{SSet^n} \ar[r] & \mathbf{SSet}}$$
$$\xymatrix{Y \mapsto \delta^* N^n(\Delta^{\boxtimes n} / Y )}$$
preserves colimits.
\end{lem}
\begin{pf}
The functor which assigns to $Y$ the expression in
(\ref{psimplicesofGrothendieck}) is colimit preserving.
\end{pf}

\begin{prop} \label{zigzag1}
For every simplicial set $X$, the canonical morphism
$$\xymatrix@1{\delta^* N^n(\Delta^{\boxtimes n} / \delta_!X ) \ar[r] &
\delta^*\delta_!X}$$ is a weak equivalence.
\end{prop}
\begin{pf}
We apply the Weak Equivalence Extension Theorem
\ref{weakequivalenceextension}. Let $\xymatrix@1{F,G\co \Delta
\ar[r] & \mathbf{SSet}}$ be defined by
$$F[m]=\delta^* N^n(\Delta^{\boxtimes n} / \delta_!\Delta[m] )$$
$$G[m]=\delta^*\delta_!\Delta[m].$$
The functor $$\xymatrix{\delta^* N^n(\Delta^{\boxtimes n} /
\delta_!\text{-})\co \mathbf{SSet} \ar[r] & \mathbf{SSet}}$$
preserves colimits by Lemma \ref{Grothendieckpreservescolimits} and
the fact that $\delta_!$ is a left adjoint. The functor
$$\xymatrix{\delta^* \delta_!\co \mathbf{SSet} \ar[r] & \mathbf{SSet}}$$
preserves colimits since both $\delta^*$ and $\delta_!$ are both
left adjoints. Thus the canonical comparison morphisms
$$\xymatrix{F^+X \ar[r] & \delta^* N^n(\Delta^{\boxtimes n} / \delta_!X )}$$
$$\xymatrix{G^+X \ar[r] & \delta^*\delta_!X}$$
are isomorphisms.

The condition on $G$ listed in Theorem
\ref{weakequivalenceextension} is easy to verify, since
$$\xymatrix{G\epsilon^0=\epsilon^0\times \epsilon^0 \co \Delta[0] \times \Delta[0] \ar[r] & \Delta[1] \times \Delta[1]}$$
$$\xymatrix{G\epsilon^1=\epsilon^1\times \epsilon^1 \co \Delta[0] \times \Delta[0] \ar[r] & \Delta[1] \times \Delta[1]}.$$

All that remains is to define natural morphisms
$$\xymatrix{\phi[m]\co \delta^* N^n(\Delta^{\boxtimes n} /
\Delta[m,\ldots,m] ) \ar[r] & \Delta[m] \times \cdots \times
\Delta[m]}$$ and to show that each is a weak equivalence of
simplicial sets. By the description in Definition
\ref{nfoldGrothendieck}, an object of $\Delta^{\boxtimes n} /
\Delta[m,\ldots,m]$ is a morphism
$$\xymatrix{y=(y_1,\ldots,y_n)\co \overline{k} \ar[r] & ([m],\ldots,[m])}$$
in $\Delta^{\times n}$. An $n$-cube $\overline{f}$ is a morphism in
$\Delta^{\times n}$ making the diagram
$$\xymatrix{\overline{k} \ar[rr]^{\overline{f}} \ar[dr]_{y} & &  \overline{k'} \ar[dl]^{y'}
\\ & ([m], \ldots, [m]) &}$$
commute. A $p$-simplex in $\delta^* N^n(\Delta^{\boxtimes n} /
\Delta[m,\ldots,m] )$ is a path
$\overline{f^1},\ldots,\overline{f^p}$ of composable morphisms in
$\Delta^{\times n}$ making the appropriate triangles commute. We see
that
$$\delta^* N^n(\Delta^{\boxtimes n} /
\Delta[m,\ldots,m] ) \cong N(\Delta/\Delta[m]) \times \cdots
N(\Delta/\Delta[m]).$$

We define $\phi[m]$ to be the product of $n$-copies of the weak
equivalence $$\xymatrix{\rho_{\Delta[m]} \co N(\Delta/\Delta[m])
\ar[r] & \Delta[m]}$$ defined on page \pageref{rhoXdefinition}.
Since $\phi[m]$ is a weak equivalence for all $m$, we conclude from
Theorem \ref{weakequivalenceextension} that the canonical morphism
$$\xymatrix@1{\phi^+ X \co \delta^* N^n(\Delta^{\boxtimes n} / \delta_!X ) \ar[r] &
\delta^*\delta_!X}$$ is a weak equivalence for every simplicial set
$X$.
\end{pf}

\begin{lem} \label{zigzag2}
There is a natural weak equivalence $\xymatrix@1{\delta^*\delta_!X &
X \ar[l]}$.
\end{lem}
\begin{pf}
In Theorem \ref{weakequivalenceextension}, let $F$ be the Yoneda
embedding and $G$ once again $\delta^*\delta_!$. The diagonal
morphism
$$\xymatrix{\Delta[m] \ar[r] & \Delta[m] \times \cdots \times \Delta[m]}$$
is a weak equivalence, as both the source and target are
contractible.
\end{pf}

\begin{prop}
There is a zig-zag of natural weak equivalences between $\delta^*
N^n(\Delta^{\boxtimes n} / \delta_!\text{-} )$ and the identity
functor on $\mathbf{SSet}$. Consequently, there is a natural
isomorphism
$$\xymatrix@1{\delta^* N^n(\Delta^{\boxtimes n} / \delta_!\text{-}
) \ar@{=>}[r] & 1_{\text{\rm Ho}\;\mathbf{SSet}}}.$$
\end{prop}
\begin{pf}
This follows from Proposition \ref{zigzag1} and Lemma \ref{zigzag2}.
\end{pf}


\section{Appendix: The Multidimensional Eilenberg-Zilber Lemma}

In Proposition \ref{multisimpliceshavefibrantconstants} we made use of the multidimensional Eilenberg-Zilber Lemma to
prove that the matching category $\mathcal{B}^i$ is either connected or empty whenever $\mathcal{B}$ is a category of
multisimplices $\Delta^{\times n} Y$. In this Appendix, we prove the multidimensional Eilenberg-Zilber Lemma.
We merely paraphrase Joyal--Tierney's proof of the two-dimensional case in \cite{joyaltierneysimplicial} in order to make the
present paper more self-contained.

\begin{prop}[Eilenberg-Zilber Lemma] \label{prop:EZSSet}
Let $Y$ be simplicial set and $y \in Y_p$. Then there exists a unique surjection $\xymatrix@1{\eta\co [p] \ar[r] & [q]}$
and a unique non-degenerate simplex $y' \in Y_p$ such that $y=\eta^*(y')$.
\end{prop}
\begin{pf}
Proofs can be found in many books on simplicial homotopy theory, for example see Lemma 15.8.4 of \cite{hirschhorn}.
\end{pf}

\begin{defn}
Let $\xymatrix@1{Y\co(\Delta^{\times n})^{\op} \ar[r] & \mathbf{Set}}$ be a
multisimplicial set. A multisimplex $y \in  Y_{\overline{p}}$ is {\it degenerate in direction $i$} if there exists a surjection
$\xymatrix@1{\eta_i \co [p_i] \ar[r] & [q_i]}$  and a multisimplex $y' \in Y_{p_1, \dots, p_{i-1}, q_i, p_{i+1}, \dots, p_n}$
such that $y=(\id_{p_1}, \dots, \id_{p_{i-1}}, \eta,\id_{p_i}, \dots, \id_{p_n})^*(y')$. A multisimplex $y \in  Y_{\overline{p}}$ is {\it non-degenerate in direction $i$} if it is not degenerate in direction $i$. A multisimplex $y \in  Y_{\overline{p}}$ is
{\it totally non-degenerate} if is it not degenerate in any direction.
\end{defn}

\begin{prop}[Multidimensional Eilenberg-Zilber Lemma] \label{prop:EZmultsimplicial}
Let $\xymatrix@1{Y\co(\Delta^{\times n})^{\op} \ar[r] & \mathbf{Set}}$ be a
multisimplicial set and $y \in  Y_{\overline{p}}$. Then there exist unique surjections $\xymatrix@1{\eta_i\co [p_i] \ar[r] & [q_i]}$
and a unique totally non-degenerate multisimplex $y_n \in Y_{\overline{q}}$ such that $y=(\overline{\eta})^* y_n$.
\end{prop}
\begin{pf}
We simply reproduce Joyal--Tierney's proof in Chapter 5 Bisimplicial sets, \cite{joyaltierneysimplicial}.

Let $y=y_0$ for the proof of existence. The Eilenberg-Zilber Lemma for $\mathbf{SSet}$, recalled in Proposition \ref{prop:EZSSet}, guarantees surjections
$\xymatrix@1{\eta_i\co [p_i] \ar[r] & [q_i]}$ and multisimplices $y_i \in Y_{q_1, \dots, q_{i-1}, q_i, p_{i+1}, \dots, p_n}$
such that $$y_{i-1}=(\id_{q_1}, \dots, \id_{q_{i-1}}, \eta_i,\id_{p_{i+1}}, \dots, \id_{p_n})^*(y_i)$$
and each $y_i$ is non-degenerate in direction $i$ for all $i=1,2, \dots, n$. Then $y=(\eta_1, \dots, \eta_n)^*(y_n)$.
The multisimplex $y_n$ is totally non-degenerate,
for if it were degenerate in direction $i$, so that
$$y_n=(\id_{q_1}, \dots, \id_{q_{i-1}}, \eta_i', \id_{q_{i+1}}, \dots, \id_{q_n})^*(y_i'),$$ we would have $y_i$ degenerate in direction $i$:
$$\aligned
y_i &= (\id_{q_1}, \dots, \id_{q_i}, \eta_{i+1}, \dots, \eta_n)^*(y_n) \\
&= (\id_{q_1}, \dots, \id_{q_i}, \eta_{i+1}, \dots, \eta_n)^*(\id_{q_1}, \dots, \id_{q_{i-1}}, \eta_i', \id_{q_{i+1}}, \dots, \id_{q_n})^*(y_i') \\
&=(\id_{q_1}, \dots, \id_{q_{i-1}}, \eta_i', \id_{p_{i+1}}, \dots, \id_{p_n})^*(\id_{q_1}, \dots, \id_{q_i}, \eta_{i+1}, \dots, \eta_n)^*(y_i').
\endaligned$$
But $y_i$ is non-degenerate in direction $i$.

For the uniqueness, suppose $\xymatrix@1{\eta_i'\co [p_i] \ar[r] & [q_i']}$
and $y_n' \in Y_{\overline{q'}}$ is another totally non-degenerate multisimplex such that $y=(\overline{\eta'})^* y_n'$.
The diagram in $\Delta^{\times n}$ associated to the $n$ pushouts in $\Delta$
$$\xymatrix{[p_i] \ar[r]^{\eta_i} \ar[d]_{\eta_i'} & [q_i] \ar[d]^{\mu_i} \\ [q_i'] \ar[r]_{\mu_i'} & [r_i]}$$
is a pushout in $\Delta^{\times n}$ ($\eta_i$ and $\eta_i'$ are all surjective). The Yoneda embedding then gives us a pushout in
$\mathbf{SSet^n}$.
$$\xymatrix@R=4pc@C=4pc{\Delta^{\times n}[\overline{p}] \ar[r]^{\Delta^{\times n}[\overline{\eta}]}
\ar[d]_{\Delta^{\times n}\left[\overline{\eta'}\right]}
& \Delta^{\times n}[\overline{q}] \ar[d]^{\Delta^{\times n}[\overline{\mu}]} \\ \Delta^{\times n}\left[\overline{q'}\right]
\ar[r]_{\Delta^{\times n}\left[\overline{\mu'}\right]} & \Delta^{\times n}[\overline{r}] }$$
Since $$(\overline{\eta'})^* y_n'=y=(\overline{\eta})^* y_n,$$
the universal property of this pushout produces a unique multisimplex $z \in Y_{\overline{r}}$
such that $$y_n'=(\overline{\mu'})^\ast(z),\;\;\;\ y_n=(\overline{\mu})^\ast(z).$$
The multisimplices $y_n$ and $y_n'$ are totally non-degenerate, so $\overline{\mu}=\overline{\id}$ and
$\overline{\mu'}=\overline{\id}$, and consequently $\overline{\eta'}=\overline{\eta}$ and $y_n'=y_n$.
\end{pf}


\begin{thebibliography}{10}

\bibitem{adamekrosicky1994}
Ji{\v{r}}{\'{\i}} Ad{\'a}mek and Ji{\v{r}}{\'\i} Rosick{\'y}.
\newblock {\em Locally presentable and accessible categories}, volume 189 of
  {\em London Mathematical Society Lecture Note Series}.
\newblock Cambridge University Press, Cambridge, 1994.

\bibitem{ehresmannone}
Andr{\'e}e Bastiani and Charles Ehresmann.
\newblock Multiple functors. {I}. {L}imits relative to double categories.
\newblock {\em Cahiers Topologie G\'eom. Diff\'erentielle}, 15(3):215--292,
  1974.

\bibitem{bergercellular}
Clemens Berger.
\newblock A cellular nerve for higher categories.
\newblock {\em Adv. Math.}, 169(1):118--175, 2002.

\bibitem{bergnersurvey}
Julia~E. Bergner.
\newblock A survey of $(\infty,1)$-categories.
\newblock In {\em Proceedings of the {IMA} {W}orkshop `$n$-Categories:
  {F}oundations and {A}pplications' {J}une 2004, {U}niversity of {M}inesota},
  To Appear. http://arxiv.org/abs/math/0610239.

\bibitem{bergnersimplicialcategories}
Julia~E. Bergner.
\newblock A model category structure on the category of simplicial categories.
\newblock {\em Trans. Amer. Math. Soc.}, 359(5):2043--2058 (electronic), 2007.

\bibitem{bergnerthreemodels}
Julia~E. Bergner.
\newblock Three models for the homotopy theory of homotopy theories.
\newblock {\em Topology}, 46(4):397--436, 2007.

\bibitem{brownhigginsgroupoidscrossedcomplexes}
Ronald Brown and Philip~J. Higgins.
\newblock The equivalence of {$\infty $}-groupoids and crossed complexes.
\newblock {\em Cahiers Topologie G\'eom. Diff\'erentielle}, 22(4):371--386,
  1981.

\bibitem{brownhigginsgroupoidscubicalTcomplexes}
Ronald Brown and Philip~J. Higgins.
\newblock The equivalence of {$\omega $}-groupoids and cubical {$T$}-complexes.
\newblock {\em Cahiers Topologie G\'eom. Diff\'erentielle}, 22(4):349--370,
  1981.

\bibitem{brownhigginscubes}
Ronald Brown and Philip~J. Higgins.
\newblock On the algebra of cubes.
\newblock {\em J. Pure Appl. Algebra}, 21(3):233--260, 1981.

\bibitem{brownhigginstensor}
Ronald Brown and Philip~J. Higgins.
\newblock Tensor products and homotopies for {$\omega$}-groupoids and crossed
  complexes.
\newblock {\em J. Pure Appl. Algebra}, 47(1):1--33, 1987.

\bibitem{brownmosa99}
Ronald Brown and Ghafar~H. Mosa.
\newblock Double categories, {$2$}-categories, thin structures and connections.
\newblock {\em Theory Appl. Categ.}, 5:No. 7, 163--175 (electronic), 1999.

\bibitem{cisinskiThomasonFix}
Denis-Charles Cisinski.
\newblock La classe des morphismes de {D}wyer n'est pas stable par retractes.
\newblock {\em Cahiers Topologie G\'eom. Diff\'erentielle Cat\'eg.},
  40(3):227--231, 1999.

\bibitem{cisinskiasterisque}
Denis-Charles Cisinski.
\newblock Les pr\'efaisceaux comme mod\`eles des types d'homotopie.
\newblock {\em Ast\'erisque}, (308):xxiv+390, 2006.

\bibitem{dawsonparepronkpaths}
R.~J.~Mac{G}. Dawson, R.~Par{\'e}, and D.~A. Pronk.
\newblock Paths in double categories.
\newblock {\em Theory Appl. Categ.}, 16:No. 18, 460--521 (electronic), 2006.

\bibitem{dawsonparefreedouble}
Robert Dawson and Robert Par{\'e}.
\newblock What is a free double category like?
\newblock {\em J. Pure Appl. Algebra}, 168(1):19--34, 2002.

\bibitem{duskinII}
John~W. Duskin.
\newblock Simplicial matrices and the nerves of weak {$n$}-categories. {II}.
  {B}icategory morphisms and simplicial maps.
\newblock {\em Preprint}, 2001.

\bibitem{duskinI}
John~W. Duskin.
\newblock Simplicial matrices and the nerves of weak {$n$}-categories. {I}.
  {N}erves of bicategories.
\newblock {\em Theory Appl. Categ.}, 9:198--308 (electronic), 2001/02.
\newblock CT2000 Conference (Como).

\bibitem{ehresmanntwo}
Andr{\'e}e Ehresmann and Charles Ehresmann.
\newblock Multiple functors. {II}. {T}he monoidal closed category of multiple
  categories.
\newblock {\em Cahiers Topologie G\'eom. Diff\'erentielle}, 19(3):295--333,
  1978.

\bibitem{ehresmannthree}
Andr{\'e}e Ehresmann and Charles Ehresmann.
\newblock Multiple functors. {III}. {T}he {C}artesian closed category {${\rm
  Cat}\sb{n}$}.
\newblock {\em Cahiers Topologie G\'eom. Diff\'erentielle}, 19(4):387--443,
  1978.

\bibitem{ehresmannfour}
Andr{\'e}e Ehresmann and Charles Ehresmann.
\newblock Multiple functors. {IV}. {M}onoidal closed structures on {${\rm
  Cat}\sb{n}$}.
\newblock {\em Cahiers Topologie G\'eom. Diff\'erentielle}, 20(1):59--104,
  1979.

\bibitem{ehresmann}
Charles Ehresmann.
\newblock Cat\'egories structur\'ees.
\newblock {\em Ann. Sci. \'Ecole Norm. Sup. (3)}, 80:349--426, 1963.

\bibitem{ehresmann2}
Charles Ehresmann.
\newblock {\em Cat\'egories et structures}.
\newblock Dunod, Paris, 1965.

\bibitem{fiore1}
Thomas~M. Fiore.
\newblock Pseudo limits, biadjoints, and pseudo algebras: categorical
  foundations of conformal field theory.
\newblock {\em Mem. Amer. Math. Soc.}, 182(860), 2006,
  http://arxiv.org/abs/math.CT/0408298.

\bibitem{fiore2}
Thomas~M. Fiore.
\newblock Pseudo algebras and pseudo double categories.
\newblock {\em J. Homotopy Relat. Struct.}, 2(2):119--170, 2007,
  http://arxiv.org/abs/math.CT/0608760.

\bibitem{fiorepaolifundamental}
Thomas~M. Fiore and Simona Paoli.
\newblock The fundamental double category of a bisimplicial set.

\bibitem{fiorepaolipronk1}
Thomas~M. Fiore, Simona Paoli, and Dorette Pronk.
\newblock Model structures on the category of small double categories.
\newblock {\em Algebr. Geom. Topol.}, 8(4):1855--1959, 2008.

\bibitem{fritschlatch1}
Rudolf Fritsch and Dana~May Latch.
\newblock Homotopy inverses for nerve.
\newblock {\em Bull. Amer. Math. Soc. (N.S.)}, 1(1):258--262, 1979.

\bibitem{fritschlatch2}
Rudolf Fritsch and Dana~May Latch.
\newblock Homotopy inverses for nerve.
\newblock {\em Math. Z.}, 177(2):147--179, 1981.

\bibitem{gabrielzisman}
P.~Gabriel and M.~Zisman.
\newblock {\em Calculus of fractions and homotopy theory}.
\newblock Ergebnisse der Mathematik und ihrer Grenzgebiete, Band 35.
  Springer-Verlag New York, Inc., New York, 1967.

\bibitem{gabrielulmer}
Peter Gabriel and Friedrich Ulmer.
\newblock {\em Lokal pr\"asentierbare {K}ategorien}.
\newblock Lecture Notes in Mathematics, Vol. 221. Springer-Verlag, Berlin,
  1971.

\bibitem{goerssjardine}
Paul~G. Goerss and John~F. Jardine.
\newblock {\em Simplicial Homotopy Theory}, volume 174 of {\em Progress in
  Mathematics}.
\newblock Birkh\"auser Verlag, Basel, 1999.

\bibitem{golasinski}
Marek Golasi{\'n}ski.
\newblock Homotopies of small categories.
\newblock {\em Fund. Math.}, 114(3):209--217, 1981.

\bibitem{golasinskiprotranslation}
Marek Golasi{\'n}ski.
\newblock Closed models on the procategory of small categories and simplicial
  schemes.
\newblock {\em Russian Math. Surveys}, 39(5):275--276, 1984.

\bibitem{golasinskipro}
Marek Golasi{\'n}ski.
\newblock Closed models on the procategory of small categories and simplicial
  schemes.
\newblock {\em Uspekhi Mat. Nauk}, 39(5(239)):239--240, 1984.

\bibitem{grandiscospans3}
Marco Grandis.
\newblock Cubical cospans and higher cobordisms ({C}ospans in {A}lgebraic
  {T}opology, {III}).
\newblock {\em J. Homotopy Relat. Struct.}, 3(1):273--308 (electronic), 2007.

\bibitem{grandiscospans1}
Marco Grandis.
\newblock Higher cospans and weak cubical categories (cospans in algebraic
  topology. {I}).
\newblock {\em Theory Appl. Categ.}, 18:No. 12, 321--347 (electronic), 2007.

\bibitem{grandisdouble1}
Marco Grandis and Robert Par\'e.
\newblock Limits in double categories.
\newblock {\em Cahiers Topologie G\'eom. Diff\'erentielle Cat\'eg.},
  40(3):162--220, 1999.

\bibitem{grandisdouble2}
Marco Grandis and Robert Par\'e.
\newblock Adjoints for double categories. {A}ddenda to: ``{L}imits in double
  categories''.
\newblock {\em Cah. Topol. G\'eom. Diff\'er. Cat\'eg.}, 45(3):193--240, 2004.

\bibitem{grandisdouble4}
Marco Grandis and Robert Par\'e.
\newblock Lax {K}an extensions for double categories ({O}n weak double
  categories, {P}art {IV}).
\newblock {\em Cah. Topol. G\'eom. Diff\'er. Cat\'eg.}, 48(3):163--199, 2007.

\bibitem{grandisdouble3}
Marco Grandis and Robert Par\'e.
\newblock Kan extensions in double categories ({O}n weak double categories,
  {P}art {III}).
\newblock {\em Theory Appl. Categ.}, 20:No. 8, 152--185, 2008.

\bibitem{heggiehomotopycofibrations}
Murray Heggie.
\newblock Homotopy cofibrations in {${\rm CAT}$}.
\newblock {\em Cahiers Topologie G\'eom. Diff\'erentielle Cat\'eg.},
  33(4):291--313, 1992.

\bibitem{heggietensorproduct}
Murray Heggie.
\newblock The left derived tensor product of {${\rm CAT}$}-valued diagrams.
\newblock {\em Cahiers Topologie G\'eom. Diff\'erentielle Cat\'eg.},
  33(1):33--53, 1992.

\bibitem{heggiehomotopycolimits}
Murray Heggie.
\newblock Homotopy colimits in presheaf categories.
\newblock {\em Cahiers Topologie G\'eom. Diff\'erentielle Cat\'eg.},
  34(1):13--36, 1993.

\bibitem{hirschhorn}
Philip~S. Hirschhorn.
\newblock {\em Model categories and their localizations}, volume~99 of {\em
  Mathematical Surveys and Monographs}.
\newblock American Mathematical Society, Providence, RI, 2003.

\bibitem{hovey}
Mark Hovey.
\newblock {\em Model categories}, volume~63 of {\em Mathematical Surveys and
  Monographs}.
\newblock American Mathematical Society, Providence, RI, 1999.

\bibitem{hoyo}
Matias Luis~del Hoyo.
\newblock On the subdivision of small categories.
\newblock 2007, http://arxiv.org/abs/0707.1718.

\bibitem{illusieII}
Luc Illusie.
\newblock {\em Complexe cotangent et d\'eformations. {II}}.
\newblock Lecture Notes in Mathematics, Vol. 283. Springer-Verlag, Berlin,
  1972.

\bibitem{jardinesummary}
J.~F. Jardine.
\newblock Categorical homotopy theory.
\newblock {\em Homology, Homotopy Appl.}, 8(1):71--144 (electronic), 2006.

\bibitem{jardineCubical}
J.F. Jardine.
\newblock Cubical homotopy theory: a beginning.
\newblock http://www.math.uwo.ca/$\sim$jardine/papers/preprints/cubical2.pdf.

\bibitem{joyalNotes}
Andr{\'e} Joyal.
\newblock Notes on quasi-categories.
\newblock In {\em Proceedings of the {IMA} {W}orkshop `$n$-Categories:
  {F}oundations and {A}pplications' {J}une 2004, {U}niversity of {M}inesota},
  To Appear.

\bibitem{joyalVolume1}
Andr{\'e} Joyal.
\newblock {\em Theory of Quasi-categories, Volume I}.

\bibitem{joyalVolume2}
Andr{\'e} Joyal.
\newblock {\em Theory of Quasi-categories, Volume II}.

\bibitem{joyaltierneyquasisegal}
Andr{\'e} Joyal and Myles Tierney.
\newblock {Q}uasi-categories vs {S}egal spaces.
\newblock http://arxiv.org/abs/math/0607820.

\bibitem{joyaltierney}
Andr{\'e} Joyal and Myles Tierney.
\newblock Strong stacks and classifying spaces.
\newblock In {\em Category theory (Como, 1990)}, volume 1488 of {\em Lecture
  Notes in Math.}, pages 213--236. Springer, Berlin, 1991.

\bibitem{joyaltierneysimplicial}
Andr{\'e} Joyal and Myles Tierney.
\newblock {\em Notes on Simplicial Homotopy Theory}.
\newblock Preprint, 2008.

\bibitem{kancss}
Daniel~M. Kan.
\newblock On c. s. s. complexes.
\newblock {\em Amer. J. Math.}, 79:449--476, 1957.

\bibitem{kellyenriched}
G.~M. Kelly.
\newblock Basic concepts of enriched category theory.
\newblock {\em Repr. Theory Appl. Categ.}, (10):vi+137 pp. (electronic), 2005.
\newblock Reprint of the 1982 original [Cambridge Univ. Press, Cambridge;
  MR0651714].

\bibitem{kocktrees}
Joachim Kock.
\newblock Polynomial functors and trees.
\newblock http://arxiv.org/abs/0807.2874.

\bibitem{lack2Cat}
Stephen Lack.
\newblock A {Q}uillen model structure for 2-categories.
\newblock {\em $K$-Theory}, 26(2):171--205, 2002.

\bibitem{lackBiCat}
Stephen Lack.
\newblock A {Q}uillen model structure for bicategories.
\newblock {\em $K$-Theory}, 33(3):185--197, 2004.

\bibitem{lackpaoli}
Stephen Lack and Simona Paoli.
\newblock 2-nerves for bicategories.
\newblock {\em $K$-Theory}, 38(2):153--175, 2008.

\bibitem{latchuniqueness}
Dana~May Latch.
\newblock The uniqueness of homology for the category of small categories.
\newblock {\em J. Pure Appl. Algebra}, 9(2):221--237, 1976/77.

\bibitem{lee}
Ming~Jung Lee.
\newblock Homotopy for functors.
\newblock {\em Proc. Amer. Math. Soc.}, 36:571--577; erratum, ibid. 42 (1973),
  648--650, 1972.

\bibitem{leinstersurvey}
Tom Leinster.
\newblock A survey of definitions of {$n$}-category.
\newblock {\em Theory Appl. Categ.}, 10:No. 1,1--70 (electronic), 2002.

\bibitem{leinsternerves}
Tom Leinster.
\newblock Nerves of algebras.
\newblock {\em Talk at CT04, Vancouver}, 2004,
  http://www.maths.gla.ac.uk/~tl/vancouver/.

\bibitem{ncategorycafenerves}
Tom Leinster, Mark Weber, and Others.
\newblock How {I} learned to love the nerve construction.
\newblock {\em The n-Category Caf\'e, A group blog on math, physics and
  philosophy, January 6th, 2008}, http://golem.ph.utexas.edu/category/2008/01/.

\bibitem{lodayfinitelymany}
Jean-Louis Loday.
\newblock Spaces with finitely many nontrivial homotopy groups.
\newblock {\em J. Pure Appl. Algebra}, 24(2):179--202, 1982.

\bibitem{lurieStableInfinity}
Jacob Lurie.
\newblock {\em Derived Algebraic Geometry I: Stable $\infty$-Categories}.
\newblock 2008, http://arxiv.org/abs/math/0608228.

\bibitem{lurieHigherToposTheory}
Jacob Lurie.
\newblock {\em Higher Topos Theory}.
\newblock 2008, http://arxiv.org/abs/math/0608040.

\bibitem{maclaneworking}
Saunders Mac~Lane.
\newblock {\em Categories for the working mathematician}, volume~5 of {\em
  Graduate Texts in Mathematics}.
\newblock Springer-Verlag, New York, second edition, 1998.

\bibitem{maysigurdsson}
J.P. May and J.~Sigurdsson.
\newblock {\em Parametrized homotopy theory}, volume 132 of {\em Mathematical
  Surveys and Monographs}.
\newblock American Mathematical Society, Providence, RI, 2006.

\bibitem{moerdijksvenssonOnshapiro}
I.~Moerdijk and J.-A. Svensson.
\newblock A {S}hapiro lemma for diagrams of spaces with applications to
  equivariant topology.
\newblock {\em Compositio Math.}, 96(3):249--282, 1995.

\bibitem{mortondouble}
Jeffrey Morton.
\newblock A double bicategory of cobordisms with corners.
\newblock http://arxiv.org/abs/math/0611930.

\bibitem{paoliinternalstructures}
Simona Paoli.
\newblock Internal categorical structures in homotopical algebra.
\newblock In {\em Proceedings of the {IMA} {W}orkshop `$n$-Categories:
  {F}oundations and {A}pplications' {J}une 2004, {U}niversity of {M}inesota},
  To Appear. http://www.maths.mq.edu.au/$\sim$simonap/.

\bibitem{pellissier}
Regis Pellissier.
\newblock Weak enriched categories - {C}ategories enrichies faibles.
\newblock 2003, http://arxiv.org/abs/math/0308246.

\bibitem{quillenI}
Daniel Quillen.
\newblock Higher algebraic {$K$}-theory. {I}.
\newblock In {\em Algebraic $K$-theory, I: Higher $K$-theories (Proc. Conf.,
  Battelle Memorial Inst., Seattle, Wash., 1972)}, pages 85--147. Lecture Notes
  in Math., Vol. 341. Springer, Berlin, 1973.

\bibitem{rezkcat}
Charles Rezk.
\newblock A model category for categories.
\newblock 2000, http://www.math.uiuc.edu/$\sim$rezk/cat-ho.dvi.

\bibitem{rezkhomotopytheory}
Charles Rezk.
\newblock A model for the homotopy theory of homotopy theory.
\newblock {\em Trans. Amer. Math. Soc.}, 353(3):973--1007 (electronic), 2001.

\bibitem{shulmanonquillenfunctors}
Michael Shulman.
\newblock Comparing composites of left and right derived functors.
\newblock http://arxiv.org/abs/0706.2868.

\bibitem{shulmanframed}
Michael Shulman.
\newblock Framed bicategories and monoidal fibrations.
\newblock http://arxiv.org/abs/0706.1286.

\bibitem{simpsonmodel}
Carlos Simpson.
\newblock A closed model structure for $n$-categories, internal ${H}om$,
  $n$-stacks and generalized {S}eifert-{V}an {K}ampen.
\newblock 1997, http://arxiv.org/abs/alg-geom/9704006.

\bibitem{simpsonhigher}
Carlos Simpson.
\newblock {H}omotopy theory of higher categories.
\newblock 2010, http://arxiv.org/abs/1001.4071.

\bibitem{street}
Ross Street.
\newblock The algebra of oriented simplexes.
\newblock {\em J. Pure Appl. Algebra}, 49(3):283--335, 1987.

\bibitem{tamsamani}
Zouhair Tamsamani.
\newblock Sur des notions de {$n$}-cat\'egorie et {$n$}-groupo\"\i de non
  strictes via des ensembles multi-simpliciaux.
\newblock {\em $K$-Theory}, 16(1):51--99, 1999.

\bibitem{thomasonhocolimit}
R.~W. Thomason.
\newblock Homotopy colimits in the category of small categories.
\newblock {\em Math. Proc. Cambridge Philos. Soc.}, 85(1):91--109, 1979.

\bibitem{thomasonCat}
R.W. Thomason.
\newblock Cat as a closed model category.
\newblock {\em Cahiers Topologie G\'eom. Diff\'erentielle}, 21(3):305--324,
  1980.

\bibitem{toenaxiomatization}
Bertrand To{\"e}n.
\newblock Vers une axiomatisation de la th\'eorie des cat\'egories
  sup\'erieures.
\newblock {\em $K$-Theory}, 34(3):233--263, 2005.

\bibitem{waldhausen}
Friedhelm Waldhausen.
\newblock Algebraic {$K$}-theory of spaces.
\newblock In {\em Algebraic and geometric topology (New Brunswick, N.J.,
  1983)}, volume 1126 of {\em Lecture Notes in Math.}, pages 318--419.
  Springer, Berlin, 1985.

\bibitem{weber}
Mark Weber.
\newblock Familial 2-functors and parametric right adjoints.
\newblock {\em Theory Appl. Categ.}, 18:No. 22, 665--732, 2007.

\bibitem{worytkiewicz2Cat}
K.~Worytkiewicz, K.~Hess, P.~E. Parent, and A.~Tonks.
\newblock A model structure \`a la {T}homason on {\bf 2-{c}at}.
\newblock {\em J. Pure Appl. Algebra}, 208(1):205--236, 2007.

\end{thebibliography}
\end{document}